\documentclass[draft]{article}
\def\today{8.1.11} 
\usepackage{amsmath,amsfonts,amsthm,amssymb,amscd}

\renewcommand{\Im}{\mathop{\rm Im}\nolimits}
\def\S{\mathhexbox278}

\newcommand{\beq}{\begin{equation}}
\newcommand{\ee}{\end{equation}}

\theoremstyle{plain} \newtheorem{theorem}{Theorem}[section]
\newtheorem{lemma}[theorem]{Lemma}
\newtheorem{proposition}[theorem]{Proposition}
 \theoremstyle{definition}
\newtheorem{definition}[theorem]{Definition} \theoremstyle{remark}
\newtheorem{remark}[theorem]{Remark}

\newcommand{\R}{{\mathbb R}} \newcommand{\U}{{\mathcal U}}

\newcommand{\Z}{{\mathbb Z}}

\newcommand{\Tr}{{\mathcal T}}
\newcommand{\Tra}{{\phi}}

\newcommand{\Ph}{{\mathcal P}}

\newcommand{\resto}{{\mathcal R}}
\def\im{{\rm i}}

\newcommand{\C}{\mathbb{C}}

\def\uno{{\kern+.3em {\rm 1} \kern -.22em {\rm l}}}

\def\norma#1{\left\| #1\right\|}

\numberwithin{equation}{section}

\begin{document}

\title{The Hamiltonian structure of the nonlinear
Schr\"odinger equation and   the  asymptotic stability of its ground
states}

\author {Scipio Cuccagna}

\date{\today}
\maketitle
\begin{abstract}
In this paper we prove that  ground states of the NLS which satisfy
the sufficient conditions  for orbital   stability of M.Weinstein,
are also asymptotically stable, for seemingly generic equations.
Here we assume that the NLS has a smooth short range nonlinearity.
We assume also the presence of a very short range and smooth linear
potential, to avoid translation invariance. The basic idea is to
perform a Birkhoff normal form argument on the hamiltonian, as in a
paper by Bambusi and Cuccagna on the stability of the 0 solution for
NLKG. But in our case, the natural coordinates arising from the
linearization  are not canonical. So we need also to apply the
Darboux Theorem. With some care though,  in order not to destroy
some nice features of the initial hamiltonian.
\end{abstract}

\section{Introduction}
\label{section:introduction} We consider the   nonlinear
Schr\"odinger equation (NLS)

\begin{equation}\label{NLS}
 \im u_{t }=-\Delta u +Vu +\beta  (|u|^2) u , \, u(0,x)=u_0(x), \, (t,x)\in\mathbb{ R}\times
 \mathbb{ R}^3
\end{equation}
with   $ -\Delta +V(x)$ a selfadjoint Schr\"odinger operator. Here
  $V(x)\neq 0$  to exclude   translation invariance. We assume
  that both $V(x)$ and $\beta (|u|^2)u$ are  short range
  and smooth.
  We assume that \eqref{NLS} has a smooth family of ground states.
  We then prove that the sufficient conditions for
 orbital stability by Weinstein \cite{W1} (which, essentially, represent the
 correct definition of linear stability, see \cite{Cu3}),    imply for a generic \eqref{NLS} that
  the ground states are not only orbitally stable, as proved in
\cite{W1} (under less restrictive hypotheses), but that their
orbits are also asymptotically stable. That is, a solution $u(t)$
of \eqref{NLS} starting sufficiently close to ground states, is
asymptotically of the form $e^{i\theta (t)} \phi _{\omega _+} (x)+
e^{it\Delta }h_+$, for $\omega _+$ a fixed number and for $h_+\in
H^1(\mathbb{R}^3)$ a small energy function. The problem of
stability of ground states has a long history.  Orbital stability
has been well understood since the 80's, see in the sequence
\cite{CL,W1,GSS1,GSS2}, and has been a very active field
afterwards. Asymptotic stability is a more recent, and less
explored, field. In the context of the NLS the first results are in
the pioneering works \cite{SW1,SW2,BP1,BP2}. Almost all references
on asymptotic stability of ground states of the NLS tackle the
problem by first linearizing at ground states, and by attempting to
deal with the resulting nonlinear problem for the error term. An
apparent problem in the linear theory is that the linearization is
a not symmetric operator. However the linearization is covered by
the scattering theory of non selfadjoint operators developed     by
T.Kato in the 60's, see his classical \cite{kato}, see also
\cite{CPV,Schlag}. Dispersive and Strichartz estimates for the
linearization, analogous to the theory for short range scalar
Schr\"odinger operators elaborated in \cite{JSS,Y1,Y2}, to name
only few of many papers, can be proved using   similar ideas, see
for example \cite{Cu1,Schlag,KS}. It is fair to say that anything
that can be proved for short range scalar Schr\"odinger operators,
can also be proved   for the linearizations. The  only notable
exception is the problem of "positive signature" embedded
eigenvalues, see \cite{Cu3}, which we conjecture not to exist (in
analogy to the absence of embedded eigenvalues for short range
Schr\"odinger operators), and which in any case are   unstable, see
\cite{CPV}. Hence it is reasonable to focus on  NLS's where these
positive signature embedded eigenvalues do not exist (in the case
of ground states, all positive eigenvalues are of positive
signature).While linear theory is not a problem in understanding
asymptotic stability, the real trouble lies in the difficult NLS
like equation one obtains for the error term. Specifically, the
linearization has discrete spectrum which, at the level of  linear
theory, tends not to decay and potentially could yield
quasiperiodic solutions.   A good analogy with more standard
problems, is that the continuous spectrum of the linearization
corresponds to stable spectrum while the discrete spectrum
corresponds to central directions. Stability cannot be established
by linear theory alone. The first intuition on how nonlinear
interactions are responsible for loss of energy of the discrete
modes, is in a paper by Sigal \cite{sigal}.
 His ideas, inspired by the classical Fermi golden rule in linear
 theory,
are later elaborated in \cite{SW3}, to study asymptotic stability of
vacuum  for the nonlinear Klein Gordon equations with a potential
with non empty discrete spectrum. This problem, easier than the one
treated in the present paper, to a large extent is solved in
\cite{bambusicuccagna}. In reality, the main ideas in  \cite{SW3}
 had already be sketched, for the problem of stability of ground
 states of NLS, in a deep
 paper by Buslaev and Perelman  \cite{BP2}, see also
the expanded version \cite{BS}. In the case when the linearization
has just one positive eigenvalue close to the continue spectrum,
\cite{SW3,BP2}, or \cite{sigal} in a different context, identify the
mechanisms for loss of energy of the discrete modes in the nonlinear
coupling of continuous and discrete spectral components.
Specifically,  in the discrete mode equation there is a key
coefficient   of the form $\langle D F,F\rangle $ for $D $ a
  positive operator and $F$ a  function. Assuming the
generic condition $\langle D F,F\rangle\not=0 $, this   gives rise
to dissipative effects leading to leaking of energy from the
discrete mode to the continuous modes, where energy disperses
because of linear dispersion, and to the ground state. After
\cite{BP2} there is strong evidence  that, generically, linearly
stable ground states, in the sense of \cite{W1}, should be
asymptotically stable. Still, it is a seemingly    technically
difficult problem to solve rigorously. After \cite{BP2,SW3}, a
number of papers analyze  the same ideas in various situations,
\cite{TY1,TY2,TY3,T,SW4,Cu2}. In the meantime, a useful series of
papers \cite{GNT,M1,M2} shows how to use endpoint Strichartz and
smoothing estimates  to prove     in energy space the result of
\cite{SW2,PiW}, generalizing the result and simplifying the
argument. The next important breakthrough is due to Zhou and Sigal
\cite{zhousigal}. They tackle for the first time the case of one
positive eigenvalue arbitrarily close to 0, developing further the
normal forms analysis of \cite{BP2} and obtaining the rate of
leaking conjectured in     \cite{SW3} p.69. The argument is improved
in \cite{cuccagnamizumachi}. The crucial coefficient is now of the
form $\langle D F,G\rangle $, with $F$ and $G$ not obviously
related. In \cite{cuccagnamizumachi} it is noticed that $\langle D
F,G\rangle <0$ is incompatible with orbital stability  (an argument
along these lines is suggested in     \cite{SW3} p.69). So, for
orbitally stable ground states, the generic condition $\langle D
F,G\rangle \neq 0$ implies positivity, and hence leaking of energy
out of the discrete
  modes. This yields a result similar to \cite{sigal,BP2,SW3}
  and in particular is a partially positive answer to a conjecture
  on p.69 in \cite{SW3}.
The case with more than one positive eigenvalue   is harder.  In
this case, due to possible cancelations, \cite{cuccagnamizumachi} is
not able to draw
  conclusions on the sign of the coefficients
under the assumption of orbital stability. But, apart from the issue
of positivity of the coefficients, \cite{cuccagnamizumachi} shows
that the rest of the proof does not depend on the number of positive
eigenvalues. Moreover, \cite{T,zhouweinstein1,Cu3}   show that if
there are many positive eigenvalues, all close to the continuous
spectrum, then the important coefficients are again of the form
$\langle D F,F\rangle $.  The reason for this lies in the
hamiltonian nature of the NLS. The above papers contain normal forms
arguments. The hamiltonian structure is somewhat lost in the above
papers. When the eigenvalues are close to the continuous spectrum,
the normal form argument consists of just one step. This single step
does not change the crucial coefficients. Then, the hamiltonian
nature of the initial system, yields information on these
coefficients (this is emphasized in \cite{Cu3}).  In the case
treated in \cite{zhousigal,cuccagnamizumachi}   though, there are
many steps in the normal form. The important coefficients are
changed in ways which look very complicated,  see \cite{Gz}  which
deals with the next two easiest cases after the easiest. The correct
way to look at this problem is introduced in \cite{bambusicuccagna},
which deals with the   problem introduced in \cite{SW3}. Basically,
the positivity can be seen by doing the normal form directly on the
hamiltonian. We give a preliminary and heuristic justification  on
why the hamiltonian structure is crucial at the end of    section
\ref{section:linearization}. \cite{bambusicuccagna} consists in a
mixture of  a Birkhoff normal forms argument, with the arguments in
\cite{cuccagnamizumachi}. For asymptotic stability of ground states
of NLS though, \cite{bambusicuccagna} is still not enough. Indeed in
\cite{bambusicuccagna} something peculiar   happens: the natural
coordinates arising by the  spectral decomposition of the
linearization  at the  vacuum solution, are also canonical
coordinates for the symplectic structure. This is no  longer true if
instead of vacuum we consider    ground states. So we need an extra
step, which consists in the search of canonical coordinates, through
the Darboux theorem. This step requires care, because we must make
sure that our problem remains similar to a semilinear NLS also in
the new system of coordinates.

In a forthcoming paper, Zhou and Weinstein \cite{zhouweinstein2}
track precisely in the setting of \cite{zhouweinstein1}  how much of
the energy of the discrete modes goes to the ground state and how
much is dispersed. For another result on asymptotic stability, that
is asymptotic stability of the blow up profile, we refer to
\cite{MR}. In some respects the situation in \cite{MR} is harder
than here, since there the additional discrete modes are
concentrated in the kernel of the linearization. There is important
work on asymptotic stability for KdV equations due to Martel and
Merle, see  \cite{MM1} and further references therein, which solve a
problem initiated by Pego and Weinstein \cite{PW}, the latter closer
in spirit to our approach to NLS.  It is an interesting question  to
see if elaboration of ideas in \cite{MM1,MMT} can be used for
alternative solutions of the problem which we consider here. Our
result does not cover important cases, like the pure power NLS, with
$\beta (|u|^2)= -|u|^{p-1}$ and $V=0$, where our result is probably
false. Indeed it  is well known that in 3D ground states are stable
for $p<7/3$ and unstable for $p\ge 7/3$. In the $p<7/3$ case there
are ground states of arbitrarily small $H^1 $ norm. They are
counterexamples to the asymptotic stability in $H^1 $ of the 0
solution. Then for $p>5/3$ the 0 solution  is asymptotically stabile
in a smaller space usually denoted by $\Sigma$, which involves also
the $\| x u\| _{L^2_x}$ norm, see in \cite{strauss} the comments
after Theorem 6 p. 55. In $\Sigma $ there are no small ground states
for $p\in (5/3,7/3 )$. Presumably one should be able to prove
asymptotic stability of ground states in $\Sigma$. To our knowledge
even the following (presumably easier) problem is not solved yet:
the asymptotic stability of 0 in $\Sigma$  when $V\neq 0$, $\sigma
_p (-\Delta +V )=\emptyset $ and $\beta (|u|^2)= -|u|^{p-1}$ with
$p\in (5/3,7/3 )$. In  the literature on asymptotic stability of
ground states like \cite{BP2,BS,zhousigal,cuccagnamizumachi}, the
case of moving solitons is left aside, because  in that set up it
appears substantially more complex. We do not treat moving solitons
here either, but it is possible that our approach might help also
with moving solitons. In the step when we perform the Darboux
Theorem, the velocity should freeze and we should reduce to the same
situation considered from section \ref{section:reformulation} on.
The extra difficulty with moving solitons is that there are more
obstructions to the fact that after Darboux we have a semilinear
NLS. But it would be  surprising if this difficulty had a really
deep nature.   In any case, the main conceptual problem stemming
from \cite{sigal,BP2,SW3}, which we solve here, is the issue of the
positive semidefiniteness of the critical coefficients. There is a
growing literature on interaction between solitons, see for example
\cite{MM2,HW,M3}, and we expect our result to be relevant.

 We do not reference all the literature on asymptotic stability of
ground states, see \cite{cuccagnatarulli} for more.  We like to
conclude observing that Sigal \cite{sigal}, Buslaev and Perelman
 \cite{BP2}  and Soffer and Weinstein \cite{SW3} had identified
 with great precision  the right mechanism of leaking
 of energy away from the discrete modes.

\section{Statement of the main result}
\label{section:statement}

We will assume the following hypotheses.

\begin{itemize}
\item[(H1)] $\beta  (0)=0$, $\beta\in C^\infty(\R,\R)$.
\item[(H2)] There exists a $p\in(1,5)$ such that for every
$k\ge 0$ there is a fixed $C_k$ with
$$\left| \frac{d^k}{dv^k}\beta(v^2)\right|\le C_k
|v|^{p-k-1} \quad\text{if $|v|\ge 1$}.$$
\item[(H3)] $V(x)$ is smooth and for any multi
index $\alpha $ there are $C_\alpha >0$ and $a_\alpha >0$ such that
 $|\partial ^\alpha _x V(x)|\le
C_\alpha e^{-a_\alpha |x|}$.
\item[(H4)]
There exists an open interval $\mathcal{O}$ such that
\begin{equation}
  \label{eq:B}
  \Delta u-Vu-\omega u+\beta(|u|^2)u=0\quad\text{for $x\in \R^3$},
\end{equation}
admits a $C^1$-family of ground states $\phi _ {\omega }(x)$ for
$\omega\in\mathcal{O}$.
\item [(H5)]
\begin{equation}
  \label{eq:1.2}
\frac d {d\omega } \| \phi _ {\omega }\|^2_{L^2(\R^3)}>0
\quad\text{for $\omega\in\mathcal{O}$.}
\end{equation}
\item [(H6)]
Let $L_+=-\Delta +V +\omega -\beta (\phi _\omega ^2 )-2\beta '(\phi
_\omega ^2) \phi_\omega^2$ be the operator whose domain is $H^2
 (\R^3)$. Then $L_+$ has exactly one negative eigenvalue and
does not have kernel.

\item [(H7)] Let $\mathcal{H}_\omega$ be the linearized operator around $e^{it\omega}\phi_\omega$
(see Section \ref{section:linearization} for the precise
definition). There is a fixed $m\ge 0$ such that
$\mathcal{H}_\omega$ has  $m$ positive eigenvalues $\lambda
_1(\omega )\le \lambda _2(\omega ) \le ...\le \lambda _m(\omega )$.
We assume there are fixed integers $m_0=0< m_1<...<m_{l_0}=m$ such that
$\lambda _j(\omega )= \lambda _i(\omega )$ exactly for $i$ and $j$
both in $(m_l, m_{l+1}]$ for some $l\le l_0$. In this case $\dim
\ker (\mathcal{H}_\omega -\lambda _j(\omega ) )=m_{l+1}-m_l$.
 We assume there exist $N_j\in \mathbb{N}$ such that
$0<N_j\lambda _j(\omega )< \omega < (N_j+1)\lambda _j(\omega )$ with
$N_j\ge 1$. We set $N=N_1$.
\item [(H8)]  There is no multi index $\mu \in \mathbb{Z}^{m}$
with $|\mu|:=|\mu_1|+...+|\mu_m|\leq 2N_1+3$ such that $\mu \cdot
\lambda =\omega$.

\item[(H9)] If $\lambda _{j_1}<...<\lambda _{j_k}$ are $k$ distinct
  $\lambda$'s, and $\mu\in \Z^k$ satisfies
  $|\mu| \leq 2N_1+3$, then we have
$$
\mu _1\lambda _{j_1}+\dots +\mu _k\lambda _{j_k}=0 \iff \mu=0\ .
$$
\item[(H10)] $\mathcal{H}_\omega$ has no other eigenvalues except for $0$ and
the $ \pm \lambda _j (\omega )$. The points $\pm \omega$ are not
resonances.

\item [(H11)]
The Fermi golden rule  Hypothesis (H11)   in subsection
\ref{subsec:FGR}, see \eqref{eq:FGR}, holds.

\end{itemize}

\begin{remark}
\label{rem:Prelim FGR}    The   novelty of this paper with
respect to \cite{cuccagnamizumachi} is that we prove that  some
crucial coefficients are of a specific form,  see \eqref{eq:FGR}. As
a consequence, see Lemma \ref{lemma:FGR8}, these coefficients are
positive semidefinite.    In the analogue of \eqref{eq:FGR}  in
\cite{cuccagnamizumachi}, see Hypothesis 5.2 p.72
\cite{cuccagnamizumachi}, except for the special case $n=1$ of just one
 eigenvalue (or of possibly many eigenvalues but all with $N_j=1$), there is no clue on the sign of the term
on the rhs of the key inequality, and the fact that it is positive
is an hypothesis.
\end{remark}

\begin{theorem}\label{theorem-1.1}
  Let $\omega_1\in\mathcal{O}$ and $\phi_{\omega_1}(x)$
be a ground state of \eqref{NLS}. Let $u(t,x)$ be a solution to
\eqref{NLS}. Assume (H1)--(H10). Then, there exist an $\epsilon_0>0$
and a $C>0$ such that if
$\epsilon:=\inf_{\gamma\in[0,2\pi]}\|u_0-e^{\im \gamma}\phi_
{\omega _1} \|_{H^1}  <\epsilon_0,$ there exist $\omega
_\pm\in\mathcal{O}$, $\theta\in C^1(\R;\R)$ and $h _\pm \in H^1$
with $\| h_\pm\| _{H^1}+|\omega _\pm -\omega_1|\le C \epsilon $
such that

\begin{equation}\label{scattering}
\lim_{t\to  \pm\infty}\|u(t,\cdot)-e^{\im \theta(t)}\phi_{\omega
_\pm}-e^{\im t\Delta }h _\pm\|_{H^1}=0 .
\end{equation}
It is possible to write $u(t,x)=e^{\im \theta(t)}\phi_{\omega  (t)}
+ A(t,x)+\widetilde{u}(t,x)$  with $|A(t,x)|\le C_N(t) \langle x
\rangle ^{-N}$ for any $N$, with $\lim _{|t|\to \infty }C_N(t)=0$,
with $\lim _{ t \to \pm \infty } \omega  (t)= \omega _\pm$, and such
that for any pair $(r,p)$ which is admissible, by which we mean that
\begin{equation}\label{admissiblepair}  2/r+3/p= 3/2\,
 , \quad 6\ge p\ge 2\, , \quad
r\ge 2,
\end{equation}
we have
\begin{equation}\label{Strichartz} \|  \widetilde{u} \|
_{L^r_t( \mathbb{R},W^{1,p}_x)}\le
 C\epsilon .
\end{equation}

\end{theorem}

We end the introduction with some notation. Given two functions
$f,g:\mathbb{R}^3\to \mathbb{C}$ we set $\langle f,g\rangle = \int
_{\mathbb{R}^3}f(x)  g(x)  dx$. Given a matrix $A$, we denote by
$A^*$, or by $^tA$,  its transpose. Given two vectors $A$ and $B$,
we denote by  $A^*B=\sum _j A_jB_j$ their inner product. Sometimes
we omit the summation symbol, and we use the convention on sum over
repeated indexes. Given two functions $f,g:\mathbb{R}^3\to \mathbb{C
}^2$ we set $\langle f,g\rangle = \int _{\mathbb{R}^3}f^*(x) g(x)
dx$. For any $k,s\in \mathbb{R}$ and any Banach space $K$ with field
$\C$
 \[ H^{ k,s}(\mathbb{R}^3,K)=\{ f:\mathbb{R}^3\to K \text{ s.t.}
 \| f\| _{H^{s,k}}:=\| \langle x \rangle ^s \| (-\Delta +1)^{k} f
 \| _{K}\| _{L^2
 }<\infty \},\]
   $(-\Delta +1)^{k} f(x)= (2\pi )^{-\frac{3}{2}}
  \int e^{ \im x \cdot \xi }(\xi ^2+1)^k \widehat{f}(\xi ) d\xi $,
$\widehat{f}(\xi )= (2\pi )^{-\frac{3}{2}}\int e^{ -\im x \cdot \xi
}  {f}(x ) dx$. In particular we set
 $L^{2,s} =H^{0,s}  $, $L^2=L^{2,0}    $,  $H^k=H^{2,0}    $.
Sometimes, to emphasize that these spaces refer to spatial
variables, we will denote them by $W^{k,p}_x$, $L^{ p}_x$, $H^k_x$,
$H^{ k,s}_x$ and $L^{2,s}_x$. For $I$ an interval and $Y_x$ any of
these spaces, we will consider Banach spaces $L^p_t( I, Y_x)$ with
mixed norm $ \| f\| _{L^p_t( I, Y_x)}:= \| \| f\| _{Y_x} \| _{L^p_t(
I )}.$ Given an operator $A$, we will denote by $R_A(z)=(A-z)^{-1}$
its resolvent. We set $\mathbb{N}_0=\mathbb{N}\cup \{0 \}$. We will
consider multi indexes $\mu \in \mathbb{N}_0^n$. For   $\mu \in
\mathbb{Z}^n$ with $\mu =(\mu _1,..., \mu _n)$ we set $|\mu |=\sum
_{j=1}^n |\mu _j|.$     For $X$ and $Y$ two Banach space, we will
denote by $B(X,Y)$ the Banach space of bounded linear operators from
$X$ to $Y$ and by $B^{\ell}(X,Y)= B (  \prod _{j=1}^\ell X ,Y)$. We
denote by $a^{\otimes \ell}$ the element $\otimes _{j=1}^\ell a$ of
$\otimes _{j=1}^\ell X$ for some $a\in X$. Given a differential form
$\alpha$, we denote by
    $d\alpha$ its exterior differential.

{\bf Acknowledgments} I wish to thank Dario Bambusi for pointing out
a gap in the proof of an earlier version of Theorem \ref{th:main}.

\section{Linearization and set up}
\label{section:linearization}

 Let $U={^t(u,\overline{u})}$. We introduce now energy $E(u)$ and mass $Q(u)$. We set

 \begin{equation} \label{eq:energyfunctional}\begin{aligned}&
 E(U)=E_K(U)+E_P(U)\\&
E_K(U)= \int _{\R ^3}
  \nabla u \cdot \nabla \overline{u} dx+ \int _{\R ^3}
  V u   \overline{u} dx \\&
E_P(U)=
 \int _{\R ^3}B( u \overline{u}) dx \end{aligned}
\end{equation}
with $B(0)=0$ and $\partial _{\overline{u}}B(|u|^2)=\beta (|u|^2)u$.
We will consider the matrices \begin{equation}
\label{eq:Pauli}\begin{aligned} &\sigma _1=
\begin{pmatrix}0 &
1  \\
1 & 0
 \end{pmatrix} \, ,
\sigma _2=\begin{pmatrix}  0 &
\im  \\
-\im & 0
 \end{pmatrix} \, ,
\sigma _3=\begin{pmatrix} 1 & 0\\0 & -1 \end{pmatrix} .
\end{aligned}
\end{equation}
We introduce the mass  \begin{equation}\label{eq:charge}Q(U)=  \int _{\R ^3}u
\overline{u} dx= \frac{1}{2}\langle U, \sigma _1 U\rangle .
\end{equation}
 Let
\begin{equation} \label{eq:function q} \Phi _\omega =\begin{pmatrix} \phi _\omega
   \\ \phi _\omega
 \end{pmatrix}   , \, q(\omega )=Q(\Phi _\omega ),
\, e (\omega )=E(\Phi _\omega ), \, d(\omega )=e (\omega )+\omega
q(\omega ) .
\end{equation}
Often we will denote $\Phi _\omega $ simply by $\Phi$.
The
\eqref{NLS} can be written as
\begin{equation}\label{eq:NLSvectorial} \im \dot U =
\begin{pmatrix} 0 &1
   \\ -1 & 0
 \end{pmatrix}   \begin{pmatrix} \partial _{u}E
   \\ \partial _{\overline{u}}E
 \end{pmatrix}  = \sigma _3 \sigma _1 \nabla E (U),
\end{equation}
with $\nabla E (U)$ defined by \eqref{eq:NLSvectorial}.
We have for $\vartheta \in \R$
\begin{equation}\label{eq:gaugeInvariance}
 E( e^{-\im \sigma _3\vartheta }  U)=E(  U)  \text{ and }
 \nabla E ( e^{-\im \sigma _3\vartheta }
  U)=e^{ \im \sigma _3\vartheta } \nabla E (  U)  .
\end{equation}
Write  for $\omega \in \mathcal{O}$
\begin{equation}
U=  e^{\im \sigma _3\vartheta }  (\Phi _\omega + R).\nonumber
\end{equation}
Then
\begin{equation}\label{system}\im \dot  U= -\sigma _3
\dot \vartheta e^{\im \sigma _3\vartheta }
 (\Phi _\omega + R) + \im \dot \omega e^{\im \sigma _3\vartheta }
  \partial _\omega \Phi _\omega +  \im    e^{\im \sigma _3\vartheta }
  \dot R
\end{equation}
and
\begin{equation}
\begin{aligned} &
 -\sigma _3
\dot \vartheta e^{\im \sigma _3\vartheta }
 (\Phi _\omega + R) + \im \dot \omega e^{\im \sigma _3\vartheta }
  \partial _\omega \Phi _\omega +  \im    e^{\im \sigma _3\vartheta }
  \dot R   =
  \sigma _3 \sigma _1 e^{-\im \sigma _3 \vartheta }\nabla E (\Phi _\omega
  +R).\end{aligned}  \nonumber
\end{equation}
Equivalently we get

\begin{equation}\label{system1}
\begin{aligned} &
  -\sigma _3
(\dot \vartheta -\omega )
 (\Phi _\omega + R) + \im \dot \omega
  \partial _\omega \Phi _\omega +  \im
  \dot R  =\\& =  \sigma _3 \sigma _1   \left ( \nabla E (\Phi _\omega
  +R)+\omega \nabla Q (\Phi _\omega
  +R)\right ) .\end{aligned}
\end{equation}
  We have $\frac{d}{dt}\sigma _3 \sigma _1   \left ( \nabla E (\Phi _\omega
  +tR)+\omega \nabla Q (\Phi _\omega
  +tR)\right ) _{\mid _{t=0}}= \mathcal{H}_\omega  R$ with

\begin{equation}  \label{eq:linearization} \begin{aligned} &
\mathcal{H}_\omega   = \sigma_3(-\Delta+V+\omega)
  +\sigma_3
\left[\beta (\phi ^2_{\omega }) +\beta ^\prime (\phi ^2_{\omega
})\phi ^2_{\omega } \right] -\im  \sigma _2 \beta ^\prime (\phi ^2
_{\omega })\phi ^2 _{\omega } .\end{aligned}
\end{equation}
The essential spectrum of $\mathcal{H}_\omega$ consists of $(-\infty
, -\omega ]\cup [ \omega,+\infty )$. It is well known (see
\cite{W2}) that by  (H5)  $0$ is an isolated eigenvalue of
$\mathcal{H}_\omega$ with $\dim N_g(\mathcal{H}_\omega)=2$ and
\begin{equation}\label{eq:Kernel} \mathcal{H}_\omega\sigma_3\Phi_\omega=0,\quad
\mathcal{H}_\omega\partial_\omega\Phi_\omega =-\Phi_\omega.
\end{equation}
Since $\mathcal{H}_\omega^*=\sigma_3\mathcal{H}_\omega\sigma_3$, we
have $N_g(\mathcal{H}_\omega^*)=\operatorname{span}\{\Phi_\omega,
\sigma_3\partial_\omega\Phi_\omega\}$. We consider
  eigenfunctions $\xi _j(\omega)$   with eigenvalue   $\lambda
_j(\omega)$:
$$
\mathcal{H}_\omega\xi _j(\omega)=\lambda _j(\omega)\xi
_j(\omega),\quad \mathcal{H}_\omega\sigma_1\xi _j(\omega)=-\lambda
_j(\omega)\sigma_1\xi _j(\omega) .$$ They can be normalized so that
$\langle \sigma_3  \xi _j(\omega),\overline{\xi} _\ell(\omega)
\rangle =\delta _{j\ell }$, this is based on Proposition 2.4
\cite{Cu3}. Furthermore, they can be chosen to be real, that is with
real entries, so $\xi _j=\overline{ \xi} _j$ for all $j$.

Both $\phi_\omega$ and $\xi _j(\omega,x)$ are smooth in
$\omega\in\mathcal{O}$ and $x\in\R^3$ and satisfy
$$\sup_{\omega\in\mathcal{K},x\in\R^3} e^{a|x|}( |\partial ^\alpha _x\phi_\omega(x)|+
\sum _{j=1}^{m}|\partial ^\alpha _x \xi _j(\omega,x)| <\infty$$ for
every $a\in(0,\inf_{\omega\in\mathcal{K}}\sqrt{\omega-\lambda
_m(\omega)})$ and every compact subset $\mathcal{K}$ of
$\mathcal{O}$.

For $\omega\in\mathcal{O}$, we have the
$\mathcal{H}_\omega$-invariant Jordan block decomposition
\begin{align}  \label{eq:spectraldecomp} &
L^2(\R^3,\C^2)=N_g(\mathcal{H}_\omega)\oplus \big (\oplus _{\pm}
\oplus _{j=1}^m \ker (\mathcal{H}_\omega\mp \lambda _j(\omega))
\big)\oplus L_c^2(\mathcal{H}_\omega),
\end{align}
  $L_c^2(\mathcal{H}_\omega):=
\left\{N_g(\mathcal{H}_\omega^\ast)\oplus \big (\oplus _{\lambda \in
\sigma _d\backslash \{ 0\}}   \ker (\mathcal{H}_\omega ^*- \lambda
 ) \big)\right\} ^\perp $ with $\sigma _d =\sigma _d
(\mathcal{H}_\omega)$. We also set $L_d^2(\mathcal{H}_\omega):=
N_g(\mathcal{H}_\omega)\oplus \big (\oplus _{\lambda \in \sigma
_d\backslash \{ 0\}}   \ker (\mathcal{H}_\omega  - \lambda (\omega))
\big ) .$ By $P_c(\mathcal{H}_{\omega})$ (resp.
$P_d(\mathcal{H}_{\omega})$), or simply by $P_c( {\omega})$ (resp.
$P_d( {\omega})$), we denote the projection on
$L_c^2(\mathcal{H}_\omega)$  (resp. $L_d^2(\mathcal{H}_\omega)$)
associated to the above direct sum. The space
$L^2_c(\mathcal{H}_{\omega})$ depends continuously on $\omega$.
  We specify the
ansatz imposing that
\begin{equation}\label{eq:anzatz}
 U=  e^{\im \sigma _3\vartheta }  (\Phi _\omega + R)
 \text{ with $\omega \in \mathcal{O}$, $\vartheta \in \R$ and $R\in N^{\perp}_g (\mathcal{H}_\omega ^*)$.}
\end{equation}
We consider  coordinates
\begin{equation}\label{eq:coordinate}U= e^{\im \sigma _3\vartheta}
 \left( \Phi
_\omega +z \cdot \xi  (\omega )+  \overline{z }\cdot \sigma_1\xi
 (\omega )+P_c(\mathcal{H}_{\omega}) f\right )    \end{equation}
where $\omega \in \mathcal{O}$, $z\in \mathbb{ C} $ and $f\in
L^2_c(\mathcal{H}_{\omega _0})$ where we fixed  $\omega _0\in
\mathcal{O}$  such that $q(\omega _0)=\| u_0\| _2^2$.
\eqref{eq:coordinate} is a system  of coordinates  if we use the
notation $ \mathcal{O}$ to denote a small neighborhood of $\omega
_1$ in Theorem \ref{theorem-1.1}.  Indeed by Lemma \ref{lem:regular}
below, then the map $P_c(\mathcal{H}_{\omega})$ is an isomorphism
from $L^2_c(\mathcal{H}_{\omega _0})$ to $L^2_c(\mathcal{H}_{\omega
})$. In particular
\begin{align}
  \label{eq:decomp2}
& R  =\sum _{j=1}^{m}z_j \xi _j(\omega )+
\sum _{j=1}^{m}\overline{z }_j\sigma_1\xi _j(\omega )
+P_c(\mathcal{H}_{\omega  } )f ,\\
\label{eq:decomp3} & R \in  N_g^\perp(\mathcal{H}_{\omega
}^*)\quad\text{and}\quad f \in L_c^2(\mathcal{H}_{\omega _0}).
\end{align}
We also set $z\cdot \xi =\sum _jz_j \xi _j$ and $\overline{z}\cdot
\sigma_1\xi =\sum _j\overline{z}_j \sigma_1\xi _j$. In the sequel we
set
\begin{equation} \label{eq:partialR} \partial _\omega R=
\sum _{j=1}^{m}z_j \partial _\omega \xi _j(\omega )+ \sum
_{j=1}^{m}\overline{z }_j\sigma_1\partial _\omega \xi _j(\omega
)+\partial _\omega P_c(\mathcal{H}_{\omega } )f.
\end{equation}
We have:
\begin{lemma}
  \label{lem:regular} We have $P_c(\mathcal{H}_{\omega } )^*=
  P_c(\mathcal{H}_{\omega } ^*)$ for all $ \omega
  \in \mathcal{O}$.
 For all $ \omega , \widetilde{\omega}\in \mathcal{O}$ the following operators
  are bounded from $H^{-k,-s}$ to $H^{ k', s'}$ for all exponents:
  \begin{equation} \label{eq:list op1} \begin{aligned} &\partial ^\ell _\omega
  P_c(\mathcal{H}_{\omega } ) \text{ for any $\ell >0$}\, ; \\&
    P_c(\mathcal{H}_{\omega } )
  -P_c(\mathcal{H}_{\omega } ^*)  \, ; \,  P_c(\mathcal{H}_{\omega } )
  -P_c(\mathcal{H}_{\widetilde{\omega}}  ).
  \end{aligned}
  \end{equation}
Consider $\omega _1$ of Theorem
  \ref{theorem-1.1}. There exists $\varepsilon _1 >0$
  such that   $(\omega _1-\varepsilon _1, \omega _1+\varepsilon _1)
  \subset \mathcal{O}$, and for any pair $\widetilde{\omega}   , \omega
  \in (\omega _1-\varepsilon _1, \omega _1+\varepsilon _1)$
we have \begin{equation} \label{isomorphism}
  \text{$P_c( \omega    ) P_c( \widetilde{\omega}
   ):   L^2_c(\mathcal{H}_{\widetilde{\omega}   }  )\to
L^2_c(\mathcal{H}_{\omega   }  ) $ is an isomorphism}\end{equation}
Furthermore,
  the following operator
  is bounded from $H^{-k,-s}$ to $H^{ k', s'}$ for all exponents:
  \begin{equation} \label{eq:list op2} \begin{aligned} &
   P_c(\mathcal{H}_{\widetilde{\omega}  }  )
\left (  1- (P_c(\mathcal{H}_{\omega  }  )
P_c(\mathcal{H}_{\widetilde{\omega}
 }  ))^{-1} \right )  P_c(\mathcal{H}_{\omega  }  )
  \end{aligned}
  \end{equation}
  where in the last line
and $(P_c( \omega    ) P_c( \widetilde{\omega}     ))^{-1}$ is the
inverse of the operator in \eqref{isomorphism}. Finally, for
$\epsilon _0 $ in Theorem \ref{theorem-1.1} sufficiently small, we
have $|\omega _0 -\omega _1|<\varepsilon _1$, with $\omega _0$
defined under \eqref{eq:coordinate}.
\end{lemma} \proof The first statement follows from the definition.
We have $P_c(\mathcal{H}_{\omega})=1-P_d(\mathcal{H}_{\omega})$
where $P_d(\mathcal{H}_{\omega})$ are finite linear combinations of
rank 1 operators $\Psi \langle \Psi ', \quad \rangle $ with
  $\Psi, \Psi ' \in H^{K,S}$ for any $(K,S)$. This implies the statement for
the second line of \eqref{eq:list op1}.  $\partial ^\ell _\omega
  P_c(\mathcal{H}_{\omega } ) $ is well defined by  the fact that in
(H4) the dependence on $\omega$ is in fact smooth (this is seen
iterating the argument in Theorem 18 \cite{shatahstrauss}). Assuming
\eqref{isomorphism}, and for  $P_c=P_c( \omega )$,  $\widetilde{P}_c
=P_c( \widetilde{\omega } )$, $P_d=P_d( \omega )$,  $\widetilde{P}_d
=P_d( \widetilde{\omega } )$, we have \begin{equation}
\begin{aligned} & \widetilde{P}_c \left (  1- (P_c \widetilde{P}_c
)^{-1} \right ) P_c = (\widetilde{P}_c
-P_c)P_c-(\widetilde{P}_c-P_c)(P_c \widetilde{P}_c)^{-1}P_c ,
\end{aligned} \nonumber
  \end{equation}
which yields   \eqref{eq:list op2}.
 We  prove \eqref{isomorphism}.
First of all the map is 1--1. Indeed if $P_c  \widetilde{P}_cf=0$,
then we have $f= P_df= (P_d-\widetilde{P}_d)f$. Then $\| f\| _{2}\le
C |\omega -\widetilde{\omega }| \| f\| _{2}$ for some fixed $C>0$.
This, for $2C \varepsilon _1<1$, is compatible only with $f=0$. If
we knew that the map in \eqref{isomorphism} is onto, then
\eqref{isomorphism} would  hold    by the open mapping theorem. So suppose
the map is not onto. Let $\mathcal{R}(P_c \widetilde{P}_c )$ be the
range of $P_c \widetilde{P}_c $. If there exists $\widetilde{g}\in
L^2_c(\mathcal{H}^*_{\omega} )$ such that $\widetilde{g}\neq 0$ and
 $\langle
\widetilde{g} ,P_c \widetilde{P}_c f\rangle =0$ for all $f \in
L^2_c(\mathcal{H}_{\widetilde{\omega } } )$, then since
$\widetilde{g}=\sigma _3 g$ for a $g \in L^2_c(\mathcal{H}_{\omega }
)$, we get $0= \langle \widetilde{g} ,P_c \widetilde{P}_c f\rangle =
\langle \widetilde{P}_c  P_c g ,\sigma _3 f\rangle$ for all $f \in
L^2_c(\mathcal{H}_{\widetilde{\omega }})$. This implies
$\widetilde{P}_c P_c g=0$, and since $g \in
L^2_c(\mathcal{H}_{\omega } )$, by the 1--1 argument this implies
$g=0$. So if the map in \eqref{isomorphism} is not onto, then
$\mathcal{R}(P_c \widetilde{P}_c )$ is dense in $
L^2_c(\mathcal{H}_{\omega } )$. We will see in a moment that
$\mathcal{R}(P_c \widetilde{P}_c )$ is closed in $
L^2_c(\mathcal{H}_{\omega } )$,   hence   concluding that the map
in \eqref{isomorphism} is also onto. To see that $\mathcal{R}(P_c
\widetilde{P}_c )$ is closed in $ L^2_c(\mathcal{H}_{\omega } )$,
let $\widetilde{f}_n\in L^2_c(\mathcal{H}_{\widetilde{\omega}  } )$
be a sequence such that $\| P_c\widetilde{f}_n -f\| _2\to 0 $ for
some $f\in L^2_c(\mathcal{H}_{\omega } )$. By $\| \widetilde{f}_n \|
_2\le \|  P_c \widetilde{f}_n \| _2 +C|\omega -\widetilde{\omega} |
\| \widetilde{f}_n \| _2$ for some fixed $C$, it follows that for
$2C\varepsilon _{1}<1$ the sequence $\| \widetilde{f}_n \| _2$ is
bounded. Then by weak compactness there is a subsequence
$\widetilde{f} _{n_j}$ weakly convergent to a $\widetilde{f}\in
L^2_c(\mathcal{H}_{\widetilde{\omega} } )$. Since $P_c
\widetilde{P}_c$ is also weakly continuous, $P_c
\widetilde{P}_c\widetilde{f}=f$.

  \noindent The final statement is elementary by \eqref{eq:1.2}.\qed

Using the system of coordinates \eqref{eq:coordinate} we rewrite
system \eqref{system1} as
\begin{equation} \label{eq:system}
\begin{aligned} &
  -\sigma _3
 (\dot \vartheta -\omega )
 (\Phi _\omega +   z\cdot \xi +
     {\overline{z}}\cdot \sigma _1\xi + P_c(\mathcal{ H}_\omega ) f)+  \\& + \im \dot \omega
 ( \partial _\omega \Phi _\omega +z\cdot \partial _\omega \xi +
    \overline{z} \cdot \sigma _1\partial _\omega \xi +\partial _\omega
P_c(\mathcal{H}_\omega ) f) +\im
  \dot z\cdot \xi +\\&  +    \im
   \dot {\overline{z}}\cdot \sigma _1\xi +\im
    P_c(\mathcal{ H}_\omega ) \dot f     = \sigma _3
    \sigma _1
  \nabla E  (  \Phi _\omega +  z\cdot \xi +
     {\overline{z}}\cdot \sigma _1\xi +
     P_c(\mathcal{ H}_\omega )  f) \\& +\omega \sigma _3
    \sigma _1 \nabla Q  (  \Phi _\omega +  z\cdot \xi +
     {\overline{z}}\cdot \sigma _1\xi +
     P_c(\mathcal{ H}_\omega )  f)    ,\end{aligned}
\end{equation}
 where $z\cdot \xi =\sum _j z_j\xi _j$ and
$\overline{z}\cdot \sigma _1\xi =\sum _j \overline{z}_j\sigma _1\xi _j$, where $\xi =\xi (\omega )$. Notice for future reference, that fixed any $\omega _0$
we also have
\begin{equation} \label{system2}
\begin{aligned} &
  -\sigma _3
 (\dot \vartheta -\omega _{0})
 (\Phi _\omega +   z\cdot \xi +
     {\overline{z}}\cdot \sigma _1\xi + P_c(\mathcal{ H}_\omega ) f)+  \\& + \im \dot \omega
 ( \partial _\omega \Phi _\omega +z\cdot \partial _\omega \xi +
    \overline{z} \cdot \sigma _1\partial _\omega \xi +\partial _\omega
P_c(\mathcal{H}_\omega ) f) +\im
  \dot z\cdot \xi +\\&  +    \im
   \dot {\overline{z}}\cdot \sigma _1\xi +\im
    P_c(\mathcal{ H}_\omega ) \dot f     = \sigma _3
    \sigma _1
  \nabla E  (  \Phi _\omega +  z\cdot \xi +
     {\overline{z}}\cdot \sigma _1\xi +
     P_c(\mathcal{ H}_\omega )  f) \\& +\omega _{0}\sigma _3
    \sigma _1 \nabla Q  (  \Phi _\omega +  z\cdot \xi +
     {\overline{z}}\cdot \sigma _1\xi +
     P_c(\mathcal{ H}_\omega )  f)    ,\end{aligned}
\end{equation}
where \eqref{system2} is the same of \eqref{eq:system} except for $\omega _0$
replacing $\omega $ in the first spot where they appear in the first and last line.

We end this section  with a short heuristic description about why
the crucial property needed to prove asymptotic stability of ground
states, is the hamiltonian nature of the \eqref{NLS}. In terms of
\eqref{eq:coordinate}, and oversimplifying, \eqref{system} splits as
\begin{equation}  \begin{aligned} &
\im \dot z -\lambda z = \sum _{\mu \nu}a _{\mu \nu}
z^{\mu}\overline{z}^{\nu}+ \sum _{\mu \nu} z^{\mu}\overline{z}^{\nu}
\langle f(t,x) ,\overline{ G} _{\mu \nu}(x,\omega )   \rangle _{L^2_x}+\cdots
\\& \im \dot f -\mathcal{H}_\omega f = \sum _{\mu \nu}
z^{\mu}\overline{z}^{\nu}    M _{\mu \nu}(x,\omega )+ \cdots
.\end{aligned} \nonumber
  \end{equation}
Here we are assuming $m=1$. We focus on positive times $t\ge 0$
only. After changes of variables, see \cite{cuccagnamizumachi}, we
obtain
\begin{equation} \label{final} \begin{aligned} &
\im \dot z -\lambda z = P(|z|^2)z+  \overline{z}^N \langle f(t,x) ,\overline{ G} _{\mu \nu}(x,\omega )   \rangle _ {L^2_x}+\cdots
\\& \im \dot f -\mathcal{H}_\omega f = z^{N+1}
     M  (x,\omega )+ \cdots
.\end{aligned}
  \end{equation}
The next step is to write, for $g$ an error term,

\begin{equation}  \begin{aligned} &  f= -z^{N+1}R_{\mathcal{H}_\omega}^{+}((N +1)\, \lambda )M +  g\\&
\\& \im \dot z-\lambda z= P(|z|^2)z- |z|^{2N }z\langle R_{\mathcal{H}_\omega}^{+}((N +1)\,
\lambda )M , \overline{G}   \rangle _{L^2_x} +...\end{aligned}\nonumber
  \end{equation}
   Then, ignoring error terms, by $$R_{\mathcal{H}_\omega}^{+}
   ((N +1)\, \lambda
) =  P.V. \dfrac{1}{{\mathcal{H}_\omega}-(N +1)\,\lambda}+\im \pi
\delta ({\mathcal{H}_\omega}-(N +1)\, \lambda )$$ the equation for
$z$ has solutions such that
  $$  \frac{d}{dt}|z|^2=-\Gamma |z|^{2N+2 }  ,
 |z(t)|=\dfrac{|z(0)|}{ (|z(0)|^{2N} \, N\, \Gamma \, t+
1)^{\frac{1}{2N}}}$$ with   (the Fourier transforms are associated
to $\mathcal{H}_\omega$)
  $$\Gamma =
2\pi \langle \delta (\mathcal{H}_\omega-(N +1)\,\lambda )M, G
\rangle = \int _{ |\xi| = \sqrt{(N +1)\, \lambda -\omega}  }\frac{\widehat{M}(\xi
) \cdot \overline{\widehat{G} (\xi )}d\sigma}{\sqrt{(N +1)\, \lambda -\omega}} .$$ If  $\Gamma
>0$, we see that $z(t)$ decays. Notice that $\Gamma
<0$ is incompatible with orbital stability, which requires $z$ to
remain small, see Corollary 4.6 \cite{cuccagnamizumachi}. The latter
indirect argument to prove positive semidefiniteness of $\Gamma$,
does not seem to work  when in \eqref{system} there are further
discrete components. So we need another way to prove that $\Gamma
\ge 0$. This is  provided by the hamiltonian structure. Indeed, if
\eqref{final} is of the form

\begin{equation} \label{ham} \begin{aligned} &
{\rm i} \dot z = \partial   _{\overline{z}} K \, , \quad {\rm i}  \dot  f = \nabla
 _{\overline{f}} K,
\end{aligned}
  \end{equation}
  then   by Schwartz lemma $(N+1)!M=\partial ^{N+1}_z\nabla _{\overline{f}}K
  =\overline{\partial ^{N }_{\overline{z}}\nabla _{f}\partial _{\overline{z}}K}=N! {G}$ at $z=0$ and $f=0$.
  So $\Gamma$ is positive semidefinite. This very simple
  idea   on system \eqref{ham}, inspired \cite{bambusicuccagna}
   and inspires the present
  paper.

\section{Gradient of the coordinates}
\label{section:modulation} We focus   on  ansatz \eqref{eq:anzatz}
and on the coordinates \eqref{eq:coordinate}. In particular we
compute  the gradient of the coordinates. Here we recall that given
a scalar valued function $F$, the relation between exterior
differential and gradient is $dF=\langle \nabla F, \quad \rangle $.
Consider the following two functions
\begin{equation} \mathcal{F}(U,\omega , \vartheta ):=\langle
e^{-\im \sigma _3\vartheta }U-\Phi _\omega , \Phi _\omega\rangle
\text{ and } \mathcal{G}(U,\omega , \vartheta ):=\langle e^{-\im
\sigma _3\vartheta }U  ,\sigma _3\partial _\omega  \Phi
_\omega\rangle . \nonumber
\end{equation}
Then ansatz \eqref{eq:anzatz} is obtained by choosing $(\omega ,
\vartheta )$ s.t. $R:=e^{-\im \sigma _3\vartheta }U-\Phi _\omega$
satisfies $R\in N_g^\perp (\mathcal{ H} _\omega ^*)$ by means of the
implicit function theorem. In particular:
\begin{equation}
\begin{aligned} & \mathcal{F}_{\vartheta} =-\im \langle \sigma _3e^{-\im \sigma _3\vartheta
}U, \Phi _\omega  \rangle =-\im \langle \sigma _3R, \Phi _\omega
\rangle ;\\& \mathcal{F}_{\omega} = -2 q'(\omega)+\langle e^{-\im
\sigma _3\vartheta }U, \partial _\omega \Phi _\omega  \rangle = -
q'(\omega)+\langle  R, \partial _\omega \Phi _\omega  \rangle
; \\& \nabla _U \mathcal{F}=e^{-\im \sigma _3\vartheta
} \Phi _\omega  \, , \, \nabla _U \mathcal{G}=e^{-\im \sigma
_3\vartheta } \sigma _3\partial _\omega \Phi _\omega ;\\&
\mathcal{G}_{\vartheta} =-\im \langle
 e^{-\im \sigma _3\vartheta }U,\partial _\omega  \Phi _\omega  \rangle
=-\im  ( q'(\omega)+   \langle  R, \partial _\omega \Phi _\omega
\rangle   ) ;\\& \mathcal{G}_{\omega} =\langle e^{-\im \sigma
_3\vartheta }U ,\sigma _3\partial _\omega ^2  \Phi _\omega\rangle
=\langle R ,\sigma _3\partial _\omega ^2  \Phi _\omega\rangle .
\end{aligned}\nonumber
\end{equation}
By $\mathcal{F}(U,\omega (U), \vartheta (U)=\mathcal{G}(U,\omega (U), \vartheta (U)=0$ we get $ \mathcal{W}_{\omega} \nabla \omega +\mathcal{W}_{\vartheta} \nabla \vartheta =-\nabla _UW$ for $\mathcal{W}=\mathcal{F},\mathcal{G}.$
By the above formulas,  if we set
\begin{equation}\label{eq:matrixA} \mathcal{A}=\begin{pmatrix} - q'(\omega)+\langle  R, \partial _\omega \Phi _\omega  \rangle &
-\im \langle \sigma _3R, \Phi _\omega \rangle \\ \langle R ,\sigma
_3\partial _\omega ^2  \Phi _\omega\rangle  & -\im  ( q'(\omega)+
\langle  R, \partial _\omega \Phi _\omega \rangle   ) \end{pmatrix}
\end{equation}
we have
\begin{equation}\label{eq:ApplmatrixA}
\mathcal{A}
\begin{pmatrix} \nabla \omega  \\ \nabla \vartheta \end{pmatrix}
=\begin{pmatrix} -e^{-\im \sigma _3\vartheta } \Phi _\omega  \\
-e^{-\im \sigma _3\vartheta } \sigma _3\partial _\omega\Phi _\omega
\end{pmatrix} .
\end{equation}
So

\begin{equation}\label{eq:GradModulation}
\begin{aligned} &
\nabla \omega =\frac{(q'(\omega)+ \langle  R, \partial _\omega \Phi
_\omega \rangle   )e^{-\im \sigma _3\vartheta } \Phi _\omega
-\langle \sigma _3R, \Phi _\omega \rangle e^{-\im \sigma _3\vartheta
} \sigma _3\partial _\omega\Phi _\omega}{(q'(\omega))^2 -\langle  R,
\partial _\omega \Phi _\omega \rangle ^2 + \langle \sigma _3R, \Phi
_\omega \rangle \langle R ,\sigma _3\partial _\omega ^2  \Phi
_\omega\rangle }
\\ & \nabla \vartheta =\frac{\langle R ,\sigma _3\partial _\omega ^2  \Phi
_\omega\rangle e^{-\im \sigma _3\vartheta } \Phi _\omega
+(q'(\omega)- \langle  R, \partial _\omega \Phi _\omega \rangle
)e^{-\im \sigma _3\vartheta }  \sigma _3\partial _\omega\Phi
_\omega} {\im \left [ q'(\omega))^2 -\langle R,
\partial _\omega \Phi _\omega \rangle ^2 + \langle \sigma _3R, \Phi
_\omega \rangle \langle R ,\sigma _3\partial _\omega ^2  \Phi
_\omega\rangle\right ]} \, .
\end{aligned}
\end{equation}
Notice that along with the decomposition \eqref{eq:spectraldecomp}
we have \begin{align}  \label{eq:dualspectraldecomp} &
L^2(\R^3,\C^2)=N_g(\mathcal{H}_\omega ^*)\oplus \big (\oplus
_{\lambda \in \sigma _d\backslash \{ 0\}}   \ker (\mathcal{H}_\omega
^*- \lambda  (\omega)) \big)\oplus L_c^2(\mathcal{H}_\omega ^*),
\end{align}
  $L_c^2(\mathcal{H}_\omega ^*):= \left\{N_g(\mathcal{H}_\omega
)\oplus\big (\oplus _{\lambda \in \sigma _d\backslash \{ 0\}}   \ker
(\mathcal{H}_\omega  - \lambda  (\omega)) \big) \right\} ^\perp.$ We
also set $L_d^2(\mathcal{H}_\omega ^*):= N_g(\mathcal{H}_\omega
^*)\oplus \big (\oplus _{\lambda \in \sigma _d\backslash \{ 0\}}
\ker (\mathcal{H}_\omega ^*  - \lambda (\omega)) \big ) .$ Notice
that $N_g(\mathcal{H}_\omega ^*) =\sigma _3N_g(\mathcal{H}_\omega
)$, $\ker (\mathcal{H}_\omega ^*- \lambda ) =\sigma _3\ker
(\mathcal{H}_\omega - \lambda ) $,  $L_c^2(\mathcal{H}_\omega
^*)=\sigma _3L_c^2(\mathcal{H}_\omega  )$ and
$L_d^2(\mathcal{H}_\omega ^*)=\sigma _3L_d^2(\mathcal{H}_\omega  )$,
  so that \eqref{eq:dualspectraldecomp} is obtained applying $\sigma
_3 $ to decomposition \eqref{eq:spectraldecomp}. We can decompose
gradients as
\begin{equation}\label{eq:gradient}\begin{aligned} &
 \nabla F(U)=e^{-\im \sigma _3\vartheta } \big [ P_{N_g(\mathcal{H}^*_\omega)}
 + \\& \sum _j  (P _{
 \ker (\mathcal{H}^*_\omega-\lambda _j)}
   +
 P _{
 \ker (\mathcal{H}^*_\omega +\lambda _j)}  )
  +
  P _c (\mathcal{H}_\omega ^*)
   \big ] e^{ \im \sigma _3\vartheta }\nabla F (U)=\\&
   \frac{\langle \nabla F(U) ,
    e^{ \im \sigma _3\vartheta }\partial _\omega \Phi \rangle }
  {q'(\omega )}
  e^{- \im \sigma _3\vartheta }  \Phi +
  \frac{\langle \nabla F (U), e^{ \im \sigma _3\vartheta }
  \sigma _3 \Phi \rangle }
  {q'(\omega )}
  e^{- \im \sigma _3\vartheta }\sigma _3\partial _\omega \Phi
  \\&  +
   \sum _j\langle \nabla F(U) ,
   e^{ \im \sigma _3\vartheta }  \xi _j \rangle e^{ -\im \sigma _3\vartheta }
   \sigma _3\xi _j  +\sum _j    \langle \nabla F(U) ,
   e^{ \im \sigma _3\vartheta }\sigma _1 \xi _j \rangle
  e^{ -\im \sigma _3\vartheta }\sigma _1\sigma _3\xi _j \\& +  e^{  -\im \sigma _3\vartheta }
   P_c(\mathcal{H}_\omega ^* )e^{  \im \sigma _3\vartheta }\nabla F(U).
 \end{aligned}
\end{equation}

Using coordinates \eqref{eq:coordinate} and
notation\eqref{eq:partialR}, at $U$ we have the following formulas
for the vectorfields
\begin{equation}\label{eq:vectorfields} \begin{aligned} &
\frac \partial {\partial  {\omega}}   =
   e^{ \im \sigma _3\vartheta } \partial _\omega ( \Phi +R)
\, ,\, \frac \partial {\partial  {\vartheta}} =\im
   e^{ \im \sigma _3\vartheta } \sigma _3 ( \Phi +R)
,\\& \frac \partial {\partial  {z_j}}  =
   e^{ \im \sigma _3\vartheta }  \xi _j   \, ,\,
   \frac \partial {\partial  {\overline{z}_j}}   =
   e^{ \im \sigma _3\vartheta }\sigma _1 \xi _j  .\end{aligned}
\end{equation}
Hence, by  $\partial _{\omega}F =dF(\frac \partial {\partial
{\omega}})=\langle \nabla F, \frac \partial {\partial  {\omega}}
\rangle $ etc.,  we have
\begin{equation}\label{eq:derivativeZ} \begin{aligned} &
\partial _{\omega}F =\langle \nabla F  ,
   e^{ \im \sigma _3\vartheta } \partial _\omega ( \Phi +R) \rangle
\, ,\, \partial _{\vartheta }F  =\im \langle \nabla F  ,
   e^{ \im \sigma _3\vartheta } \sigma _3 ( \Phi +R) \rangle
,\\&
\partial _{z_j}F =
\langle \nabla F  ,
   e^{ \im \sigma _3\vartheta }  \xi _j \rangle \, ,\,
   \partial _{\overline{z}_j}F = \langle \nabla F  ,
   e^{ \im \sigma _3\vartheta }\sigma _1 \xi _j\rangle .\end{aligned}
\end{equation}

  \begin{lemma} \label{lem:gradient z} We have the following
formulas:

\begin{eqnarray} \label{ZOmegaTheta}&
\nabla z_j =  - \langle \sigma _3 \xi _j, \partial _\omega R \rangle
  \nabla \omega - \im \langle \sigma _3 \xi _j, \sigma _3 R \rangle
  \nabla \vartheta +  e^{-\im \sigma _3\vartheta }
  \sigma
_3\xi _j\\  \label{barZOmegaTheta}&\nabla \overline{z}_j =  -
\langle \sigma _1\sigma _3 \xi _j,
\partial _\omega R \rangle \nabla \omega-
   \im \langle \sigma _1\sigma _3 \xi _j, \sigma _3 R \rangle
   \nabla \vartheta  +  e^{-\im \sigma _3\vartheta }
  \sigma _1\sigma
_3\xi _j.\end{eqnarray}
\end{lemma}
\proof   Equalities   $ \frac{\partial z_j}{\partial z_\ell }  =\delta _{j\ell}$,
$ \frac{\partial z_j}{\partial \overline{z}_\ell }  = \frac{\partial z_j}{\partial \omega }= \frac{\partial z_j}{\partial \vartheta } =0$
 and $\nabla _f z_j=0$  are equivalent to
\begin{equation}\label{eq:indentitiesGradZ} \begin{aligned} & \langle \nabla z_j,
e^{ \im \sigma _3\vartheta } \xi _\ell \rangle =\delta _{j\ell},
\langle \nabla z_j, e^{ \im \sigma _3\vartheta } \sigma _1\xi _\ell
\rangle \equiv 0 =\langle \nabla z_j, e^{ \im \sigma _3\vartheta }
\sigma _3(\Phi +R) \rangle \\& \langle \nabla z_j, e^{ \im \sigma
_3\vartheta } \partial _\omega (\Phi +R) \rangle =0\equiv \langle
\nabla z_j, e^{ \im \sigma _3\vartheta } P _c( \omega ) P _c( \omega
_0)g \rangle \, \forall g\in L^2_c(\mathcal{H}_{\omega _0}).
\end{aligned}
\end{equation}
Notice that the last identity implies $ P _c( \mathcal{H}_{\omega
_0}^{*} ) P _c( \mathcal{H}_{\omega  }^{*} )e^{ \im \sigma
_3\vartheta }\nabla z_j=0$ which in turn implies $   P _c(
\mathcal{H}_{\omega  }^{*} )e^{ \im \sigma _3\vartheta }\nabla
z_j=0$. Then , applying \eqref{eq:gradient} and using the product
row column, we get for some pair of numbers $(a,b)$
\begin{equation} \begin{aligned} & \nabla z_j=a
e^{- \im \sigma _3\vartheta }  \Phi
+b e^{- \im \sigma _3\vartheta }\sigma _3\partial _\omega \Phi +
e^{-\im \sigma _3\vartheta } \sigma _3\xi _j
\\& =  (a,b)  \begin{pmatrix} e^{- \im \sigma _3\vartheta }  \Phi
\\ e^{- \im \sigma _3\vartheta }\sigma _3\partial _\omega \Phi \end{pmatrix} +   e^{-\im \sigma _3\vartheta }
  \sigma
_3\xi _j  =-(a,b)  \mathcal{A} \begin{pmatrix} \nabla \omega
\\ \nabla \vartheta \end{pmatrix} +   e^{-\im \sigma _3\vartheta }
  \sigma
_3\xi _j,\end{aligned}\nonumber
\end{equation}
where in the last line we used \eqref{eq:ApplmatrixA}. Equating the two extreme sides and  applying to the formula $\langle \quad , \frac{\partial}{\partial \omega}\rangle $ and $\langle \quad , \frac{\partial}{\partial \vartheta}\rangle $, by $\langle \nabla z_j , \frac{\partial}{\partial \omega}\rangle =\langle \nabla z_j , \frac{\partial}{\partial \vartheta}\rangle =\langle \nabla \vartheta , \frac{\partial}{\partial \omega}\rangle =\langle \nabla \omega , \frac{\partial}{\partial \vartheta}\rangle =0 $, by
$ \langle \nabla \vartheta , \frac{\partial}{\partial \vartheta}\rangle =\langle \nabla \omega , \frac{\partial}{\partial \omega}\rangle =1 $ and by
  \eqref{eq:vectorfields} and \eqref{eq:indentitiesGradZ},  we get
\begin{equation}  \mathcal{A}^* \begin{pmatrix} a  \\ b \end{pmatrix} =
\begin{pmatrix} \langle \sigma _3 \xi _j, \partial _\omega R \rangle
  \\ \im \langle \sigma _3 \xi _j, \sigma _3 R \rangle \end{pmatrix}.\nonumber
\end{equation}
This implies
\begin{equation} \begin{aligned} &
\nabla z_j =  -(\langle \sigma _3 \xi _j, \partial _\omega R \rangle
  , \im \langle \sigma _3 \xi _j, \sigma _3 R \rangle )\begin{pmatrix} \nabla \omega
\\ \nabla \vartheta \end{pmatrix}+  e^{-\im \sigma _3\vartheta }
  \sigma
_3\xi _j.\end{aligned}\nonumber
\end{equation}
This yields \eqref{ZOmegaTheta}. Similarly
\begin{equation} \nabla \overline{z}_j=a e^{- \im \sigma _3\vartheta }  \Phi
+b e^{- \im \sigma _3\vartheta }\sigma _3\partial _\omega \Phi +
e^{-\im \sigma _3\vartheta } \sigma _1\sigma _3\xi _j , \nonumber
\end{equation}
where
\begin{equation}  \mathcal{A}^* \begin{pmatrix} a  \\ b \end{pmatrix} =
\begin{pmatrix} \langle \sigma _1\sigma _3 \xi _j, \partial _\omega R \rangle
  \\ \im \langle \sigma _1\sigma _3 \xi _j, \sigma _3 R \rangle \end{pmatrix}.\nonumber
\end{equation}
 \qed

\begin{lemma} \label{lem:gradient f} Consider the map
$f(U)=f$ for $U$ and $f$ as in \eqref{eq:coordinate}. Denote by
$f'(U) $ the Frech\'et  derivative of this map. Then

\begin{equation} \begin{aligned} &
f'(U)= (P_c(\omega )P_c(\omega _0))^{-1} P_c(\omega )\left [
 -   \partial _\omega R \, d\omega -\im
  \sigma _3 R  \, d\vartheta +  e^{-\im
\sigma _3\vartheta }\uno \right ] .
\end{aligned} \nonumber
\end{equation}
               \end{lemma}
\proof  We have
\begin{equation}\label{eq:indentitiesGradf} \begin{aligned} &
f'(U) e^{ \im \sigma _3\vartheta } \xi _\ell \equiv f'(U) e^{ \im
\sigma _3\vartheta } \sigma _1\xi _\ell
 \equiv 0= f'(U) e^{ \im \sigma _3\vartheta }
\sigma _3(\Phi +R) =  \\&  f'(U) e^{ \im \sigma _3\vartheta }
\partial _\omega (\Phi +R) \rangle   \text { and } f'(U)
e^{ \im \sigma _3\vartheta } P_c(\omega  )g=g \, \forall g\in
L^2_c(\mathcal{H}_{\omega _0}).
\end{aligned}
\end{equation}
This implies that for a pair of vectors valued functions $A$ and $B$
and with the inverse of $P_c(\mathcal{H}_{\omega }
)P_c(\mathcal{H}_{\omega _0}):L^2_c(\mathcal{H}_{\omega _0})\to
L^2_c(\mathcal{H}_{\omega  })$,
\begin{equation} \begin{aligned} & f'=   (A,B)
\begin{pmatrix} \langle e^{- \im \sigma _3\vartheta }  \Phi , \quad
\rangle
\\ \langle e^{- \im \sigma _3\vartheta }\sigma _3
\partial _\omega \Phi , \quad \rangle  \end{pmatrix} +   (P_c(\omega )P_c(\omega _0))^{-1} P_c(\omega
)e^{-\im \sigma _3\vartheta } =\\&  -(A,B)\mathcal{A}
\begin{pmatrix}d\omega
\\ d \vartheta    \end{pmatrix} +
 (P_c(\omega )P_c(\omega _0))^{-1}  P_c(\omega
)e^{-\im \sigma _3\vartheta } .\end{aligned}\nonumber
\end{equation}
By \eqref{eq:indentitiesGradf} we  obtain that $A$ and $B$ are
identified by the following equations (treating the last
$(P_c(\omega )P_c(\omega _0))^{-1} P_c(\omega )$ like a scalar):
\begin{equation}  \mathcal{A}^* \begin{pmatrix} A  \\ B \end{pmatrix} =
(P_c(\omega )P_c(\omega _0))^{-1}  P_c(\omega ) \begin{pmatrix}
\partial _\omega R
  \\ \im    \sigma _3 R
\end{pmatrix}.\nonumber
\end{equation}
\qed

\section{Symplectic structure}
\label{section:symplectic}

 Our ambient space is    $H^1( \R ^3, \C )\times H^1( \R ^3, \C )$.
We focus only on points with $\sigma _1U=\overline{U}$. The natural
symplectic structure for our problem is
\begin{equation}\label{eq:SymplecticForm}
  \Omega (X,Y)=\langle X, \sigma _3\sigma _1 Y \rangle .
\end{equation}
We will see that the coordinates we introduced in
\eqref{eq:coordinate}, which arise naturally from the linearization,
are not canonical   for \eqref{eq:SymplecticForm}. This is the main
difference with \cite{bambusicuccagna}. In this section we exploit
the work in section \ref{section:modulation} to compute the Poisson
brackets for pairs of coordinates. We end the section with a crucial
property for $Q$,
  Lemma
\ref{lem:InvarianceQ}.

The hamiltonian vector field $X_G$ of a scalar function $G$   is
defined by the equation $\langle X_G, \sigma _3\sigma _1 Y
\rangle=-\im \langle \nabla G ,  Y \rangle$ for any vector   $Y$ and
is $X_G=-\im \sigma _3\sigma _1 \nabla G  $. At $U=e^{\im \sigma
_3\vartheta } (\Phi _\omega + R)$ as in \eqref{eq:anzatz} we have by
\eqref{eq:gradient}
\begin{equation}\label{eq:HamVectorfield}\begin{aligned} &
 X_G(U)=\im
  \frac{\langle \nabla G (U), e^{ \im \sigma _3\vartheta }
  \sigma _3 \Phi \rangle }
  {q'(\omega )}
  e^{ \im \sigma _3\vartheta }\partial _\omega \Phi -\im
\frac{\langle \nabla G(U) ,
    e^{ \im \sigma _3\vartheta }\partial _\omega \Phi \rangle }
  {q'(\omega )}
  e^{ \im \sigma _3\vartheta }\sigma _3 \Phi
  \\&  +\im
   \sum _j\partial _{z_j}G(U) e^{ \im \sigma _3\vartheta }
   \sigma _1\xi _j -\im \sum _j\partial _{\overline{z}_j}G(U)
  e^{ \im \sigma _3\vartheta }\xi _j -\\&  -\im   e^{  \im \sigma _3\vartheta }\sigma _3\sigma _1
   P_c(\mathcal{H}_\omega ^* )e^{  \im \sigma _3\vartheta }\nabla G(U).
 \end{aligned}
\end{equation}
We call Poisson bracket of a pair  of scalar valued  functions $F$
and $G$  the  scalar valued   function
\begin{equation}\label{eq:PoissonBracket}
  \{ F,G \} = \langle \nabla F , X_ G   \rangle  =
  -\im \langle \nabla F ,\sigma _3\sigma _1\nabla  G   \rangle
= \im \Omega ( X_F,   X_G )
   .
\end{equation}
By $0= \im \frac{d}{dt} Q (U(t)) =\langle
 \nabla  Q (U(t)),\sigma _3\sigma _1\nabla E  (U(t)) \rangle  $ we have the
  commutation
\begin{equation}\label{eq:PoissonCommutation}
\{  Q,E \} =0.
\end{equation}
In terms of spectral components we have
\begin{equation}\label{eq:PoissonBracketComponent}
\begin{aligned} & \im \{  F,G \} (U)=\langle
  \nabla F (U) ,\sigma _3\sigma _1\nabla  G (U)   \rangle =
  (q') ^{-1}
  \times \\&
\big [ \langle \nabla F  ,
   e^{  \im \sigma _3\vartheta }\sigma _3 \Phi \rangle
   \langle \nabla G  ,
   e^{  \im \sigma _3\vartheta }\partial _\omega \Phi \rangle
- \langle \nabla F , e^{  \im \sigma _3\vartheta }\partial _\omega
\Phi \rangle
   \langle \nabla G  ,
   e^{  \im \sigma _3\vartheta }\sigma _3 \Phi \rangle
\big ]
    \\&  + \sum _j  \big [
 \partial _{z_j} F  \partial _{\overline{z}_j} G
-\partial _{\overline{z}_j} F  \partial _{ {z}_j} G
 \big ] +\\& +
   \langle \sigma _3 e^{ - \im \sigma _3\vartheta }P_c(\mathcal{H}_\omega ^*)
   e^{  \im \sigma _3\vartheta } \nabla F ,
   \sigma _1  e^{  -\im \sigma _3\vartheta }
   P_c(\mathcal{H}_\omega ^*) e^{  \im \sigma _3\vartheta }
   \nabla G  \rangle
 . \end{aligned}
\end{equation}

\begin{lemma}
\label{lem:PoissBrackCoord} Let $F(U)$ be a scalar function. We have
the following equalities:
\begin{eqnarray} & \{ \omega ,\vartheta  \} = \frac{  q'}{(q'
)^2 -\langle R,
\partial _\omega \Phi   \rangle ^2 + \langle \sigma _3R, \Phi
  \rangle \langle R ,\sigma _3\partial _\omega ^2  \Phi
 \rangle }\label{omegatheta};\\ & \{ z_j ,F \} =\langle \sigma _3
 \xi _j, \partial _\omega R \rangle \{ F, \omega \} +\im
  \langle \sigma _3
 \xi _j, \sigma _3 R \rangle \{ F, \vartheta \} -\im
 \partial _{\overline{z}_j}F ;\label{ZF} \\ & \{ \overline{z}_j ,F
 \} =\langle \sigma _1\sigma _3
 \xi _j, \partial _\omega R \rangle \{ F, \omega \} +\im
  \langle \sigma _1\sigma _3
 \xi _j, \sigma _3 R \rangle \{ F, \vartheta \} +\im
 \partial _{z_j}F .\label{barZF}
 \end{eqnarray}
In particular we have:

\begin{equation} \begin{aligned}
     &
\{  z_j, \omega   \} =\im
  \langle \sigma _3
 \xi _j, \sigma _3 R \rangle \{ \omega, \vartheta \}
\, ; \,  \{  \overline{z}_j, \omega   \} =\im
  \langle \sigma _1\sigma _3
 \xi _j, \sigma _3 R \rangle \{ \omega, \vartheta \} ; \\&
 \{  z_j, \vartheta   \} =
  \langle \sigma _3
 \xi _j, \partial  _\omega R \rangle \{ \vartheta , \omega  \}
\, ; \,  \{  \overline{z}_j, \vartheta  \} =
  \langle \sigma _1\sigma _3
 \xi _j,  \partial  _\omega R \rangle \{ \vartheta , \omega  \}
 ; \\
&
 \{ z_k,z_j  \}  =\im (
 \langle \sigma _3 \xi _k , \partial _\omega R \rangle \langle \sigma _3
 \xi _j, \sigma _3 R\rangle -\langle \sigma _3 \xi _j , \partial _\omega R \rangle \langle \sigma _3
 \xi _k, \sigma _3 R\rangle ) \{ \omega, \vartheta \} ;\\&
 \{ \overline{z}_k,\overline{z}_j  \}    =\im (
 \langle \sigma _1\sigma _3 \xi _k ,
 \partial _\omega R \rangle \langle \sigma _1\sigma _3
 \xi _j, \sigma _3 R\rangle -\langle
 \sigma _1\sigma _3 \xi _j ,
 \partial _\omega R \rangle \langle \sigma _1\sigma _3
 \xi _k, \sigma _3 R\rangle ) \{ \omega, \vartheta \} ;\\&
 \{ z_k,\overline{z}_j  \}
   =-\im \delta _{jk}+
 \im (\langle  \sigma _3 \xi _k , \partial _\omega R \rangle
 \langle \sigma _1\sigma _3
 \xi _j, \sigma _3 R\rangle -\langle \sigma _1\sigma _3 \xi _j ,
  \partial _\omega R \rangle \langle
 \xi _k,   R\rangle ) \{ \omega, \vartheta \}.
 \end{aligned}\nonumber\end{equation}
\end{lemma}
\proof  By \eqref{eq:GradModulation} and  \eqref{eq:PoissonBracketComponent}
we have $ \im  \{ \omega ,\vartheta \} = $ \begin{equation}\begin{aligned}&
 (q') ^{-1}
\big [ \langle \nabla \omega  ,
   e^{  \im \sigma _3\vartheta }\sigma _3 \Phi \rangle
   \langle \nabla \vartheta  ,
   e^{  \im \sigma _3\vartheta }\partial _\omega \Phi \rangle
- \langle \nabla \omega , e^{  \im \sigma _3\vartheta }\partial _\omega
\Phi \rangle
   \langle \nabla \vartheta  ,
   e^{  \im \sigma _3\vartheta }\sigma _3 \Phi \rangle
\big ] =\\& \frac{- \langle \sigma _3R, \Phi _\omega \rangle q' \langle R ,\sigma _3\partial _\omega ^2  \Phi
_\omega\rangle q'-[(q'(\omega))^2 -\langle R,
\partial _\omega \Phi _\omega \rangle ^2] (q' )^2  }{q' \im \left [ q'(\omega))^2 -\langle R,
\partial _\omega \Phi _\omega \rangle ^2 + \langle \sigma _3R, \Phi
_\omega \rangle \langle R ,\sigma _3\partial _\omega ^2  \Phi
_\omega\rangle\right ]^2}
   .\end{aligned}\nonumber
\end{equation}
This yields
\eqref{omegatheta}. For \eqref{ZF}, substituting \eqref{ZOmegaTheta}
in \eqref{eq:PoissonBracket}, we get $\{ z_j ,F \} =$
\begin{equation}\begin{aligned}& \langle \nabla z_j, X_F\rangle = - \langle \sigma _3 \xi _j, \partial _\omega R \rangle
  \{ \omega ,F \} - \im \langle \sigma _3 \xi _j, \sigma _3 R \rangle
  \{ \vartheta ,F \} +  \langle e^{-\im \sigma _3\vartheta }
  \sigma
_3\xi _j, X_F\rangle
.\end{aligned}\nonumber
\end{equation}
When we substitute $X_F$ with the decomposition in  \eqref{eq:HamVectorfield}, the last term in the above sum  becomes  $\langle e^{-\im \sigma _3\vartheta }
  \sigma
_3\xi _j, X_F\rangle =-\im    \partial _{\overline{z}_j}F\langle e^{-\im \sigma _3\vartheta }
  \sigma
_3\xi _j, e^{ \im \sigma _3\vartheta }\xi _j\rangle  =-\im    \partial _{\overline{z}_j}F  . $ This yields \eqref{ZF}.  \eqref{barZF} can be derived
by first replacing  $F$ with  $\overline{F}$  in \eqref{ZF} and by taking the
complex conjugate of the resulting equation:
 \begin{equation} \{ z_j ,F \} =\langle \sigma _3
 \xi _j, \partial _\omega \overline{R} \rangle \{ F, \omega \} -\im
  \langle \sigma _3
 \xi _j, \sigma _3 \overline{R} \rangle \{ F, \vartheta \} +\im
 \partial _{ {z}_j}F.
 \nonumber
\end{equation}
 Then \eqref{barZF} follows by  using that $\overline{R}=\sigma _1 R$ and $\sigma _1\sigma _3=-\sigma _3\sigma _1$.  The remaining formulas in the
 statement follow from  \eqref{ZF}--\eqref{barZF}.
  \qed

\begin{definition}\label{def:PoissonFunct}
Given a  function  $\mathcal{G}(U)$
 with
values in $L^2 _c(\mathcal{H} _{\omega _0}) $, a symplectic form
$\Omega$ and a scalar function $F(U)$, we define
\begin{equation}\label{fF}   \{ \mathcal{G}, F\} :=
\mathcal{G}'(U)X_F(U)
\end{equation}
with  $X_F$ the hamiltonian vector field associated to $F$. We set
  $\{ F,\mathcal{G} \} :=-\{ \mathcal{G}, F\}$.
\end{definition}
We have:

\begin{lemma}
\label{lem:fF} For $f(U)$ the functional in Lemma \ref{lem:gradient
f}, we have:
\begin{equation}\label{fF1} \{ f, F\} =
 (P_c(\omega )P_c(\omega _0))^{-1} P_c(\omega ) \left [ \{ F,\omega \}
 \partial _\omega R +\im \{ F, \vartheta   \} \sigma _3 R -\im
 e^{-\im \sigma _3 \vartheta} \sigma _3\sigma _1\nabla F\right ] .
\end{equation}
In particular we have:
\begin{equation}\label{fomegaZ} \begin{aligned}& \{ f, \omega\} =\im  \{ \omega, \vartheta
\} (P_c(\omega )P_c(\omega _0))^{-1} P_c(\omega )   \sigma _3 R
 ; \\&\{ f, \vartheta \} = \{   \vartheta ,\omega
\}(P_c(\omega )P_c(\omega _0))^{-1} P_c(\omega )   \partial _\omega
R;
 \\&  \{ f, z_j\} =
(P_c(\omega )P_c(\omega _0))^{-1} P_c(\omega ) \left [ \{ z_j,\omega
\}
 \partial _\omega R +\im \{ z_j, \vartheta   \} \sigma _3 R  \right ] ; \\&
 \{ f, \overline{z}_j\} =
(P_c(\omega )P_c(\omega _0))^{-1} P_c(\omega ) \left [ \{
\overline{z}_j,\omega \}
 \partial _\omega R +\im \{ \overline{z}_j, \vartheta   \} \sigma _3 R  \right ] .
 \end{aligned}
\end{equation}

\end{lemma}
\proof Using Lemma  \ref{lem:gradient f} and  by \eqref{eq:ApplmatrixA}
\begin{equation}\begin{aligned} & f' \sigma _3\sigma _1\nabla
F =-(A,B)\mathcal{A} \begin{pmatrix} \langle \nabla \omega ,\sigma
_3\sigma _1\nabla F \rangle
\\ \langle \nabla \vartheta ,\sigma
_3\sigma _1\nabla F \rangle \end{pmatrix}  \\& +(P_c(\omega
)P_c(\omega _0))^{-1} P_c(\omega )
 e^{-\im \sigma _3 \vartheta} \sigma _3\sigma _1\nabla F.
\end{aligned} \nonumber
\end{equation}
By Lemma \ref{lem:gradient f}  we have
\begin{equation}\begin{aligned}   (A,B)\mathcal{A}
 \begin{pmatrix}   \{ \omega ,    F \}
\\ \{   \vartheta ,   F \} \end{pmatrix}  =
 (P_c(\omega )P_c(\omega _0))^{-1}
P_c(\omega ) (  \partial _\omega R
 ,   \im \sigma _3 R) \begin{pmatrix}   \{ \omega ,    F \}
\\ \{   \vartheta ,   F \} \end{pmatrix}.
\end{aligned} \nonumber
\end{equation}
 \qed

The following result is important in the sequel.
\begin{lemma}
  \label{lem:InvarianceQ} Let $Q$ be the function defined in
  \eqref{eq:charge}.
Then, we have the following formulas:
\begin{eqnarray}
& \{ Q,\omega  \} =0 ; \label{Qomega}\\& \{ Q,\vartheta  \} =1;
\label{Qtheta}  \\& \{ Q,z_j  \} =\{ Q,\overline{z}_j  \}=0;
\label{QZ}\\& \{ Q,f  \}=0 . \label{Qf}
\end{eqnarray}
Denote by $X_Q$ the hamiltonian vectorfield of $Q$. Then

\begin{equation} \label{eq:Ham.VecFieldQ} X_Q=-
\frac{\partial}{\partial \vartheta} .
\end{equation}
\end{lemma}
\proof   We have by \eqref{eq:PoissonBracketComponent},
\eqref{eq:GradModulation} and $\nabla Q(U)=\sigma _1U$,
\begin{equation}
\begin{aligned} & \im q'\{ Q,\omega \}  =
  \langle \nabla Q  ,
   e^{  \im \sigma _3\vartheta }\sigma _3 \Phi \rangle
   \langle \nabla \omega ,
   e^{  \im \sigma _3\vartheta }\partial _\omega \Phi \rangle
- \langle \nabla Q , e^{  \im \sigma _3\vartheta }\partial _\omega
\Phi \rangle
   \langle \nabla \omega  ,
   e^{  \im \sigma _3\vartheta }\sigma _3 \Phi \rangle
    \\&
=  q'\frac{  -\langle R, \sigma _3 \Phi \rangle (q'(\omega)+ \langle
R,
\partial _\omega \Phi _\omega \rangle   )
-(q'(\omega)+ \langle R,
\partial _\omega \Phi _\omega \rangle   ) (-1)\langle R, \sigma _3 \Phi \rangle
}{(q'(\omega))^2 -\langle  R,
\partial _\omega \Phi _\omega \rangle ^2 + \langle \sigma _3R, \Phi
_\omega \rangle \langle R ,\sigma _3\partial _\omega ^2  \Phi
_\omega\rangle } =0. \end{aligned}\nonumber
\end{equation}
Similarly,
\begin{equation}
\begin{aligned} & \im q'\{ Q,\vartheta \}  =
  \langle \nabla Q  ,
   e^{  \im \sigma _3\vartheta }\sigma _3 \Phi \rangle
   \langle \nabla \vartheta ,
   e^{  \im \sigma _3\vartheta }\partial _\omega \Phi \rangle
- \langle \nabla Q , e^{  \im \sigma _3\vartheta }\partial _\omega
\Phi \rangle
   \langle \nabla \vartheta  ,
   e^{  \im \sigma _3\vartheta }\sigma _3 \Phi \rangle
    \\&
= q'\frac{  -\langle R, \sigma _3 \Phi \rangle \langle R, \sigma _3
\partial _\omega ^2\Phi \rangle -(q'(\omega)+ \langle R,
\partial _\omega \Phi _\omega \rangle   )  (q'(\omega)- \langle R,
\partial _\omega \Phi _\omega \rangle   )
}{\im [(q'(\omega))^2 -\langle  R,
\partial _\omega \Phi _\omega \rangle ^2 + \langle \sigma _3R, \Phi
_\omega \rangle \langle R ,\sigma _3\partial _\omega ^2  \Phi
_\omega\rangle  ]} =q' \im . \end{aligned}\nonumber
\end{equation}
By \eqref{ZF},\eqref{Qomega} and  \eqref{Qtheta} we have
\begin{equation}\label{eq:ZQ}
\begin{aligned} & \im \{ z_j ,Q \} = -
  \langle
 \xi _j,   R \rangle   +
 \partial _{\overline{z}_j}Q \\& \im \{ \overline{z}_j ,Q \} =
  \langle
 \xi _j,   \sigma _1 R \rangle   -
 \partial _{ {z}_j}Q. \end{aligned}
\end{equation}
By
\begin{equation} \label{eq:Qcoord} Q(U)=q +\frac{1}{2}
\langle z\cdot \xi +\overline{z}\cdot \sigma _1\xi +P_c(\omega
)f,\sigma _1(z\cdot \xi +\overline{z}\cdot \sigma _1\xi +P_c(\omega
)f)\rangle
\end{equation}
we have
\begin{equation} \label{eq:QZderivatives}
\partial _{ {z}_j}Q=\langle
 \xi _j,   \sigma _1 R \rangle \, , \quad
\partial _{\overline{z}_j}Q=\langle
 \xi _j,   R \rangle .
\end{equation}
So both lines in \eqref{eq:ZQ} are 0 and yield \eqref{QZ}.
 Finally \eqref{Qf} follows by \eqref{fF},
  Lemma \ref{lem:fF},  \eqref{Qomega} ,
  \eqref{Qtheta} and by
\begin{equation} \begin{aligned} &
 \{ f,Q \}  =
(P_c(\omega )P_c(\omega _0))^{-1}P_c(\omega ) \left [ \im \{ Q,
\vartheta   \} \sigma _3 R -\im
 e^{-\im \sigma _3 \vartheta} \sigma _3\sigma _1\nabla Q\right ]
 \\& = (P_c(\omega )P_c(\omega _0))^{-1}P_c(\omega )
 \left [ \im \sigma _3 R -\im \sigma _3\Phi -\im \sigma _3 R \right ] =0.
 \end{aligned}\nonumber
\end{equation}
\eqref{eq:Ham.VecFieldQ} is an immediate consequence of the
definition  of $X_Q$ and of \eqref{Qomega}--\eqref{Qf}.
 \qed
\section{Hamiltonian riformulation of the system}
\label{section:Hamiltonian riformulation}

  \eqref{eq:system} is how the problem is framed in the literature.
Yet \eqref{eq:system} hides the  crucial hamiltonian nature of the
problem. In the coordinate system \eqref{eq:coordinate}   can be written
as follows:
\begin{equation} \label{eq:SystPoiss} \begin{aligned} &
  \dot \omega  = \{ \omega , E \}   \, , \quad    \dot f= \{f, E
\}   \, , \\&  \dot z_j  = \{ z_j , E \} \, , \quad    {\dot
{\overline{z}}_j }=  \{ \overline{z}_j , E \}  \, ,
 \\ &     \dot \vartheta = \{ \vartheta , E  \}. \end{aligned}
\end{equation}
For the scalar coordinates the equations in \eqref{eq:SystPoiss} are
due to the hamiltonian nature of \eqref{eq:NLSvectorial}. Exactly for the same reasons we have  the equation $ \dot f= \{f, E \}  $, which we now derive in
the following standard way. Multiplying \eqref{eq:system} by $e^{\im
\vartheta\sigma _3}$ one can rewrite \eqref{eq:system} by
\eqref{eq:vectorfields} and \eqref{eq:gaugeInvariance},
 as \begin{equation}
\label{eq:sys1}
\begin{aligned} &
  -\im
 \dot \vartheta
  \frac{\partial}{\partial \vartheta }+  \im \dot \omega
 \frac{\partial}{\partial \omega }  +  \im
  \sum _j \dot z_j\frac{\partial}{\partial z _j } +  \im
  \sum _j   \dot {\overline{z}}_j\frac{\partial}{\partial
  \overline{z} _j } \\& +\im
   e^{\im \vartheta\sigma
_3} P_c(\mathcal{ H}_\omega ) \dot f  =  \sigma _3
    \sigma _1
  \nabla E  (  U)    .\end{aligned}
\end{equation}
When we apply the  derivative $f'(U)$ to \eqref{eq:sys1} the first
line cancels, so that
\begin{equation}
\label{eq:sys2}
\begin{aligned} & \dot f= f'(U)
   e^{\im \vartheta\sigma
_3} P_c(\mathcal{ H}_\omega ) \dot f  =- f'(U)\im \sigma _3
    \sigma _1
  \nabla E  (  U) =f'(U)X_E(U)=\{ f,E\} .
\end{aligned}\nonumber
\end{equation}
where the first equality is \eqref{eq:indentitiesGradf}, the third
the definition of hamiltonian field two lines above
\eqref{eq:HamVectorfield} and the last equality is Definition
\ref{def:PoissonFunct}.

 We now introduce a new hamiltonian. For $u_0$ the initial datum in
\eqref{NLS}, set
\begin{equation}    \label{eq:K}  \begin{aligned} &
K(U)=E(U)+\omega (U)   Q(U)-\omega (U)\| u_0\| _{L^2_x}^{2}.
\end{aligned}
\end{equation}
By Lemma \ref{lem:InvarianceQ}  the solution of the initial value
problem in \eqref{NLS}    solves also
\begin{equation} \label{eq:SystPoissK} \begin{aligned} &
  \dot \omega  = \{ \omega , K \}   \, , \quad    \dot f= \{f, K
\}  \, ,  \\&  \dot z_j  = \{ z_j , K \} \, , \quad    {\dot
{\overline{z}}_j }=  \{ \overline{z}_j , K \}  \, ,
 \\ &      \dot \vartheta  -\omega  = \{ \vartheta , K  \}. \end{aligned}
\end{equation}
By $ \frac{\partial}{\partial \vartheta} K  =0$  the right hand
sides in the equations \eqref{eq:SystPoissK} do not depend on
$\vartheta$. Hence, if we look at the new system

\begin{equation} \label{eq:SystK} \begin{aligned} &
  \dot \omega  = \{ \omega , K \}   \, , \quad    \dot f= \{f, K
\}  \, ,  \\&   \dot z_j  = \{ z_j , K \} \, , \quad   {\dot
{\overline{z}}_j }=  \{ \overline{z}_j , K \}  \, ,
 \\ &       \dot \vartheta  = \{ \vartheta , K  \}, \end{aligned}
\end{equation}
the evolution of the crucial variables $(\omega , z, \overline{z},
f)$ in \eqref{eq:SystPoiss} and \eqref{eq:SystK} is the same.
Therefore, to prove Theorem \ref{theorem-1.1} it is sufficient to
consider system \eqref{eq:SystK}.

\section{Application of the   Darboux Theorem}
\label{section:Darboux}

Since the main obstacle at reproducing the Birkhoff normal forms
argument of \cite{bambusicuccagna} for \eqref{eq:SystK} is that the
coordinates \eqref{eq:coordinate} are not canonical, we change
coordinates. That is,  we apply the Darboux Theorem. We warn the
reader not to confuse the variable $t\in [0,1]$ of this section with
the time of the evolution equation of the other sections.

 We introduce the  2-form, for
$q=q(\omega )=\| \phi _\omega \| ^{2}_{L^2_x}$ and summing on repeated indexes,
\begin{equation} \label{eq:Omega0} \Omega _0=\im d\vartheta \wedge
dq +   dz_j\wedge d\overline{z}_j+\langle f' (U )\quad , \sigma _3
\sigma _1 f' (U )\quad \rangle ,
\end{equation}
with $f (U)$ the function in Lemma \ref{lem:gradient f},  $f '(U)$
its Frech\'et derivative and the last term in \eqref{eq:Omega0}
acting on pairs $(X,Y)$ like $\langle f' (U )X  , \sigma _3 \sigma
_1 f' (U )Y\rangle $.
  It is an elementary exercise to show that $\Omega _0  $  is a
closed and non degenerate 2 form. In Lemma \ref{lem:OmegaOmega0} we
check that $\Omega _0(U) =\Omega  (U)$ at $U=e^{\im \sigma
_3\vartheta}
  \Phi _{\omega_{0}}$. Then the proof of the Darboux
Theorem
  goes as follows. One first considers
\begin{equation} \label{eq:Omegat} \Omega _t =(1-t)\Omega _0+t
\Omega =\Omega _0 +t\widetilde{\Omega} \text{ with
$\widetilde{\Omega}
 :=\Omega -\Omega _0$.}  \end{equation}
Then one considers a 1- differential form $\gamma (t,U)$ such that
 (external differentiation will always be on the $U$ variable only)
 $\im d \gamma (t,U) = \widetilde{\Omega} $ with $\gamma (U) =0$
 at $U=e^{\im
\sigma _3\vartheta}
  \Phi _{\omega_{0}}$.  Finally one considers the vector field
   $\mathcal{Y}^t$ such that
 $i_{\mathcal{Y}^t}\Omega _t=-\im \gamma $ (here for $\Omega$ a 2 form and $Y$ a vector field, $i_Y\Omega $ is the 1 form defined by $i_Y\Omega (X):=\Omega (Y,X))$
 and   the
 flow $\mathfrak{F}_t$ generated by $\mathcal{Y}^t$, which near
 the points
$e^{\im \sigma _3\vartheta}
  \Phi _{\omega_{0}}$ is defined up to time 1, and show that
$ \mathfrak{F}_1^*\Omega =\Omega _0$  by
\begin{equation} \label{eq:dartheorem}\begin{aligned} &\frac{d}{dt}
\left ( \mathfrak{F}_t^*\Omega _t\right )
  =\mathfrak{F}_t^*\left ( L_{\mathcal{Y}_t} \Omega _t\right )
  +  \mathfrak{F}_t^*\frac{d}{dt}\Omega _t =
  \\& = \mathfrak{F}_t^*d\left (
  i_{\mathcal{Y}^t} \Omega _t\right )
  +  \mathfrak{F}_t^*   \widetilde{\Omega }  =
  \mathfrak{F}_t^*\left ( -\im d \gamma
  +  \widetilde{\Omega}\right ) =0
 .
\end{aligned}
\end{equation}
For $\Omega _0$, the coordinates \eqref{eq:coordinate} are
canonical. But if one does not choose the 1 form $\gamma$ carefully,
then the new hamiltonian $\widetilde{K}= K\circ \mathfrak{F}_1$ will
not yield a semilinear NLS for coordinates \eqref{eq:coordinate},
which is what we need to perform the argument of
\cite{bambusicuccagna,cuccagnamizumachi}.   In the sequel of this
section  all the work is finalized to the correct choice if
$\gamma$. In Lemma \ref{lem:1forms} we compute explicitly a
differential form $\alpha$ and we make the preliminary choice
$\gamma =-\im \alpha$. This is not yet the right choice. By the
computations in Lemma \ref{lem:linearAlgebra} and Remark
\ref{rem:correction}, we find the obstruction to the fact that
$\widetilde{K}$ is of the desired type. Lemmas
\ref{lem:HamThetaOmega}--\ref{lem:flow Htheta} are necessary to find
an appropriate solution   $F$ of a differential equation in Lemma
\ref{lem:correction alpha}. Then $\gamma =-\im \alpha +  dF$ is the
right choice of $\gamma$. In Lemma \ref{lem:flow1} we collect a
number of useful estimates for $\mathfrak{F}_1$. Finally, Lemma
\ref{lem:flow2} is valid independently of the precise $\gamma$
chosen and contains information necessary for
\eqref{eq:newH}--\eqref{eq:SystK1}.

For any vector $Y\in T _U L^2$ we set
\begin{equation} \label{eq:Y}  \begin{aligned} &
Y=Y _{\vartheta}\frac{\partial}{\partial \vartheta}+Y
_{\omega}\frac{\partial}{\partial \omega} +\sum Y
_{j}\frac{\partial}{\partial z_j}+\sum Y
_{\overline{j}}\frac{\partial}{\partial \overline{z}_j} +e^{\im
\sigma _3 \vartheta} P_c (\omega )Y _{f }
\end{aligned}\end{equation}
for
\begin{equation} \label{eq:Y1}  \begin{aligned} & Y_{\vartheta}
=d\vartheta (Y)\, , \quad  Y_{\omega} =d\omega (Y)\, , \quad Y_{j}
=dz_j (Y)\\& \quad Y_{\overline{j}} =d\overline{z}_j (Y) \, , \quad
\quad Y_{f } =f' (U)  Y .
\end{aligned}\end{equation}
Similarly, a  differential 1-form $\gamma $ decomposes as
\begin{equation} \label{eq:gamma}  \begin{aligned} &
\gamma=\gamma ^{\vartheta}d \vartheta +\gamma ^{\omega}d \omega
+\sum \gamma ^{j}d z_j +\sum \gamma ^{\overline{j}}d\overline{z}_j
+\langle \gamma ^{f }, f'\quad \rangle  ,
\end{aligned}\end{equation}
where: $\langle \gamma ^{f }, f'\quad \rangle $ acts on a vector $Y$
as $\langle \gamma ^{f }, f'Y \rangle $, with here $\gamma ^f \in
L^2_c(\mathcal{H}_{\omega _0}^*)$; $\gamma ^{\vartheta}$, $\gamma
^{\omega}$, $\gamma ^{j}$ and  $\gamma ^{\overline{j}}$ are in $\C$.
Notice that we are reversing the standard notation on super and
 subscripts for forms and vector fields. In    the sequel, given a differential
  1 form $\gamma $ and
 a point $U$, we will denote by $\gamma _U$ the value of
$\gamma $ at $U$.

Given a function $\chi $, denote its hamiltonian vector field with
respect to $\Omega _t$
 by $X^t_\chi$ :  $i _{X^t_\chi}\Omega _t=-\im \, d\chi$.
By \eqref{eq:Omega0}  we have the following   hamiltonian
vectorfield associated to $q(\omega )$ (this is important in Lemma
\ref{lem:flow2} later):
\begin{equation} \label{eq:HamVect0q} X_{q(\omega )}^{0}=
-\frac{\partial}{\partial \vartheta}.
\end{equation}
We have the following preliminary observation:

\begin{lemma}
  \label{lem:OmegaOmega0} At $U=e^{\im \sigma _3\vartheta}
  \Phi _{\omega_{0}}$, for any $\vartheta $, we have $\Omega _0(U)=\Omega
(U)$.
\end{lemma}
\proof  Using the following partition of the identity
\begin{equation}\label{eq:partitionUno}\uno =
e^{\im \sigma _3 \vartheta}
[ P_{N_g(\mathcal{H}_{\omega})} +\sum _{\lambda \in \sigma
(\mathcal{H}_{\omega}) \backslash \{ 0 \} } P _{\ker
(\mathcal{H}_{\omega}-\lambda )} +P_{c} (\mathcal{H}_{\omega})  ]
e^{-\im \sigma _3 \vartheta}
\end{equation}
we get, summing on repeated indexes,
\begin{equation}\label{eq:OmegaComponent1}
\begin{aligned} & \Omega (X,Y)=\langle
 X ,\sigma _3\sigma _1Y  \rangle =
   \\&
\frac{1}{q'}\big [ \langle  X ,  e^{ - \im \sigma _3\vartheta
}\sigma _3\partial _\omega \Phi \rangle
   \langle   Y  ,e^{ - \im \sigma _3\vartheta }
     \Phi \rangle -\langle    X ,
  e^{ - \im \sigma _3\vartheta }  \Phi \rangle
   \langle  Y  ,e^{ - \im \sigma _3\vartheta }
   \sigma _3 \partial _\omega \Phi \rangle
  \big ] +
    \\&      \big [
 \langle  X ,   e^{ - \im \sigma _3\vartheta }\sigma _3  \xi _j
 \rangle \langle  Y , e^{ - \im \sigma _3\vartheta }\sigma _1
 \sigma _3 \xi _j
 \rangle -\langle  X , e^{ - \im \sigma _3\vartheta }\sigma _1
  \sigma _3 \xi _j
 \rangle \langle  Y ,   e^{ - \im \sigma _3\vartheta }\sigma _3 \xi _j
 \rangle
 \big ]  \\& +
   \langle  P_c(\mathcal{H}_\omega  )
   e^{  -\im \sigma _3\vartheta } X ,
   \sigma _3 \sigma _1
   P_c(\mathcal{H}_\omega  ) e^{ - \im \sigma _3\vartheta }
   Y  \rangle
 . \end{aligned}
\end{equation}
By \eqref{eq:ApplmatrixA} we have
\begin{equation}\label{eq:Omega1}
\begin{aligned} & \langle  \quad ,  e^{ - \im \sigma _3\vartheta
}\sigma _3\partial _\omega \Phi \rangle
 \wedge  \langle   \quad  ,e^{ - \im \sigma _3\vartheta }
     \Phi \rangle =\det \mathcal{A}\,  d\omega \wedge  d\vartheta .
\end{aligned}
\end{equation}
Substituting \eqref{ZOmegaTheta}--\eqref{barZOmegaTheta} we get
\begin{equation}\label{eq:Omega2}
\begin{aligned} & \langle  \quad ,  e^{ - \im \sigma _3\vartheta
}\sigma _3\xi _j\rangle
 \wedge  \langle   \quad  ,e^{ - \im \sigma _3\vartheta }
    \sigma _1 \sigma _3\xi _j \rangle =\\&    (dz_j + \langle \sigma _3 \xi _j, \partial _\omega R \rangle
 d \omega + \im \langle \sigma _3 \xi _j, \sigma _3 R \rangle
  d \vartheta ) \\&  \wedge  ( d\overline{z}_j  +\langle \sigma _1\sigma _3 \xi _j,
\partial _\omega R \rangle d \omega +
   \im \langle \sigma _1\sigma _3 \xi _j, \sigma _3 R \rangle
   d \vartheta ) .
\end{aligned}
\end{equation}
By Lemma \ref{lem:gradient f} we have
\begin{equation}\label{eq:Omega3}
\begin{aligned} &  \langle  P_c(\mathcal{H}_\omega  )
   e^{  -\im \sigma _3\vartheta }\quad ,
   \sigma _3 \sigma _1
   P_c(\mathcal{H}_\omega  ) e^{ - \im \sigma _3\vartheta }
   \quad  \rangle =\\&    \langle P_c(\omega )P_c(\omega _0)f'\quad +  P_c(\omega )\partial _\omega R \, d\omega +\im
  P_c(\omega )\sigma _3 R  \, d\vartheta ,\\&  \sigma _3 \sigma _1 (P_c(\omega )P_c(\omega _0)f'\quad  +  P_c(\omega )\partial _\omega R \, d\omega +\im
  P_c(\omega )\sigma _3 R  \, d\vartheta)  \rangle
    .
\end{aligned}
\end{equation}

Then by \eqref{eq:OmegaComponent1}--\eqref{eq:Omega3}  we have

\begin{equation}\label{eq:OmegaComponent2}
\begin{aligned} & \Omega  =   ( \im q'+a_1) d\vartheta \wedge d\omega
+   dz_j\wedge d\overline{z}_j+ \\& +     dz_j\wedge \left ( \langle
\sigma _1 \sigma _3 \xi _j  , \partial _\omega R  \rangle \, d\omega
+ \im \langle \sigma _1 \sigma _3 \xi _j ,\sigma _3  R \rangle \,
d\vartheta \right ) \\& -    d\overline{z}_j\wedge \left ( \langle
\sigma _3 \xi _j  , \partial _\omega R \rangle \, d\omega + \im
\langle   \sigma _3 \xi _j ,\sigma _3  R \rangle \, d\vartheta
\right ) +\\& +\langle P_c(\omega ) P_c(\omega _0)f' \quad , \sigma
_3 \sigma _1 P_c(\omega ) P_c(\omega _0) f'
 \quad \rangle +\\& +
 \langle P_c(\omega ) P_c(\omega _0)f'  \quad , \sigma
_3 \sigma _1 P_c(\omega  )
 \partial _\omega R \rangle \wedge  d\omega + \\&
+  \im \langle P_c(\omega ) P_c (\omega _0)f'  \quad , \sigma _3
\sigma _1 P_c(\omega )
 \sigma _3  R \rangle \wedge  d\vartheta
 , \end{aligned}
\end{equation}
where

\begin{equation}\label{eq:Omega4} \begin{aligned} & \im q'+a_1= \frac{\det \mathcal{A}}{q'}  + \langle P_c(\omega ) \partial _\omega R, \sigma _3 \sigma _1 P_c(\omega )\sigma _3  \im R \rangle \\& + \langle \sigma _3\xi _j,
\partial _\omega R \rangle  \langle \sigma _1\sigma _3\xi _j,\im \sigma _3R
  \rangle  -
  \langle \sigma _1\sigma _3\xi _j,  \partial _\omega R
  \rangle
   \langle \sigma _3\xi _j, \im \sigma _3R
 \rangle .
\end{aligned}
\end{equation}
In particular we have
\begin{equation}\label{eq:a1}
\begin{aligned} & a_1:= -\im q'+\frac{\det \mathcal{A}}{q'} +
 \langle P_{N^\perp _g(\mathcal{H}^{*}_{\omega})} \im\sigma _3
R, \sigma _3 \sigma _1    \partial _\omega R\rangle .
\end{aligned}
\end{equation}
Notice that $a_1=a_1(\omega , z,f)$ is smooth in the arguments
$\omega \in \mathcal{O}$, $z\in \C ^n$ and $f\in H^{-K',-S'}$ for any pair
$(K', S')$  with, for  $( z,f)$ near 0,
\begin{equation}\label{eq:bounda1}
\begin{aligned} & |a_1| \le C (K',S') (|z|+\| f\| _{H^{-K',-S'} } )^2 .
\end{aligned}
\end{equation}
At points $U=e^{\im \sigma _3 \vartheta} \Phi _\omega$, that is for
$R=0$,  we have
\begin{equation}\label{eq:OmegaComponent3}\begin{aligned} &
\Omega  =\im d\vartheta \wedge dq +    dz_j\wedge
d\overline{z}_j+\langle P_c(\omega ) P_c(\omega _0)f' \quad  ,
\sigma _3 \sigma _1 P_c(\omega ) P_c(\omega _0) f'
 \quad \rangle
 . \end{aligned} \nonumber
\end{equation}

At $\omega =\omega _0$ we get $\Omega =\Omega _0$.
 \qed

\begin{lemma}
  \label{lem:1forms} Consider the following  forms:
\begin{equation} \label{eq:1forms}\begin{aligned} &
\beta (U)Y:=\frac{1}{2}\langle \sigma _1\sigma _3U , Y\rangle ;\\&
\beta _0(U):=-\im qd\vartheta - \sum  _j\frac{\overline{z}_j dz_j -
{z}_j d\overline{z}_j}{2} +\frac{1}{2}\langle f (U),\sigma _3\sigma
_1f '(U)\quad \rangle.
\end{aligned}
\end{equation}
Then \begin{equation} \label{eq:1formsExtDiff} d\beta _0=\Omega _0\,
, \quad d\beta =\Omega   .\end{equation} Set
\begin{equation} \label{eq:alpha1} \alpha (U)=\beta (U)-\beta _0(U)
+d\psi (U)\text{ where } \psi (U):=\frac{1}{2}\langle \sigma _3\Phi
  , R\rangle .
\end{equation}
We have  $\alpha  = \alpha ^{\vartheta} d\vartheta +\alpha ^{\omega}
   d\omega
 + \langle \alpha ^f,f'\rangle$ with:

\begin{equation} \label{eq:alpha2} \begin{aligned}
\alpha ^{\vartheta} +\frac{\im}{2}\| f\| _2^2   =&
 -\frac{\im}{2}     \| z\cdot \xi +\overline{z}\cdot
 \sigma _1\xi \| _2^2-\im \langle z\cdot \xi +\overline{z}\cdot
 \sigma _1\xi , \sigma _1P_c(\omega ) f\rangle
\\& - \frac{{\im}}{2} \langle (P_c(\omega )
 -P_c(\omega _0 )) f, \sigma _1(P_c(\omega )
 +P_c(\omega _0 ))f \rangle ;
\\    \alpha ^{\omega}  =&   -\frac{1}{2}
 \langle \sigma _1R, \sigma _3\partial _\omega
R\rangle  ;  \\    \alpha ^{f} = &\frac{1}{2}\sigma_1 \sigma _3P_c(\omega
_0)\left ( P_c(\omega ) -P_c(\omega _0)\right )f.
\end{aligned}
\end{equation}

\end{lemma}
\proof Everything  is straightforward   except for
\eqref{eq:alpha2}, which we now prove. We will sum over repeated
indexes. We substitute $U$ using \eqref{eq:coordinate}
getting
\begin{equation} \label{eq:beta} \begin{aligned} &
\beta =\frac{1}{2}\langle e^{-\im \sigma _3\vartheta}\sigma _1\sigma
_3 \Phi  , \quad \rangle  + \frac{1}{2}\langle e^{-\im \sigma
_3\vartheta}\sigma _1\sigma _3P_c(\omega )f, \, \rangle +\\&
\frac{1}{2}  \left [z_j \langle e^{-\im \sigma _3\vartheta}\sigma
_1\sigma _3\xi _j, \quad
 \rangle -\overline{z}_j \langle e^{-\im
\sigma _3\vartheta} \sigma _3\xi _j, \quad \rangle \right ] .
\end{aligned}
\end{equation}
When we decompose $\frac{1}{2}  e^{-\im \sigma _3\vartheta}\sigma
_1\sigma _3 \Phi $ like $\nabla F$ in \eqref{eq:gradient}, we obtain
\begin{equation} \label{eq:betatilde} \begin{aligned} &
 \frac{1}{2}\langle e^{-\im \sigma _3\vartheta}\sigma _1\sigma
_3 \Phi  , \quad \rangle =  -\frac{q}{q'}  \langle e^{-\im \sigma
_3\vartheta} \sigma _3 \partial _\omega \Phi  , \quad \rangle
\\& -\frac{1}{2}   \langle \sigma _3\Phi ,   \xi _j
\rangle  \left ( \langle e^{-\im \sigma _3\vartheta} \sigma _3\xi
_j, \quad \rangle  -\langle e^{-\im \sigma _3\vartheta} \sigma
_1\sigma _3\xi _j, \quad \rangle \right ) \\& - \frac{1}{2}\langle
e^{-\im \sigma _3\vartheta}P_c (\mathcal{H}_\omega ^{*})\sigma _3
\Phi , \quad \rangle \,
\end{aligned}
\end{equation}
with by \eqref{eq:ApplmatrixA}
\begin{equation} \label{eq:bettild}-\frac{q}{q'}\langle e^{-\im \sigma
_3\vartheta} \sigma _3 \partial _\omega \Phi  , \quad \rangle
=\frac{q}{q'}\langle R, \sigma _3 \partial _\omega ^2 \Phi \rangle
\, d \omega -\im \, \frac{q}{q'}   \, (q'   +\langle R,   \partial
_\omega \Phi \rangle ) \, d \vartheta .
\end{equation}
  Substituting slightly manipulated versions of   the formulas in Lemmas \ref{lem:gradient z}--\ref{lem:gradient f}, in particular using $\sigma _3P_c(\omega )= P_c(\omega )^*\sigma _3$, $\sigma _1P_c(\omega )= P_c(\omega ) \sigma _1$ and $\sigma _1\sigma _3=-\sigma _1\sigma _3 $,
  and summing over repeated indexes, we get
\begin{equation} \label{eq:beta0} \begin{aligned} &
\beta _0=-\im q\, d\vartheta +  \frac{1}{2} \overline{z}_j (  \langle \sigma _1 \sigma _1 \xi _j, \sigma _3\partial _\omega R \rangle
  d \omega + \im \langle   \xi _j,  R \rangle
  d \vartheta -  \langle e^{-\im \sigma _3\vartheta } \sigma _3 \xi _j, \quad \rangle ) \\&  + \frac 12
{z}_j \, (
\langle \sigma _1 \xi _j, \sigma _3
\partial _\omega R \rangle d \omega +
   \im \langle   \xi _j, \sigma _1 R \rangle
   d \vartheta  + \langle  e^{-\im \sigma _3\vartheta }
  \sigma _1\sigma
_3\xi _j, \quad \rangle ) \\& +\frac{1}{2}\langle f  ,\sigma
_3\sigma _1 (1-P_c(\omega ) P_c(\omega _0)) f ' \quad \rangle  + \frac{1}{2}\langle f  ,\sigma
_3\sigma _1P_c(\omega ) e^{-\im
\sigma _3\vartheta } \quad \rangle \\& +  \frac{1}{2}\langle \sigma _1 P_c(\omega ) f  ,\sigma
_3 \partial _\omega R \rangle  d\omega
+  \frac{ \im}{2}\langle P_c(\omega )f  ,  \sigma _1 R \rangle  d\vartheta  .
\end{aligned}\nonumber
\end{equation}
Hence
\begin{equation} \label{eq:beta0} \begin{aligned} &
 \beta _0=\im \left ( -q +\frac{1}{2}\langle R, \sigma _1R\rangle \right
) \, d \vartheta    + \frac{1}{2}\langle \sigma _1R, \sigma
_3\partial _\omega R\rangle \, d \omega +
\\& +
  \frac{1}{2}\langle \sigma _1\sigma _3  \left ( 1 -  P_c(\omega
 _0)P_c(\omega )  \right )   f,  f'\, \rangle + \\& + \frac{1}{2}  \left (z_j
\langle e^{-\im \sigma _3\vartheta}\sigma _1\sigma _3\xi _j, \quad
 \rangle -\overline{z}_j \langle e^{-\im
\sigma _3\vartheta} \sigma _3\xi _j,    \quad \rangle \right )+
\\& +
\frac{1}{2}\langle e^{-\im \sigma _3\vartheta} \sigma _1\sigma
_3P_c(\omega )f,\quad \rangle    .
\end{aligned}
\end{equation}
By \eqref{eq:coordinate} we have

    \begin{equation} \label{eq:dPsi0} \begin{aligned} & d\psi =
\frac{1}{2}\langle \sigma _3\Phi , \partial _\omega R \rangle
d\omega +\frac{1}{2}  \langle \sigma _3\Phi , \xi _j \rangle \left (
dz_j -d\overline{z}_j\right ) +\frac{1}{2}\langle \sigma _3\Phi ,
P_c(\omega ) f' \quad \rangle   .
\end{aligned}
\end{equation}
Applying to \eqref{eq:dPsi0}   Lemmas
 \ref{lem:gradient z}--\ref{lem:gradient f},  the fact that, in particular, we have
 \begin{equation} \begin{aligned} &
P_c(\omega )f'(U)= P_c(\omega )P_c(\omega _0)f'(U)=  P_c(\omega )\left [
 -   \partial _\omega R \, d\omega -\im
  \sigma _3 R  \, d\vartheta +  e^{-\im
\sigma _3\vartheta }\uno \right ] ,
\end{aligned} \nonumber
\end{equation}
and the identities
\eqref{eq:cancel0}--\eqref{eq:cancel1} below,
we get $d\psi =$
\begin{equation} \label{eq:dPsi} \begin{aligned} &   = \frac{1}{2}  \langle \sigma
_3\Phi ,   \xi _j \rangle  \left ( \langle e^{-\im \sigma
_3\vartheta} \sigma _3\xi _j, \quad \rangle -\langle e^{-\im \sigma
_3\vartheta} \sigma _1\sigma _3\xi _j, \quad \rangle \right )
\\& + \frac{1}{2}\langle e^{-\im \sigma _3\vartheta}P_c
(\mathcal{H}_\omega ^{*})\sigma _3 \Phi , \quad \rangle
\\& +\frac{q}{q'}  \langle \sigma _3\partial _\omega\Phi , \partial
_\omega R \rangle d\omega   -\frac{\im }{2} \langle \sigma _3 \Phi ,
P _{N_g^\perp (\mathcal{H}^*_\omega )}\sigma _3 R \rangle d\vartheta
 .
\end{aligned}
\end{equation}
To get the last line of \eqref{eq:dPsi} we have used:

\begin{equation} \label{eq:cancel0} \begin{aligned} & \frac{1}{2} \langle \sigma _3 \Phi , \partial
_\omega R \rangle - \frac{1}{2}    \langle \sigma _3 \Phi , \xi _j \rangle  \langle \sigma _3\xi _j , \partial
_\omega R \rangle \\& - \frac{1}{2}    \langle \sigma _3 \Phi , \sigma _1\xi _j \rangle  \langle \sigma _1\sigma _3\xi _j , \partial
_\omega R \rangle - \frac{1}{2}      \langle  \sigma _3\Phi , P_c(\omega )\partial
_\omega R \rangle   =\frac{1}{2} \langle \sigma _3 \Phi , \partial
_\omega R \rangle \\& - \frac{1}{2}\left [ \langle \sigma _3 \Phi , \partial
_\omega R \rangle  -\frac{1}{q'} \langle \sigma _3 \Phi , \sigma _3 \Phi  \rangle  \langle \sigma _3\partial _\omega\Phi , \partial
_\omega R \rangle   \right ] =\frac{2q}{2q'}\langle \sigma _3\partial _\omega\Phi , \partial
_\omega R \rangle ;
\end{aligned}
\end{equation}

\begin{equation} \label{eq:cancel1} \begin{aligned} &   - \frac{\im }{2}    \langle \sigma _3 \Phi , \xi _j \rangle  \langle \sigma _3\xi _j , \sigma _3 R \rangle   - \frac{\im }{2}    \langle \sigma _3 \Phi , \sigma _1\xi _j \rangle  \langle \sigma _1\sigma _3\xi _j ,\sigma _3 R \rangle \\& - \frac{\im }{2}      \langle  \sigma _3\Phi , P_c(\omega )\sigma _3 R \rangle    =-\frac{\im }{2} \langle \sigma _3 \Phi ,P _{N_g^\perp (\mathcal{H}^*_\omega )} \sigma _3 R \rangle .
\end{aligned}
\end{equation}

 Let us consider the sum \eqref{eq:alpha1}. There are various cancelations.
 The first and second   (resp.
the first term of the third) line of \eqref{eq:dPsi} cancel  with
the second and third lines of \eqref{eq:betatilde} (resp. the first
term of the rhs of
  \eqref{eq:bettild}). The last three terms in
rhs\eqref{eq:beta} cancel with the last two lines of \eqref{eq:beta0}.
The $-\im qd\vartheta $ term in the rhs of \eqref{eq:beta0}  cancels
with the $-\im qd\vartheta $ term in \eqref{eq:bettild}. Adding
the second term of the third line of \eqref{eq:dPsi} with the last
term in the  rhs of \eqref{eq:bettild} we get  the product of $\im  $ times the
following quantities:

\begin{equation} \label{eq:cancel} \begin{aligned} &
 -\frac{1}{2}   \langle \sigma _3
\Phi , P_{N^\perp _g(\mathcal{H}^*_\omega)} \sigma _3R\rangle -
\frac{q}{q'} \langle R,   \partial _\omega \Phi \rangle =-
\frac{1}{2}   \langle   \Phi ,  R\rangle \\& +\frac{1}{2} \langle \sigma
_3 \Phi , P_{N  _g(\mathcal{H} _\omega)} \sigma _3R\rangle   -
\frac{q}{q'} \langle R,   \partial _\omega \Phi \rangle \\& = -
\frac{1}{2}  \langle   \Phi ,
R \rangle+
 \frac{1}{2q'} \langle  \sigma _3R,
  \Phi \rangle  \langle  \sigma _3\Phi ,
\partial _\omega \Phi \rangle \\& +
 \frac{1}{2q'} \langle  \sigma _3R,
 \sigma _3 \partial _\omega \Phi \rangle  \langle  \sigma _3\Phi ,
\sigma _3 \Phi \rangle - \frac{q}{q'} \langle R,   \partial _\omega
\Phi \rangle =0,
\end{aligned}
\end{equation}
where for the second equality we have used
\begin{equation} P_{N  _g(\mathcal{H} _\omega)}=\frac{1}{q'} \sigma _3 \Phi
\langle \sigma _3 \partial _\omega \Phi , \quad \rangle +\frac{1}{q'}\partial _\omega \Phi
\langle  \Phi , \quad \rangle .\nonumber
\end{equation}
The last equality in \eqref{eq:cancel} can be seen as follows. The two terms in the third line in \eqref{eq:cancel}  are both equal to 0. Indeed,  $\langle  \sigma _3\Phi ,
\partial _\omega \Phi \rangle =0$ and, by $
 R\in N^{\perp}_g (\mathcal{H}_\omega ^*)$
  and $\Phi \in   N _g (\mathcal{H}_\omega ^*) $,
    $ \langle R,     \Phi \rangle =0$. The two terms in
the fourth line in \eqref{eq:cancel}  cancel each other. Then we get formulas for $\alpha ^{\omega} $ and $\alpha ^{f}$. We get $\alpha ^{\vartheta} $ also by
\begin{equation} \| P_c(\omega )f\| _2^2=\|  f\| _2^2 + \langle (P_c(\omega )
 -P_c(\omega _0 )) f, \sigma _1(P_c(\omega )
 +P_c(\omega _0 ))f \rangle .\nonumber
\end{equation}  \qed

 We have, summing over repeated indexes (also on  $j$
 and $\overline{j}$):
\begin{lemma}
  \label{lem:linearAlgebra} We have
\begin{equation} \label{eq:linAlg0}\begin{aligned} &
i _{Y  }\Omega _0=\im q'  Y_\vartheta d\omega  -\im q'  Y_\omega
d\vartheta + ( Y_j d\overline{z}_{j}-Y_{\overline{j}}dz_j)+
\langle \sigma _1 \sigma _3 Y_f,f' \quad \rangle .
\end{aligned}\end{equation}
For $a_1$ given by \eqref{eq:a1}, and for $  \Gamma  =
i_Y\widetilde{\Omega}$, we have
\begin{equation} \label{eq:linAlg1}  \begin{aligned}     \Gamma
_\omega =& a_1  Y_\vartheta +  \langle \sigma _1\sigma _3\xi _j,
\partial _\omega R \rangle Y_j-  \langle  \sigma _3\xi _j,
\partial _\omega R \rangle Y_{\overline{j}} \\& +  \langle
Y_f,\sigma _3\sigma _1P_c(\omega )  \partial _\omega R \rangle ;\\
-  \Gamma _\vartheta   =&  a_1  Y_\omega -\im  \,   \langle \sigma
_1\sigma _3\xi _j, \sigma _3 R \rangle Y_j+\im \,   \langle \sigma
_3\xi _j, \sigma _3 R \rangle Y_{\overline{j}}
\\& - \im  \,   \langle Y_f,\sigma _3\sigma _1P_c(\omega )  \sigma _3 R \rangle ;
\\    -  \Gamma
_j =&    \langle \sigma _1\sigma _3\xi _j,
\partial _\omega R \rangle Y_\omega +\im \,    \langle \sigma
_1\sigma _3\xi _j, \sigma _3 R \rangle Y_{\vartheta}
   ; \\      \Gamma_{\overline{j}} =&
\langle  \sigma _3\xi _j, \partial _\omega R \rangle Y_\omega +\im
\,    \langle \sigma _3\xi _j, \sigma _3 R \rangle Y_{\vartheta}
  ;
\\     \sigma _3 \sigma _1\Gamma _{f}= &
 ( P_c(\omega _0)  P_c(\omega  )-1 ) Y_{f} \\& +  Y_\omega
 P_c(\omega  _0 )P_c(\omega  )  \partial _\omega R
     +\im \, Y_{\vartheta}  P_c(\omega _0 )
 P_c(\omega  )  \sigma _3 R
.
\end{aligned}\end{equation}
In particular, for $  \gamma  = i_{Y^{t}} {\Omega}_t=i_{Y^{t}}
{\Omega}_0+t\, i_{Y^{t}} \widetilde{{\Omega} } $ we have

\begin{equation} \label{eq:linAlg2}  \begin{aligned}    \gamma
_\omega =&   (\im q'+t a_1)   ({Y}^t)_\vartheta + t  \langle \sigma
_1\sigma _3\xi _j,
\partial _\omega R \rangle  ({Y}^t)_j- t  \langle  \sigma _3\xi _j,
\partial _\omega R \rangle  ({Y}^t)_{\overline{j}} \\& + t \langle
 ({Y}^t)_f,\sigma _3\sigma _1P_c(\omega _0) P_c(\omega  )\partial _\omega R
\rangle ;\\   -  \gamma _\vartheta   =&  (\im q'+t a_1)
 ({Y}^t)_\omega -\im \, t \,   \langle \sigma _1\sigma _3\xi _j,
\sigma _3 R \rangle  ({Y}^t)_j+\im \, t \, \langle \sigma _3\xi _j,
\sigma _3 R \rangle  ({Y}^t)_{\overline{j}}
\\& -\im \,
t \,  \langle  ({Y}^t)_f,\sigma _3\sigma _1P_c(\omega _0) P_c(\omega
 )\sigma _3 R \rangle ;
\\    -  \gamma
_j =&   ({Y}^t) _{\overline{j}}  +t\langle \sigma _1\sigma _3\xi _j,
\partial _\omega R \rangle  ({Y}^t)_\omega +\im \,
t \,    \langle \sigma _1\sigma _3\xi _j, \sigma _3 R \rangle
 ({Y}^t)_{\vartheta}
   ; \\      \gamma_{\overline{j}} =&   ({Y}^t)_{ {j}}+t
\langle  \sigma _3\xi _j, \partial _\omega R \rangle  ({Y}^t)_\omega
+\im \, t \, \,    \langle \sigma _3\xi _j, \sigma _3 R \rangle
 ({Y}^t)_{\vartheta}
  ;
\\    \sigma _3 \sigma _1\gamma _{f}= &  ({Y}^t)_{f}+
t
 ( P_c(\omega _0 )P_c(\omega  )   -1 )      ({Y}^t)_{f} +\\& +t \,
  ({Y}^t)_\omega    \,
 P_c(\omega _0 )P_c(\omega  )  \partial _\omega R
     +t  \, \im \,    ({Y}^t)_{\vartheta}\,  P_c(\omega _0 )
 P_c(\omega  )  \sigma _3 R
    \,
     .
\end{aligned}\end{equation}
\end{lemma}
\proof \eqref{eq:linAlg0} is trivial. \eqref{eq:linAlg2} follows
immediately from \eqref{eq:linAlg0}--\eqref{eq:linAlg1}. In the
following formulas we denote $P_c=P_c(\omega )$, $P_c^0=P_c(\omega
_0)$  and we sum on repeated indexes. We can split
$\widetilde{\Omega} =\widehat{\Omega}+ \widehat{\Omega}_1$ with, see
\eqref{eq:OmegaComponent2},
\begin{equation}   \begin{aligned} & \widehat{\Omega}_1
=\langle    (P_c^0P_c -1)  f' \, , \sigma _3 \sigma _1f' \, \rangle
,\\& \widehat{\Omega} = a_1 d\vartheta \wedge d\omega +dz_j\wedge (
\langle \sigma _1 \sigma _3 \xi _j, \partial _\omega R \rangle
d\omega + \im \langle \sigma _1 \sigma _3 \xi _j, \sigma _3 R
\rangle d\vartheta ) \\& -d\overline{z}_j\wedge ( \langle   \sigma
_3 \xi _j,
\partial _\omega R \rangle d\omega +  \im \langle   \sigma _3 \xi
_j, \sigma _3  R \rangle d\vartheta ) +\\& \langle P_cP_c^0f' \, ,
\sigma _3 \sigma _1P_c \partial _\omega R \rangle \wedge d\omega
+\im \langle P_cP_c^0f' \, , \sigma _3 \sigma _1P_c  \sigma _3 R
\rangle \wedge d\vartheta .
\end{aligned} \nonumber
\end{equation}
Then
\begin{equation} i_Y\widehat{\Omega}_1 =   \langle \sigma _1 \sigma _3
  (P_c^0P_c -1)Y_f , f' \, \rangle   \nonumber
\end{equation}
and
\begin{equation}\begin{aligned} & i_Y\widehat{\Omega} = \big [
 a_1Y_\vartheta +Y_j \langle \sigma _1
  \sigma _3 \xi _j, \partial _\omega R \rangle
 - Y_{\overline{j}} \langle
  \sigma _3 \xi _j, \partial _\omega R \rangle
  +\langle   Y_f
  , \sigma _3 \sigma _1  P_c  \partial _\omega R \rangle \big ]
  d\omega +\\& \big [
 -a_1Y_\omega +\im Y_j \langle \sigma _1
  \sigma _3 \xi _j, \sigma _3 R \rangle
 - \im Y_{\overline{j}} \langle
  \sigma _3 \xi _j, \sigma _3 R \rangle  +\im \langle   Y_f
  , \sigma _3 \sigma _1  P_c  \sigma _3 R   \rangle \big ]
  d\vartheta \\& -(
\langle \sigma _1 \sigma _3 \xi _j, \partial _\omega R \rangle
Y_\omega  + \im \langle \sigma _1 \sigma _3 \xi _j, \sigma _3  R
\rangle Y_ \vartheta ) dz_j\\& +( \langle   \sigma _3 \xi _j,
\partial _\omega R \rangle Y_\omega  + \im \langle   \sigma
_3 \xi _j, \sigma _3  R \rangle Y_ \vartheta ) d\overline{z}_j
\\& -\langle f'\, ,   Y_\omega \sigma _3 \sigma _1 P_c^0P_c
  \partial _\omega R + \im Y_\vartheta \sigma _3 \sigma _1 P_c^0P_c
 \sigma _3 R\rangle .
 \end{aligned} \nonumber
\end{equation}
\qed

\begin{remark} \label{rem:correction} If we choose $\gamma =-\im \alpha$
in Lemma \ref{lem:linearAlgebra} with the $\alpha$ of
\eqref{eq:alpha1}, and if $\mathcal{F}_t$  is the flow of $Y^t$,
then  $  ({Y}^t)_\vartheta \neq 0$   is an obstruction to the
fact that, for $0<t\le 1$,   $K\circ \mathcal{F}_t$  is the
hamiltonian of the  sort of semilinear NLS that \eqref{eq:SystPoiss}
is.  Indeed $ ({Y}^t)_f = - t\im  ({Y}^t)_{\vartheta}\,  P_c(\omega _0 )
 P_c(\omega  )  \sigma _3 R + \mathcal{S}(\R ^3, \C ^2).$ Then if we substitute
 $f$ with $f-  \im  ({Y}^1)_{\vartheta}\,  P_c(\omega _0 )
 P_c(\omega  )  \sigma _3 R +\dots  $  in $\langle \mathcal{H}_\omega f, \sigma _3 \sigma _1f\rangle $ we obtain a term of the form $({Y}^1)_{\vartheta} ^2 \langle \mathcal{H}_\omega f, \sigma _3 \sigma _1f\rangle .  $
To avoid terms like this, we want flows defined from fields with  $ ({Y}^t)_\vartheta =0
$. To this effect we add  a correction to $\alpha$.
\end{remark}
We first   consider the hamiltonian fields of $\vartheta$ and
$\omega$.
\begin{lemma}
  \label{lem:HamThetaOmega} Consider the vectorfield
   $X^t_\vartheta $ (resp. $X^t_\omega $) defined by
   $i_{X^t_\vartheta }\Omega _t=-\im d\vartheta$
   (resp. $i_{X^t_\omega }\Omega _t=-\im d\omega$). Then we have  (here $P_c=P_c(\mathcal{H}_\omega )$
   and $P_c^0=P_c(\mathcal{H}_{\omega _0 })$):
\begin{equation}\label{hamiltonians1}\begin{aligned}& X^t_\vartheta
= (X^t_\vartheta ) _\omega \big [ \frac{\partial}{\partial\omega}  -
t
 \langle
\sigma _3 \xi _j , \partial _\omega R\rangle
\frac{\partial}{\partial z_j} - t
 \langle \sigma _1
\sigma _3 \xi _j , \partial _\omega R\rangle
\frac{\partial}{\partial \overline{z}_j} \\& -tP_c^0(1 +t
P_c-tP_c^0
 )^{-1} P_c^0 P_c\partial _\omega R \big ] ,
\\& X^t_\omega
= (X^t_\omega ) _\vartheta \big [ \frac{\partial}{\partial\vartheta}
- \im t
 \langle
 \xi _j ,   R\rangle
\frac{\partial}{\partial z_j} +\im  t
 \langle
\sigma _1 \xi _j ,   R\rangle \frac{\partial}{\partial
\overline{z}_j} \\& -\im tP_c^0(1 +t  P_c-tP_c^0
 )^{-1} P_c^0
P_c\sigma _3 R \big ] ,
\end{aligned}
\end{equation}
where, for the $a_1$  of \eqref{eq:a1}, we have

\begin{equation}\label{hamiltonians2}\begin{aligned}&
 (X^t_\vartheta ) _\omega =\frac{ \im}{ \im q'+ta_1+ ta_2}=-
 (X^t_\omega ) _\vartheta
\end{aligned}
\end{equation}

   \begin{equation}\label{a2}\begin{aligned}  a_2:=&  \im t
   \langle
\sigma _3 \xi _j , \partial _\omega  R\rangle \langle \sigma _1\xi
_j ,
  R\rangle -\im t \langle \sigma _1 \sigma _3 \xi _j ,
 \partial _\omega  R\rangle \langle    \xi _j ,   R\rangle +\\&
+\im t  \langle P_c^0(1 +t P_c-tP_c^0
 )^{-1} P_c^0 P_c\partial _\omega R,
\sigma _3 \sigma _1 P_c \sigma _3 R\rangle .
\end{aligned}
\end{equation}

\end{lemma}

\proof By \eqref{eq:linAlg2} for $\gamma =-\im \,d\vartheta$,
$X^t_\vartheta $ satisfies
\begin{equation}\label{HamTheta1}\begin{aligned}
& (X^t_\vartheta ) _\vartheta =0; \\ &  \im  = (\im q' + ta_1)
(X^t_\vartheta ) _\omega -\im t \langle \sigma _1 \sigma _3 \xi _j ,
\sigma _3 R\rangle (X^t_\vartheta ) _j+\\& + \im t \langle \sigma _3
\xi _j , \sigma _3 R\rangle (X^t_\vartheta ) _{\overline{j}}-\im t
\langle (X^t_\vartheta ) _{f},  \sigma _3 \sigma _1   P_c \sigma _3 R\rangle ; \\&
(X^t_\vartheta ) _{f} =t (1-P_c^0 P_c )(X^t_\vartheta ) _{f}  -   t
(X^t_\vartheta ) _\omega P_c^0P_c
\partial _\omega R ;\\&
   (X^t_\vartheta ) _{\overline{j}} = -t (X^t_\vartheta ) _\omega
 \langle \sigma _1
\sigma _3 \xi _j , \partial _\omega R\rangle ; \,
  (X^t_\vartheta ) _{ {j}}=- t (X^t_\vartheta ) _\omega
 \langle
\sigma _3 \xi _j , \partial _\omega R\rangle .
\end{aligned}
\end{equation}
This yields \eqref{hamiltonians1} for $X^t_\vartheta $ and the first
equality in \eqref{hamiltonians2}. By \eqref{eq:linAlg2} for $\gamma
=-\im \,d\omega$, $X^t_\omega $ satisfies
\begin{equation}\label{HamOmega1}\begin{aligned}
& (X^t_\omega ) _\omega =0; \\ & -\im \, -\im \, q' (X^t_\omega )
_\vartheta = ta_1 (X^t_\omega ) _\vartheta +  t \langle \sigma _1
\sigma _3 \xi _j , \partial _\omega R\rangle (X^t_\omega ) _j-\\& -
t \langle \sigma _1 \sigma _3 \xi _j , \partial _\omega  R\rangle
(X^t_\omega ) _{\overline{j}}+  t \langle (X^t_\omega ) _{f}, \sigma
_3 \sigma _1 P_c \partial _\omega  R\rangle ;
\\&  (X^t_\omega ) _{f} = t   (1-P_c^0P_c)
 (X^t_\omega ) _{f}-  \im \,  t (X^t_\omega ) _\vartheta P_c^0P_c  \sigma _3 R ;\\&
   (X^t_\omega ) _{\overline{j}} = -\im \,
   t (X^t_\omega ) _\vartheta
 \langle \sigma _1
\sigma _3 \xi _j , \sigma _3  R\rangle ;\quad
  (X^t_\omega ) _{ {j}}=- \im \, t (X^t_\omega ) _\vartheta
 \langle
\sigma _3 \xi _j , \sigma _3 R\rangle .
\end{aligned}
\end{equation}
This yields the rest of
\eqref{hamiltonians1}--\eqref{hamiltonians2}.   \qed

The following lemma is an immediate consequence of the formulas in Lemma \ref{lem:HamThetaOmega} and of \eqref{eq:bounda1}.
\begin{lemma}
  \label{lem:HamBounds}
For any $ (K', S' ,K , S)$ we have
\begin{equation}\label{HamTheta2}\begin{aligned} &
|1-   (X^t_\vartheta ) _\omega  \, q'|\lesssim \| R \| _{H^{-K',-S'}}^2
\\& |(X^t_\vartheta ) _j| +|(X^t_\vartheta )
_{\overline{j}}| + \| (X^t_\vartheta ) _{f}  \| _{H^{ K , S
}}\lesssim \| R \| _{H^{-K',-S'}} .
 \end{aligned}
\end{equation}
and
\begin{equation}\label{HamOmega2}\begin{aligned} &
|1+   (X^t_\omega ) _\vartheta \,  q'|\lesssim \| R \| _{H^{-K',-S'}}^2\,
,
\\& |(X^t_\omega ) _j| +|(X^t_\omega )
_{\overline{j}}| + \| (X^t_\omega ) _{f}  \| _{H^{ -K',-S'
}}\lesssim \| R \| _{H^{-K',-S'}} . \end{aligned}
\end{equation}
\end{lemma}
Set $H_c^{K,S}(\omega )=P_c(\omega )H ^{K,S}$ and denote
\begin{equation}\label{eq:PhaseSpace} \widetilde{\Ph}^{K,S}=\mathbb{C }^m\times H_c^{K,S}(\omega
_0)\, , \quad \Ph^{K,S}=\mathbb{R}^2 \times \widetilde{\Ph}^{K,S}
\end{equation}
with elements $(\vartheta , \omega , z, f)\in \Ph^{K,S}$ and $(  z,
f)\in \widetilde{\Ph}^{K,S}$.

\begin{lemma}
  \label{lem:flow Htheta} We consider  $\forall$ $t\in [0,1]$
the hamiltonian field $X^t_\vartheta $ and
   the flow
\begin{equation}\label{FlowTheta1}
\frac{d}{ds}\Phi _s(t,U)=X^t_\vartheta  (\Phi _s(t,U))\, , \, \Phi
_0(t,U)=U.\end{equation} \begin{itemize}
\item[(1)] For any $(K', S')$ there is a $s_0>0$
and a neighborhood $\mathcal{U}$ of $\mathbb{R}\times \{  (\omega
_0,0,0)\}$ in $\Ph^{-K',-S'}$ such that the map $ (s,t,U)\to \Phi
_s(t,U)$ is smooth

\begin{equation}\label{FlowTheta2} (-s_0,s_0) \times
 [0,1]\times \left (\mathcal{U}\cap \{ \omega =\omega _0\}
 \right )\to   \Ph^ {-K',-S'} .
\end{equation}

\item[(2)] $\mathcal{U}$ can be chosen so that for any $t\in [0,1]$
there  is  another neighborhood $\mathcal{V}_t$ of $\mathbb{R}\times
\{ (\omega _0,0,0)\}$ in $\Ph^{-K',-S'}$ s.t.  the above map
establishes a diffeomorphism

\begin{equation}\label{FlowTheta5} (-s_0,s_0) \times
 \left ( \mathcal{U}\cap \{ \omega =\omega _0\}
 \right )\to   \mathcal{V}_t .
\end{equation}

 \item[(3)]
  $ f(\Phi _s(t,U))
-f (U )=G (t, s , z, f )  $ is a smooth map for all $(K,S)$

\begin{equation}(-s_0,s_0) \times
 [0,1]\times \left ( \mathcal{U}\cap \{ \omega =\omega _0\}
 \right )\to   H^{K,S} \nonumber
\end{equation}
with $ \| G (t, s , z, f ) \| _{H^{K,S}} \le C |s| (|z| +\| f\|
_{H^{-K',-S'}} ).$

\end{itemize}
\end{lemma}\proof
Claims (1)--(2) follow by Lemma \ref{lem:HamThetaOmega} which
implies $X^t_\vartheta \in C^\infty (\mathcal{U}, \Ph ^{K,S})$ for
all $(K,S)$. Let $\zeta $ be any coordinate $z_j$ or $f$.
 Then, for $\zeta $ a scalar coordinate, we have
\begin{equation}\label{FlowTheta31} \begin{aligned} &
 |\zeta (\Phi _s(t,U))- \zeta (U)|\le
 \int _{-s}^s |(X^t_\vartheta  )_{\zeta} (\Phi _{s'}(t,U))|  ds'
 \\& \le C  |s|  \sup _{|s'|\le s} (|z(\Phi _{s'}(t,U))|
 +\| f(\Phi _{s'}(t,U))\|
_{H^{-K',-S'}} ) . \end{aligned}
  \end{equation}
For $\zeta =f$   we have
\begin{equation}\label{FlowTheta32} \begin{aligned} &
 \|f (\Phi _s(t,U))- f (U)\|_{H^{K,S}}\le
 \int _{-s}^s \|(X^t_\vartheta  )_{f} (\Phi _{s'}(t,U))\| _{H^{K,S}}  ds'
   \le  \text{rhs\eqref{FlowTheta31}.}  \end{aligned}\nonumber
  \end{equation}
The above two formulas imply the following, which
yields  claim (3),
\begin{equation}\label{FlowTheta33} \begin{aligned} &
 \|f (\Phi _s(t,U))- f (U)\|_{H^{K,S}}
   \le     C |s| (|z| +\| f\|
_{H^{-K',-S'}} ) , \\&  |z (\Phi _s(t,U))- z (U)|
   \le     C |s| (|z| +\| f\|
_{H^{-K',-S'}} ) .\end{aligned}
  \end{equation}
\qed
\begin{lemma}
  \label{lem:correction alpha} We consider a scalar function
  $F(t,U)$ defined as follows:
  \begin{equation}\label{defCorrect} F (t,\Phi _s(t,U))= \im \,
\int _0^s\alpha   _{\Phi _{s'}(t,U)}\left ( X^t_\vartheta  (\Phi
_{s'}(t,U))\right ) ds'\, ,  \text{ where $\omega (U)=\omega _0$ .}
\end{equation}
We have $F
  \in C^{ \infty} ( [0,1]\times \mathcal{U}, \mathbb{R})$ for
   a neighborhood $\mathcal{U}$ of $\mathbb{R}\times \{  (\omega
_0,0,0)\}$ in $\Ph^{-K',-S'}$. We have
\begin{equation} \label{estCorrect}
|F (t,U)| \le C (K',S')  |\omega -\omega _0|\, \left (  |z|+ \| f \|
_{H^{-K',-S'}}\right )^2.
   \end{equation}
We have (exterior differentiation only in $U$)
\begin{equation}\label{Vectorfield21}
(\alpha +\im \, d F )(X^t_\vartheta ) =0.
\end{equation}

\end{lemma}
\proof    $F$  is smooth by \eqref{eq:alpha2} and Lemma
\ref{lem:flow Htheta}.  \eqref{Vectorfield21}  follows by   \eqref{FlowTheta1}
and by  \begin{equation}\alpha   _{ U}\left (X^t_\vartheta (U)\right ) +
  \im  \frac{d}{ds} _{|s=0}   F (t,\Phi _s(t,U)) =0 . \end{equation}    By
\eqref{eq:alpha2} and \eqref{HamTheta2}  we have
\begin{equation}\label{BoundCorrect1} \begin{aligned} &
|\alpha (X^t_\vartheta )|\le | \alpha ^\omega | \, | (X^t_\vartheta
) _\omega | + | \langle \alpha ^f, (X^t_\vartheta ) _f \rangle |
\lesssim  \left ( |z|+\| f \|  _{H^{-K',-S'}} \right ) ^2.
\end{aligned}
\end{equation}
Then \eqref{estCorrect} follows by $|s|\approx |\omega (\Phi
_s(t,U))-\omega _0|.$ \qed

\begin{lemma}
  \label{lem:vectorfield}  Denote by $\mathcal{X}^t$    the vector field which solves

\begin{equation}\label{Vectorfield1}
 i_{\mathcal{X}^t} \Omega _t=-\alpha - \im \, d   F (t) .
   \end{equation}
Then the following properties hold.

 \begin{itemize}
\item[(1)]  There is a neighborhood $\mathcal{U}$ of
$ \mathbb{{R}}\times \{ (\omega _0, 0,0)  \} $ in $\Ph ^{1,0}$ such
that $ \mathcal{X}^t ( U)
  \in C^{ \infty} ( [0,1]\times \mathcal{U}, \Ph ^{1,0})$.

  \item[(2)]  We have
$
 (\mathcal{X}^t)_\vartheta  \equiv 0.
   $

   \item[(3)] For constants $C(K,S,K',S')$ we have

    \begin{equation}\label{Vectorfield3}\begin{aligned} &
\left | (\mathcal{X}^t)_\omega + \frac{\|f\| _2^2 }{2q'(\omega )}
\right | \lesssim (|z|+\| f \| _{H^{-K',-S'}})^2;  \\&
|(\mathcal{X}^t  )_
 {j} |  +|(\mathcal{X}^t  )_
 {\overline{j }} | +\|
  (\mathcal{X}^t  )_{f} \| _{H^{ K , S }}
  \lesssim  (|z|+\| f \| _{H^{-K',-S'}}) \times  \\&
  \times
(|\omega -\omega _0| +|z|+\| f \| _{H^{-K',-S'}} +\|f\| _{L^2}^2).
   \end{aligned}\end{equation}

\item[(4)]    We have
\begin{equation} \label{eq:L_X^t} L _{\mathcal{X}^t}
\frac{\partial}{\partial \vartheta} :=\left [ \mathcal{X}^t,
\frac{\partial}{\partial \vartheta} \right ] =0.
\end{equation}
\end{itemize}
\end{lemma}

\proof Claim (1) follows from the regularity properties of $\alpha$,
$F$ and $\Omega _t$ and from equations \eqref{Vectorfield5} and
\eqref{Vectorfield7} below. \eqref{Vectorfield21} implies (2) by
\begin{equation}  \im (\mathcal{X}^t)_\vartheta =
\im d\vartheta (\mathcal{X}^t)=- i _{X^t_\vartheta}\Omega _t
(\mathcal{X}^t )= i _{\mathcal{X}^t}\Omega _t (X^t_\vartheta
)=-(\alpha +\im \, dF )(X^t_\vartheta ) =0. \nonumber
\end{equation}
 We have $\im (\mathcal{X}^t)_\omega =\im d\omega
  (\mathcal{X}^t)=- i _{X^t_\omega}\Omega _t
(\mathcal{X}^t )$, so by \eqref{Vectorfield1} and \eqref{hamiltonians1} we get
\begin{equation}\label{Vectorfield5} \begin{aligned}
 & \im (\mathcal{X}^t)_\omega = i _{\mathcal{X}^t}\Omega
_t (X^t_\omega )= -(X^t_\omega )_{\vartheta} \big [\alpha ^\vartheta
  +t
\partial _j F \,  \langle
 \xi _j ,   R\rangle  -t
\partial _{\overline{j}}F   \langle
\sigma _1 \xi _j ,   R\rangle \\ & +t \langle \nabla _{f}F + \im
\alpha ^{f} , P_c^0(1 +t  P_c-tP_c^0
 )^{-1} P_c^0
P_c\sigma _3 R \rangle \big ].\end{aligned}
\end{equation}
Then by \eqref{eq:alpha2}, \eqref{hamiltonians2}, \eqref{eq:a1} and
\eqref{a2}, we get  the first inequality in \eqref{Vectorfield3}:
\begin{equation}\label{Vectorfield51} \begin{aligned}
&\left | (\mathcal{X}^t)_\omega +\frac{\| f\| _2^2  }{2q'(\omega )}
    \right |
  \le C  \left ( |z| +\|  f\| _{H^{-K',-S'}}  \right ) ^2.\end{aligned}
\end{equation}
By \eqref{eq:linAlg2} we have the following equations
\begin{equation}\label{Vectorfield7} \begin{aligned}
      \im \, \partial _j F  & =(\mathcal{X}^t  )_
 {\overline{j }} +
 t\langle \sigma _1 \sigma _3 \xi _j, \partial _\omega R \rangle
(\mathcal{X}^t)_\omega \, \\  - \im \,
\partial _{\overline{j }}F   & =(\mathcal{X}^t  )_{j} +
 t\langle  \sigma _3 \xi _j, \partial _\omega R \rangle
(\mathcal{X}^t)_\omega \, \\   \sigma _3\sigma _1(\alpha ^{f}+\im \,
\nabla _{f} F ) & =-(\mathcal{X}^t  )_{f} -
 t   (P_c ^0  P_c-1)(\mathcal{X}^t  )_{f}
  \\& - t (\mathcal{X}^t)_\omega  P_c^0P_c   \partial _\omega
 R   .\end{aligned}
\end{equation}
Formulas \eqref{Vectorfield7} imply

\begin{equation}\label{Vectorfield8} \begin{aligned}
  &     |(\mathcal{X}^t_\omega )_
 {\overline{j }} |\le    |\partial _{j} F | +C
 \left ( |z| +\|  f\| _{H^{-K',-S'}}  \right )
|(\mathcal{X}^t)_\omega |\\ & |(\mathcal{X}^t_\omega )_
 {j} |\le    |\partial _{\overline{j}}
F | +C\left ( |z| +\|  f\| _{H^{-K',-S'}}  \right )
|(\mathcal{X}^t)_\omega | \\& \|
  (\mathcal{X}^t_\omega )_{f} \| _{H^{ K , S }} \le
\|
  \alpha ^{f} \| _{H^{ K , S }}+ \|
 \nabla
_{f} F  \| _{H^{ K , S }} +C\left ( |z| +\|  f\| _{H^{-K',-S'}}
\right ) |(\mathcal{X}^t)_\omega |   \end{aligned}\nonumber
\end{equation}
which with \eqref{Vectorfield51}, \eqref{eq:alpha2} and Lemma
\eqref{estCorrect} imply \eqref{Vectorfield3}.  \eqref{eq:L_X^t} is a consequence of the following equalities, which we will justify below:
\begin{equation}\label{lieder}0= L_{\frac{\partial}{\partial \vartheta}} \left (
i_{\mathcal{X}^t} \Omega _t\right ) = i_{[\frac{\partial}{\partial
\vartheta},\mathcal{X}^t]} \Omega _t+i_{\mathcal{X}^t}
L_{\frac{\partial}{\partial \vartheta}}\Omega _t
=i_{[\frac{\partial}{\partial \vartheta},\mathcal{X}^t]} \Omega _t .
\end{equation}
The first equality is a consequence of \eqref{Vectorfield1} and
  $L_{\frac{\partial}{\partial \vartheta}} \left (\alpha + \im
d F \right ) =0  $. The latter is a consequence of $ L_{\frac{\partial}{\partial \vartheta}} \alpha =0$  and $\frac{\partial}{\partial \vartheta}F=0$. Notice that
 $\frac{\partial}{\partial \vartheta}F=0$ can be proved observing that
   \eqref{Vectorfield21}, \eqref{eq:alpha2} and  Lemma \ref{lem:HamThetaOmega} imply $ X^t _\vartheta \frac{\partial}{\partial \vartheta}F=0$ and that
   on $\omega =\omega _0$ we have  $\frac{\partial}{\partial \vartheta}F=0$.
   $ L_{\frac{\partial}{\partial \vartheta}} \alpha =0$  is a consequence of
     the Cartan "magic"      formula $L_X \gamma = (i_Xd+ di_X)\gamma $,
    of  the definition \eqref{eq:alpha1}
     and of  following equalities:

  \begin{equation}\label{eq:magic} \begin{aligned} &
  L_{\frac{\partial}{\partial \vartheta}} \beta =  d i_{\frac{\partial}{\partial \vartheta}} \beta +   i_{\frac{\partial}{\partial \vartheta}} d\beta = -\frac{\im }{2}d\langle \sigma _1U,U\rangle + \im \langle \sigma _1U,\quad \rangle =0;\\&
   L_{\frac{\partial}{\partial \vartheta}} \beta _0=  d i_{\frac{\partial}{\partial \vartheta}} \beta _0+   i_{\frac{\partial}{\partial \vartheta}} d\beta _0=-\im dq -\im  i_{\frac{\partial}{\partial \vartheta}} (dq\wedge d\vartheta )=-\im dq +\im  dq\wedge  i_{\frac{\partial}{\partial \vartheta}} d\vartheta =0; \\&
  L_{\frac{\partial}{\partial \vartheta}}d\psi = d   i_{\frac{\partial}{\partial \vartheta}}d\psi = \frac{1}{2}d \frac{\partial}{\partial \vartheta}\langle \sigma _3\Phi
  , R\rangle =0.
  \end{aligned}
\end{equation}
   The second equality in \eqref{lieder} follows  by the product rule for the Lie
derivative. Finally, the third  equality in \eqref{lieder} follows  by
$L_{\frac{\partial}{\partial \vartheta}}\Omega _t= (1-t)L_{\frac{\partial}{\partial \vartheta}}\Omega _0 + tL_{\frac{\partial}{\partial \vartheta}}\Omega  =0$, consequence of
$ L_{\frac{\partial}{\partial \vartheta}}\Omega   =0$ (resp. $ L_{\frac{\partial}{\partial \vartheta}}\Omega  _0 =0$), in turn consequence of the first (resp. second) line in \eqref{eq:magic} and of the identity $L_X d\gamma =   dL_X \gamma $.
\qed

We have:
\begin{lemma}
  \label{lem:flow1}    Consider the vectorfield $\mathcal{X}^t$
in Lemma \ref{lem:correction alpha} and denote by $\mathcal{F}_t(U)$
the corresponding flow. Then the flow $\mathcal{F}_t(U)$ for $U$
near $ e^{\im \sigma _3 \vartheta}\Phi _{\omega _0}$ is defined for
all $t\in [0,1]$. We have $\vartheta \circ \mathcal{F}_1=\vartheta
$. We have for $\ell =j,\overline{j}$,
\begin{equation}\label{flow5}\begin{aligned} &
q\left ( \omega   (\mathcal{F}_1 (U))\right ) = q\left (\omega (U)
\right ) -\frac{\| f \|_2^2 }{2}  + \mathcal{E}_{\omega} (U)
\\& z_\ell (\mathcal{F}_1 (U))=z_\ell  (U)+ \mathcal{E}_{\ell }
(U)\\& f (\mathcal{F}_1 (U))= f (U)   + \mathcal{E}_{f}  (U)
\end{aligned}
\end{equation}

with
\begin{eqnarray}
  \label{flow6}  & |\mathcal{E}_{
\omega  }(U)| \lesssim   (|\omega -\omega _0| +|z|+\| f \|
_{H^{-K',-S'}} )^2,
\\& \label{flow7}
|\mathcal{E}_{  \ell }(U)| +\| \mathcal{E}_f(U)\| _{H^{ K , S
}}\lesssim (|\omega -\omega _0| +|z|+\| f \| _{H^{-K',-S'}}+\| f \|
^2_{L^{2}} ) \\& \times (|\omega -\omega _0| +|z|+\| f \|
_{H^{-K',-S'}} ) .\nonumber
\end{eqnarray}
For each $\zeta =\omega , z_\ell , f$ we have
\begin{equation}\label{flow71}\begin{aligned} & \mathcal{E}_\zeta(U)
=\mathcal{E}_\zeta ( \| f\| _{L^2}^2, \omega , z, f)
\end{aligned}
\end{equation}
with, for a neighborhood $\mathcal{U}^{-K',-S'}$ of
$\mathbb{R}\times \{ (\omega _0,0,0)\}$ in $\Ph^{-K',-S'}$ and for
some  fixed $a_0>0$
\begin{equation}\label{flow8}\begin{aligned} & \mathcal{E}_\zeta
( \varrho , \omega , z, f)\in C^\infty ( (-a_0,a_0)\times
\mathcal{U}^{-K',-S'}, \mathbb{C})
\end{aligned}
\end{equation}
for  $\zeta =\omega , z_\ell  $ and with
\begin{equation}\label{flow9}\begin{aligned} & \mathcal{E}_f
( \varrho , \omega , z, f)\in C^\infty ( (-a_0,a_0)\times
\mathcal{U}^{-K',-S'}, H ^{ K , S }).
\end{aligned}
\end{equation}

\end{lemma}
\proof We add a new variable $\varrho$. We define a new field by
\begin{equation}\label{Vectorfield52} \begin{aligned}
 & \im (Y^t)_\omega =   -(X^t_\omega )_{\vartheta} \big [
 \alpha ^\vartheta +\im \frac{\| f\|_2^2-\rho }{2}
  +t
\partial _j F \,  \langle
 \xi _j ,   R\rangle  -t
\partial _{\overline{j}}F   \langle
\sigma _1 \xi _j ,   R\rangle \\ & +t \langle \nabla _{f}F + \im
\alpha ^{f} , P_c^0(1 +t  P_c-tP_c^0
 )^{-1} P_c^0
P_c\sigma _3 R \rangle \big ], \end{aligned}
\end{equation}
  by \begin{equation}\label{Vectorfield71} \begin{aligned}
      \im \, \partial _j F  & =(Y^t  )_
 {\overline{j }} +
 t\langle \sigma _1 \sigma _3 \xi _j, \partial _\omega R \rangle
(Y^t)_\omega \, \\  - \im \,
\partial _{\overline{j }}F   & =(Y^t  )_{j} +
 t\langle  \sigma _3 \xi _j, \partial _\omega R \rangle
(Y^t)_\omega \, \\   \sigma _3\sigma _1(\alpha ^{f}+
\im \, \nabla
_{f} F ) & =(Y^t  )_{f} +
 t   (P_c ^0  P_c-1)(Y^t  )_{f}
  \\& - t (Y^t)_\omega  P_c^0P_c   \partial _\omega
 R   .\end{aligned}
\end{equation}
and by $Y^t_\rho =2\langle (Y^t)_f,\sigma _1 f\rangle $. Then   $Y^t=Y^t(
\omega , \rho , z, f)$ defines a new flow $\mathcal{G}_t (\rho , U)
$, which reduces to $\mathcal{F}_t ( U) $ in the invariant manifold
defined by $\rho = \| f\| _2^2.$ Notice that by $\rho (t)= \rho (0)+\int _0^t
Y_{\rho}^{s}ds$ it is easy to conclude $\rho  (\mathcal{G}_1 (\rho ,U) )=  ( U) +O(\text{rhs\eqref{flow6}})$.   Using \eqref{HamOmega2}, \eqref{eq:alpha2} and
\eqref{Vectorfield52} it is then easy to get
\begin{equation} \begin{aligned} &q(\omega (t))= q(\omega (0))+\int _0^t q'(\omega (s))Y_{\omega}^{s}ds =q(\omega (0))-\int _0^t \frac{\rho(  s ) }{2} ds +O(\text{rhs\eqref{flow6}}). \end{aligned}\nonumber
\end{equation}
By standard arguments, see for example the proof of Lemma 4.3 \cite{bambusicuccagna},  we get
\begin{equation}\label{flow51}\begin{aligned} &
q\left ( \omega   (\mathcal{G}_1 (\rho ,U))\right ) = q\left (
\omega (U)\right )   -\frac{\rho }{ 2}  + \mathcal{E}_{\omega} (\rho
,U)\, ,
\\& z_\ell (\mathcal{G}_1 (\rho ,U))=z_\ell  (U)+ \mathcal{E}_{\ell
} (\rho ,U)\, , \\& f (\mathcal{G}_1 (\rho ,U))= f (U)   +
\mathcal{E}_{f} (\rho ,U) \, ,
\end{aligned}
\end{equation}
with $\mathcal{E}_{\zeta }  (\rho ,U)$ satisfying \eqref{flow8} for
$\zeta =\omega  , z_\ell$ and \eqref{flow9} for $\zeta =f$. We have
$\mathcal{E}_{\zeta }  (  U)=\mathcal{E}_{\zeta }  (\| f\| _2 ,U)$
  satisfying \eqref{flow6} for
$\zeta =\omega   $ and \eqref{flow7} for $\zeta =z_\ell ,f$. \qed

 We have:
\begin{lemma}
  \label{lem:flow2}
Consider  the flow $\mathcal{F}_t$ of Lemma \ref{lem:flow1}. Then we
have
\begin{equation} \label{eq:Darboux} \mathcal{F}_t^*\Omega _t=\Omega
_0 .
\end{equation}
  We have
\begin{equation} \label{eq:QcircF}
Q\circ \mathcal{F}_1=q.
\end{equation}
If $\chi$ is a function with $\partial _\vartheta \chi \equiv 0$,
then $\partial _\vartheta (\chi \circ \mathcal{F}_t) \equiv 0$.
\end{lemma}
\proof \eqref{eq:Darboux} is Darboux Theorem, see
\eqref{eq:dartheorem}.
  Let $ \mathcal{G}_t= (\mathcal{F}_t) ^{-1}$.  Then
$\mathcal{G}_t^*\Omega _0=\Omega _t$. We have
$\mathcal{G}_t^*X^0_{q(\omega)} =X^t_{q(\omega) \circ
\mathcal{G}_t}$ by
\begin{equation}\begin{aligned} & i _{\mathcal{G}_t^*X^0_{q(\omega)}}
 \Omega _t= i _{\mathcal{G}_t^*X^0_q(\omega)}
\mathcal{G}_t^*\Omega _0= \mathcal{G}_t^*i _{X^0_q(\omega)}\Omega
_0=-\im d(q(\omega) \circ \mathcal{G}_t ) = i _{X^t_{q(\omega) \circ
\mathcal{G}_t}}\Omega _t.
\end{aligned} \nonumber
\end{equation}
 Then by  $\left [ \mathcal{X}^t,
\frac{\partial}{\partial \vartheta} \right ] =0$  for all $t$
\begin{equation}\frac{d}{dt}X^t_{q(\omega) \circ \mathcal{G}_t}
 = \frac{d}{dt}
\mathcal{G}_t^*X^0_{q(\omega)} = - \frac{d}{dt}
\mathcal{G}_t^*\frac{\partial}{\partial \vartheta} =-
\mathcal{G}_t^*\left [ \mathcal{X}^{1-t}, \frac{\partial}{\partial
\vartheta} \right ] =0.  \nonumber
\end{equation}
So   $X_{{q(\omega)}\circ \mathcal{G}_1}^1= X_{q(\omega)}^0$. Since
by \eqref{eq:Ham.VecFieldQ} and \eqref{eq:HamVect0q} this implies $d
(q\circ \mathcal{G}_1) =dQ$ and since there are points with $q\circ
\mathcal{G}_1 (U)= Q (U) $,  we obtain \eqref{eq:QcircF}. Finally,
the last statement of Lemma \ref{lem:flow2} follows by
\eqref{eq:L_X^t} and by

\begin{equation} \frac{\partial}{\partial \vartheta}\mathcal{F}_t^*
\chi = \left (\mathcal{F}_t^*\frac{\partial}{\partial
\vartheta}\right )\left (\mathcal{F}_t^* \chi \right )=
\mathcal{F}_t^* \left ( \frac{\partial}{\partial \vartheta}\chi
\right ) . \nonumber
\end{equation}
 \qed

\section{Reformulation of \eqref{eq:SystK} in the new coordinates}
\label{section:reformulation}
 We set
\begin{equation}    \label{eq:newH}  \begin{aligned} &
H=K\circ \mathcal{F}_1.
\end{aligned}
\end{equation}
In the new coordinates \eqref{eq:SystK}  becomes
\begin{equation} \label{eq:SystK1} \begin{aligned} &
q' \dot \omega  =  \frac{\partial H}{\partial \vartheta}\equiv 0 \,
, \quad    q'  \dot \vartheta  =
 -\frac{\partial H}{\partial \omega}
  \end{aligned}
\end{equation}
and

\begin{equation} \label{eq:SystK2} \begin{aligned} &
 \im \dot z_j  =   \frac{\partial H}{\partial   \overline{z}_j  }
\, , \quad
 \im {\dot {\overline{z}}}_j  = - \frac{\partial H}{\partial    {z}_j  }
\\&   \im \dot f=    \sigma _3    \sigma _1 \nabla _f  H. \end{aligned}
\end{equation}
Recall that we are solving the initial value problem \eqref{NLS} and
that we have chosen $\omega _0$ with $q(\omega _0)=\| u_0\|
_{L^2_x}^{2}.$ Correspondingly it is enough to focus on
\eqref{eq:SystK2} with $\omega =\omega _0$.  For system
\eqref{eq:SystK2} we prove :
\begin{theorem}\label{theorem-1.2}
  Then there exist $\varepsilon  >0$
and $C>0$ such that for   $ |z(0)|+\| f (0) \| _{H^1 }\le \epsilon
<\varepsilon  $  the corresponding  solution of \eqref{eq:SystK2} is
globally defined and there are    $f_\pm  \in H^1$ with $\| f_\pm\|
_{H^1 }\le C \epsilon $   such that
\begin{equation}\label{scattering1}\lim _{t\to \pm \infty } \| e^{ \im  \vartheta (t) \sigma _3     } f(t) -
  e^{ \im t \sigma _3  \Delta    } f_\pm \| _{H^1 }=0
\end{equation}
  where  $\vartheta (t)$  is the variable associated to $U^T(t)=(u(t), \overline{u}(t))$ in \eqref{eq:anzatz} and \eqref{eq:coordinate}.
  We also have
\begin{equation}\label{decay}\lim _{t\to   \infty } z(t)=0.
\end{equation}
In particular,  it is possible to write
$R(t,x)=A(t,x)+\widetilde{f}(t,x)$ with $|A(t,x)|\le C_N(t) \langle
x \rangle ^{-N}$ for any $N$, with $\lim _{t\to \infty }C_N(t)=0$
and such that for any admissible  pair $(r,p)$, i.e.
\eqref{admissiblepair}, we have
\begin{equation}\label{Strichartz1}
 \|  \widetilde{f} \| _{L^r_t(  \mathbb{R},
W^{ 1,p}_x)}\le
 C\epsilon .
\end{equation}
\end{theorem}
By Lemma \ref{lem:flow1},     Theorem \ref{theorem-1.2} implies
Theorem \ref{theorem-1.1}. Indeed, if we denote $(\omega , z', f')$ the initial coordinates, and $(\omega _0 , z , f )$ the coordinates in \eqref{eq:SystK2},
we have $z'=z+O(|z|+ \| f \| _{L _x ^{2,-2}} )$ and  $f'=f+O(|z|+ \| f \| _{L _x ^{2,-2}} )$. The two error terms $O$ converge to 0 as $t\to \infty$. Hence the asymptotic behavior of $(  z', f')$ and of $(  z , f )$ is the same. We also have
$
q\left ( \omega   (t) \right ) = q\left (\omega _0 \right )
-\frac{\| f(t) \|_2^2 }{2}    + O( |z(t)|+\| f(t) \| _{L^{2,-2 }_x} )
$  which implies, say at $+\infty$
\begin{equation}   \begin{aligned}
 &   \lim _{t \to +\infty }
q\left ( \omega   (t) \right ) = \lim _{t \to +\infty  }
  \left ( q\left (\omega _0 \right ) -\frac{\|  e^{ {\rm i}t
\sigma _3   \Delta } f_+ \|_2^2 }{2}\right ) = q\left (\omega _0 \right )
-\frac{\|   f_+ \|_2^2 }{2}=q(\omega _+)
\end{aligned} \nonumber
\end{equation}
for somewhere  $\omega _+$ is the unique element near $\omega _0$ for which the last inequality holds. So $ \lim _{t
\to +\infty } \omega (t) = \omega _+.$

In the rest of the paper we focus on Theorem \ref{theorem-1.2}. The
main idea  is that \eqref{eq:SystK2} is basically like the system
considered in \cite{bambusicuccagna}. Therefore Theorem
\ref{theorem-1.2} follows by the Birkhoff normal forms argument of
\cite{bambusicuccagna}, supplemented  with the various dispersive
estimates in \cite{cuccagnamizumachi}.

\subsection{Taylor expansions}
\label{subsec:Taylor expansions}

Consider $U=  e^{\im \sigma _3\vartheta }  (\Phi _\omega + R)$ as in
\eqref{eq:anzatz}. Decompose $R$ as in \eqref{eq:decomp2}. Set
$u=\varphi +  u_c$ with $^t (P_c(\omega )f)=(u_c,\overline{u_c})$.
We have

\begin{equation}\label{ExpEP0}\begin{aligned}
 &   B (|u  |^2)=B \left (   |u_c |^2
   \right ) +\int _0^1 \left [   \frac{\partial }{\partial u}
   B (|u|^2 )_{|u=u_c+t\varphi} \varphi + \frac{\partial }
   {\partial \overline{u}}
   B (|u|^2 )_{|u=u_c+t\varphi} \overline{\varphi} \right ]
dt\\& = B \left (   |u_c |^2
   \right ) +\int _0^1  dt  \sum _{i+j\le 4}  \frac{1}{i!j!}
 \partial _{u}^{i+1} \partial _{\overline{u}}^j
 B \left ( | u | ^2 \right )_{|u= t\varphi} u_c^i\overline{u_c}^j \varphi +\\&
 \int _0^1  dt  \sum _{i+j\le 4}  \frac{1}{i!j!}
 \partial _{u}^{i } \partial _{\overline{u}}^{j+1}
 B \left ( | u | ^2 \right )_{|u= t\varphi}  u_c^i\overline{u_c}^j
 \overline{\varphi} +\\&
 5\int _{[0,1]^2}  dt ds  (1-s)^4\sum _{i+j= 5} \frac{1}{i!j!}
 \partial _{u}^{i+1 } \partial _{\overline{u}}^{j }
 B \left ( | u| ^2 \right ) _{|u= t\varphi +su_c}  u_c^i\overline{u_c}^j \varphi
+\\&
 5\int _{[0,1]^2}  dt ds  (1-s)^4\sum _{i+j= 5} \frac{1}{i!j!}
 \partial _{u}^{i } \partial _{\overline{u}}^{j+1}
 B \left ( | u| ^2 \right ) _{|u= t\varphi +su_c} u_c^i\overline{u_c}^j
  \overline{\varphi} .
\end{aligned}\end{equation}

\begin{lemma}
\label{lem:K} The  following statements  hold.
\begin{equation}  \label{eq:ExpK} \begin{aligned} & K =d(\omega )-\omega
\| u_0\| _2^2+K_2+K_P  \\& K_2=\sum _j\lambda _j (\omega  ) |z_j|^2+
\frac{1}{2} \langle \sigma _3 \mathcal{H}_{\omega  } f, \sigma _1
f\rangle \\& K_P = \sum _{|\mu +\nu |= 3} \langle a_{\mu \nu
}(\omega ,z ) , 1 \rangle z^\mu \overline{z}^\nu +\sum _{|\mu +\nu
|= 2} z^\mu \overline{z}^\nu \langle  G_{\mu \nu }(\omega ,z
),\sigma _3\sigma _1P_c(\omega )f\rangle   \\&   + \sum _{d=2}^4
\langle B_{d } ( \omega , z  ), (P_c(\omega )f)^{\otimes d} \rangle
+   \langle B_6 (\omega , f) , 1 \rangle   +\int _{\mathbb{R}^3}
B_5(x,\omega, z , f(x) )  f^{\otimes 5}(x) dx,
\end{aligned}\nonumber
\end{equation}
for $ B_6 (x, \omega , f)=B \left ( \frac{|P_c(\omega )f(x)|^2
  }{2}\right ),$
where we have what follows.
\begin{itemize}
\item[(1)] $a_{\mu \nu
}( \cdot , \omega ,z ) \in C^\infty ( \mathrm{U},
H^{K,S}_x(\mathbb{R}^3,\mathbb{C})) $ for any pair $(K,S)$ and a
small neighborhood $\mathrm{U}$ of $(\omega _0,0)$ in
$\mathcal{O}\times \mathbb{C}^m$.

\item[(2)] $G_{\mu \nu
}( \cdot , \omega ,z ) \in C^\infty ( \mathrm{U},
H^{K,S}_x(\mathbb{R}^3,\mathbb{C}^2)) $, for $\mathrm{U}$ like in
(1), possibly smaller.

\item[(3)] $B_{d
}( \cdot , \omega ,z ) \in C^\infty ( \mathrm{U},
H^{K,S}_x(\mathbb{R}^3, B   (
 (\mathbb{C}^2)^{\otimes d},\mathbb{C} ))) $, for $2\le d \le 4$
for $\mathrm{U}$   possibly smaller.

\item[(4)] Let $^t\eta = (\zeta , \overline{\zeta}) $ for
$ \zeta \in \C$. Then for
  $B_5(\cdot ,\omega  , z , \eta  )$   we have
\begin{equation} \label{5power2}\begin{aligned} &\text{for any $l$ ,
}  \| \nabla _{ \omega ,z,\overline{z} ,\zeta,\overline{\zeta}   }
^lB_5( \omega ,z,\eta  ) \| _{H^{K,S}_x (\mathbb{R}^3,   B   (
 (\mathbb{C}^2)^{\otimes 5},\mathbb{C} )} \le C_l .
 \end{aligned}\nonumber \end{equation}
\item[(5)] We have $a_{\mu \nu }=\overline{a}_{\nu \mu   }$, $G_{\mu \nu }
=-\sigma _1\overline{G}_{\nu \mu   } $.

 \end{itemize}
\end{lemma}
\proof The expansion for $K$ is a consequence of well know
cancelations. (1)--(4) follow from \eqref{ExpEP0}  and elementary
calculus.   (5) follows from the fact that $K(U)$ is real valued for
$\overline{U}=\sigma _1U$.\qed

We set $\delta _j$ be for $j\in \{ 1,...m \}$ the multi index
$\delta _j=( \delta _{1j}, ..., \delta _{mj}).$ Let $\lambda
_j^0=\lambda _j(\omega _0)$ and $\lambda ^0 = (\lambda _1^0, \cdots,
\lambda _m^0)$.

\begin{lemma}
  \label{lem:ExpH} Let $H=K\circ \mathcal{F} _1$. Then, at $
 e^{i\sigma _3 \vartheta}  \Phi _{\omega _0}$
 we have the
expansion
\begin{equation}  \label{eq:ExpH1} \begin{aligned} & H =d(\omega _0)
-\omega _0\| u_0\| _2^2+ \psi (\|f\| _2^2) +H_2 ^{(1)}+\resto ^{(1)}
\end{aligned}
\end{equation}   for $\omega =\omega _0$, where the following holds.
\begin{itemize}
\item[(1)]
 We have for $r=1$

\begin{equation}  \label{eq:ExpH2} H_2 ^{(r)}=
\sum _{\substack{ |\mu +\nu |=2\\
\lambda ^0\cdot (\mu -\nu )=0}}
 a_{\mu \nu}^{(r)}( \| f\| _2^2 )  z^\mu
\overline{z}^\nu + \frac{1}{2} \langle \sigma _3 \mathcal{H}_{\omega
_0} f, \sigma _1 f\rangle .
\end{equation}

\item[(2)] We have $\resto ^{(1)}=\widetilde{\resto ^{(1)}} +
 \widetilde{\resto ^{(2)}}  $, with $\widetilde{\resto ^{(1)}}=$
\begin{equation}  \label{eq:ExpH2resto} \begin{aligned} &
 =\sum _{\substack{ |\mu +\nu |=2\\
\lambda ^0\cdot (\mu -\nu )\neq 0  }} a_{\mu \nu }^{(1)}(\| f\| _2^2
)z^\mu \overline{z}^\nu +\sum _{|\mu +\nu |  = 1} z^\mu
\overline{z}^\nu \langle \sigma _1 \sigma _3G_{\mu \nu }(\| f\| _2^2
),f\rangle ,\\& \widetilde{\resto ^{(2)}}= \sum _{|\mu +\nu |  = 3 }
z^\mu \overline{z}^\nu \int _{\mathbb{R}^3}a_{\mu \nu }
(x,z,f,f(x),\| f\| _2^2 ) dx\\&  + \sum _{|\mu +\nu | =2  }z^\mu
\overline{z}^\nu \int _{\mathbb{R}^3} \left [ \sigma _1 \sigma
_3G_{\mu \nu } (x,z,f,f(x), \| f\| _2^2 )\right ]^*f(x) dx  \\&
+\sum _{j=2}^5 \resto ^{(1)} _j   +
    \int _{\mathbb{R}^3}B (|f(x)|^2/2 ) dx + \widehat{\resto}
    ^{(1)} _2(z,f,  \| f\| _2^2 )  \\& \text{with }
\resto ^{(1)} _j    =\int _{\mathbb{R}^3} F_j(x,z ,f, f(x),\| f\|
_2^2) f^{\otimes j}(x) dx .
\end{aligned}
\end{equation}
\item[(3)]   We have  $F_2(x,0 ,0, 0,0)=0 $.

\item[(4)]   $\psi (s)$ is smooth with $\psi (0)=\psi ' (0)=0$.

\item[(5)] At $\| f\| _2=0$ with $r=1$
\begin{equation}  \label{eq:ExpHcoeff1} \begin{aligned} &
a_{\mu \nu }^{(r)}( 0 ) =0 \text{ for $|\mu +\nu | = 2$  with $(\mu
, \nu )\neq (\delta _j, \delta _j)$ for all $j$,} \\& a_{\delta _j
\delta _j }^{(r)}( 0 ) =\lambda _j (\omega _0)  , \text{ where
$\delta _j=( \delta _{1j}, ..., \delta _{mj}),$}
\\& G_{\mu \nu }(  0 ) =0 \text{ for $|\mu +\nu | = 1$ }
\end{aligned}
\end{equation}
These $a_{\mu \nu }^{(r)}( \varrho )$ and $G_{\mu \nu }( x,\varrho
)$ are smooth in all variables with $G_{\mu \nu }( \cdot ,\varrho )
\in C^\infty ( \mathbb{R}, H^{K,S} _x(\mathbb{R}^3,\mathbb{C}^2))$
for all $(K,S)$.

\item[(6)] We have for all indexes and for $r=1$
\begin{equation}  \label{eq:ExpHcoeff2} \begin{aligned} & a_{\mu \nu }^{(r)} =
\overline{a}_{\nu \mu   }^{(r)}\, , \quad a_{\mu \nu }
=\overline{a}_{\nu \mu   }\, , \quad  G_{\mu \nu } =-\sigma
_1\overline{G}_{\nu \mu } .
\end{aligned}
\end{equation}

\item[(7)] Let $^t\eta = (\zeta , \overline{\zeta}) $ for
$ \zeta \in \C$. For all $(K,S, K',S')$ there is a neighborhood $
\mathcal{U}^{-K',-S'}$ of $  \{ (  0, 0) \}$ in
$\widetilde{\Ph}^{-K',-S'}$, see \eqref{eq:PhaseSpace}, such that
we have,  for
 $a_{\mu \nu } (x, z,f,\eta , \varrho )$ with
 $(z,f,\zeta , \varrho )\in \mathcal{U}^{-K',-S'}\times \C \times \R$
 \begin{equation}\label{eq:coeff a2} \| \nabla _{
z,\overline{z},\zeta,\overline{\zeta},f,\varrho} ^l a_{\mu \nu } \|
_{ H^{K,S}_x(\mathbb{R}^3,\mathbb{C})} \le C_l \text{ for all $l$}.
\end{equation}
\item[(8)]  Possibly restricting $\mathcal{U}^{-K',-S'}$, we have also, for
 $G_{\mu \nu } (x , z,f,g ,
\varrho ) $,
\begin{equation}\label{eq:coeff G2}
\| \nabla _{ z,\overline{z},\zeta,\overline{\zeta},f,\varrho} ^l
G_{\mu \nu } \| _{ H^{K,S}_x(\mathbb{R}^3,\mathbb{C}^2)} \le C_l
  \text{ for all $l$}.\end{equation}

\item[(9)] Restricting $\mathcal{U}^{-K',-S'}$  further, we have also,
for $F_j (x ,z,f,g,\varrho )$,
  \begin{equation} \begin{aligned} &
\| \nabla _{ z,\overline{z},\zeta,\overline{\zeta},f,\varrho} ^l F_j
\| _{ H^{K,S}_x(\mathbb{R}^3,B   (
 (\mathbb{C}^2)^{\otimes j},\mathbb{C} ))} \le C_l
  \text{ for all $l$}
.
 \end{aligned}\nonumber \end{equation}

\item[(10)]  Restricting $\mathcal{U}^{-K',-S'}$
 further, we have $\widehat{\resto} ^{(1)}
 _2(z,f,  \varrho ) \in C^\infty ( \mathcal{U}^{-K',-S'}\times
  \R , \R)
   $ with  \begin{equation} \begin{aligned} & |\widehat{\resto} ^{(1)}
 _2(z,f,  \varrho ) | \le C (|z|+|\varrho |+ \| f \| _{ H^{-K',-S'}})
\| f \| _{ H^{-K',-S'}}^2.
 \end{aligned}\nonumber \end{equation}
\end{itemize}
\end{lemma}
\proof By $\mathcal{F}_1(\Phi _{\omega _0})= \Phi _{\omega _0}$,
  $K'(\Phi _{\omega _0})=0$ and   $\|
\mathcal{F}_1(U) -U\| _{\Ph ^{K,S}}\lesssim \| R \| _{L^2}^2$ we
conclude  $H'(\Phi _{\omega _0})=0$ and  $H''(\Phi _{\omega _0})
=K''(\Phi _{\omega _0})$. In particular, this yields the formula for
$H_2^{(1)}$   for $\| f\| _2=0$. The other terms are obtained by
substituting in \eqref{eq:ExpK} the formulas  \eqref{flow5}. By
$\langle \sigma _3   f, \sigma _1 f\rangle =0$ we have $\langle
\sigma _3 \mathcal{H}_{\omega _0 +\delta \omega   } f, \sigma _1
f\rangle =\langle \sigma _3 \mathcal{H}_{\omega  _0   } f, \sigma _1
f\rangle  + \widetilde{F}_2$ where $\widetilde{F}_2$ can be absorbed
in $j=2$ in \eqref{eq:ExpH2resto}. $\psi (\| f\| _2^2)$ arises from
$d(\omega \circ \mathcal{F}_1)-\omega \circ \mathcal{F}_1\| u_0\|
_2^2$.  Other terms  coming from the latter end up in
\eqref{eq:ExpH2resto}: in particular there are no monomials $\| f\|
_2^j z^\mu \overline{z}^\nu \langle G, f\rangle ^i$ with $|\mu +\nu
|+i=1$, because of \eqref{flow6} (applied for $\omega =\omega _0$).

 \qed

 \section{Canonical transformations}
\label{sec:Canonical} Our goal in this section is to prove the
following result.

\begin{theorem}
\label{th:main}   For any integer $r\ge 2$   there are a
neighborhood $ \mathcal{U}^{1,0}$ of $  \{ (  0, 0) \}$ in
$\widetilde{\Ph}^{1,0}$, see \eqref{eq:PhaseSpace}, and a smooth
canonical transformation $\Tr_r: \mathcal{U}^{1,0}\to
\widetilde{\Ph}^{1 ,0}$ s.t.
\begin{equation}
\label{eq:bir1} H^{(r)}:=H\circ \Tr_r=d(\omega _0) -\omega _0\|
u_0\| _2^2+ \psi (\|f\| _2^2)+H_2^{(r)}+Z^{(r)}+\resto^{(r)},
\end{equation}
     where:
\begin{itemize}
\item[(i)] $H_2^{(r)}=H_2^{(2)}$ for $r\ge 2$, is of the
form \eqref{eq:ExpH2}
where  $a_{\mu \nu}^{(r)}(\| f\| _2)$  satisfy
\eqref{eq:ExpHcoeff1}--\eqref{eq:ExpHcoeff2};

\item[(ii)]$Z^{(r)}$ is in normal form, in the sense of Definition \ref{def:normal form} below, with monomials of
degree $\le r$  whose coefficients satisfy \eqref{eq:ExpHcoeff2};
\item[(iii)] the transformation   $\Tr _r$ is of the form (see below)
 \eqref{lie.11.a}--
  \eqref{lie.11.b} and satisfies \eqref{lie.11.f}--\eqref{lie.11.c}
   for $M_0=1$;
\item[(iv)] we have $\resto^{(r)} = \sum _{d=0}^6\resto^{(r)}_d$
with the following properties:
\begin{itemize}
\item[(iv.0)] for all $(K,S, K',S')$ there is a neighborhood
$  \mathcal{U}^{-K',-S'}$ of $   \{ (  0, 0) \}$ in
$\widetilde{\Ph}^{-K',-S'}$
 such that

\begin{equation} \resto^{(r)}_0=
\sum _{|\mu +\nu |  = r+1 } z^\mu \overline{z}^\nu \int
_{\mathbb{R}^3}a_{\mu \nu }^{(r)}(x,z,f,f(x),\| f\| _2^2 ) dx
\nonumber
\end{equation}
and for   $a_{\mu \nu }^{(r)}(z,f,\eta , \varrho )$   with $^t\eta =
(\zeta , \overline{\zeta} )$, $\zeta \in \C$ we have for $(z,f)\in
\mathcal{U}^{-K',-S'}$ and $|\varrho |\le 1$

\begin{equation}\label{eq:coeff a} \|\nabla _{ z,\overline{z},
\zeta,\overline{\zeta},f,\varrho} ^l a_{\mu \nu }^{(r)}(\cdot
,z,f,\eta , \varrho )\| _{H^{K,S}(\R ^3, \C  )} \le C_l \text{ for
all $l$};
\end{equation}

\item[(iv.1)]  possibly taking $\mathcal{U}^{-K',-S'}$
  smaller, we have
\begin{equation} \resto^{(r)}_1=
\sum _{|\mu +\nu |  = r  }z^\mu \overline{z}^\nu \int
_{\mathbb{R}^3} \left [ \sigma _1 \sigma _3G_{\mu \nu
}^{(r)}(x,z,f,f(x), \| f\| _2^2 )\right ]^*f(x) dx  \nonumber
\end{equation}

\begin{equation}\label{eq:coeff G}\text{with }  \|\nabla _{ z,\overline{z},
\zeta,\overline{\zeta},f,\varrho} ^l G_{\mu \nu }^{(r)}(\cdot
,z,f,\eta , \varrho )\| _{H^{K,S}(\R ^3, \C ^2)} \le C_l \text{ for
all $l$};
\end{equation}

\item[(iv.2--5)]  possibly  taking $\mathcal{U}^{-K',-S'}$
  smaller, we have for $2\le d \le 5$,

\begin{equation} \resto^{(r)}_d=
 \int
_{\mathbb{R}^3} F_d^{(r)}(x, z ,f,f(x),\| f\| _2^2)  f^{\otimes
d}(x) dx +  \widehat{\resto}^{(r)}_d,\nonumber
\end{equation} with   for
any $l$

\begin{equation}\label{eq:coeff F}   \|\nabla _{ z,\overline{z},
\zeta,\overline{\zeta},f,\varrho} ^l F_d^{(r)} (\cdot ,z,f,\eta ,
\varrho )\| _{H^{K,S} (\mathbb{R}^3,   B   (
 (\mathbb{C}^2)^{\otimes d},\mathbb{C} )} \le C_l,
\end{equation}
  with $F_2^{(r)}(x, 0 ,0,0,0)=0$  and
  with $\widetilde{\resto}^{(r)}_d (z ,f, \| f\| _2^2)$ s.t.

\begin{equation}\label{eq:Rhat}\begin{aligned} &
\widehat{\resto}^{(r)}_d (z ,f, \varrho )\in C^\infty
(\mathcal{U}^{-K',-S'}\times \R ,\R ), \\& |
\widehat{\resto}^{(r)}_d (z ,f, \varrho )|\le C \| f \|
_{H^{-K',-S'}}^d , \\&    | \widehat{\resto}^{(r)}_2 (z ,f, \varrho
)|  \le C (|z|+|\varrho |+ \| f \| _{ H^{-K',-S'}}) \| f \| _{
H^{-K',-S'}}^2;
\end{aligned}\end{equation}

\item[(iv.6)] $ \resto^{(r)}_6=
 \int _{\mathbb{R}^3}  B ( | f (x)| ^2/2) dx$.

\end{itemize}
\end{itemize}
\end{theorem}

We develop the proof in the following subsections. The basic  ideas
are classical. However we need to develop a number of tools, along
the lines of \cite{bambusicuccagna}. The situation here   is more
complicated than in \cite{bambusicuccagna} because of the dependence of the coefficients on $\| f\| _2$.

\subsection{Lie transform}
\label{subsec:LieTransf} We   consider   functions
\begin{equation}
\label{chi.1}    \chi  =   \sum_{|\mu +\nu |=M_0 +1}b_{\mu   \nu
}(\| f\| _2^2) z^{\mu} \overline{z}^{\nu}   + \sum_{|\mu +\nu |=M_0
}z^{\mu} \overline{z}^{\nu}
 \langle   \sigma _1\sigma _3B_{\mu   \nu
}(\| f\| _2^2)
  , f \rangle   \end{equation}
where    $ b_{\mu   \nu }(\varrho )\in C^{\infty
}(\mathbb{R}_\varrho , \mathbb{C})$ and $ B_{\mu   \nu }(x,\varrho)
\in C^{\infty }(\mathbb{R} , P_c(\omega _0)H^{k,s}_x(\mathbb{R}^3,
\mathbb{C}^2))$ for all $k$ and $s$. Assume
\begin{equation} \label{chi.11} b_{\mu   \nu } =\overline{b}_{\nu
\mu   }\text{ and }\sigma _1 B_{\mu   \nu }=-\overline{B}_{ \nu \mu}
\text{ for all indexes}. \end{equation} We set for  $K>0$ and $S>0$
fixed and large  set
\begin{equation} \label{chi.12} \|  \chi \| = \|  \chi (\|  f\|
_2^2)\| =
\sum |b_{\mu   \nu }(\|  f\| _2^2)|+\sum \|B_{\mu   \nu }(\|  f\|
_2^2)\| _{H^{ K , S }}. \end{equation} Denote by $\phi ^t$ the flow
of the Hamiltonian vector field $X_{\chi}$ ( from now on with
respect to $\Omega _0$ and only in $(z,f)$). The {\it Lie transform}
$\phi =
 \phi^t\big|_{t=1}$ is defined in a sufficiently small neighborhood
 of the origin and is a canonical transformation.

\begin{lemma}\label{lie_trans}
Consider the $\chi$ in \eqref{chi.1} and its Lie transform $\phi$.
Set $(z',f')=\phi (z,f)$.   Then there are $ \mathcal{G}(z,f,
\varrho)$, $\Gamma (z,f,\varrho)$, $\Gamma _0 (z,f,\rho)$ and
$\Gamma _1 (z,f,\rho)$ with the following properties.

\begin{itemize}
\item[(1)]  $\Gamma \in C^\infty
(  \U ^{-K',-S'}, \mathbb{C}^{ m}) $, $\Gamma _0, \Gamma _1\in
C^\infty ( \U ^{-K',-S'}, \mathbb{R}) $, with $\U ^{-K',-S'}\subset
\mathbb{C}^{ m} \times H^{-K',-S'}_c(\omega _0)\times \mathbb{R}$ an
appropriately small neighborhood of the origin.
\item[(2)] $  \mathcal{G}\in C^\infty
(\U ^{-K',-S'} ,  H^{K,S}_c(\omega _0) )  $  for any $K,S$.
\item[(3)] The transformation $\phi$ is of the following form:
\begin{eqnarray}
\label{lie.11.a}&    z'  =   z  +
 \Gamma (z,f,\|  f\| _2^2) ,\\
\label{lie.11.b} & f' = e^{ \im \Gamma _0 (z,f,\|  f\|
_2^2)P_c(\omega _0) \sigma _3}f +\mathcal{G}(z,f,\|  f\| _2^2)  .
\end{eqnarray}

\item[(4)] We have
\begin{eqnarray} \label{lie.11.h}  & \| f'\| _2^2 =
\| f   \| _2^2  +
 \Gamma _1(z,f,\|  f\| _2^2),
\\& \nonumber
\left | \Gamma _1(z,f,\|  f\| _2^2) \right | \le \\&
\label{lie.11.f}  C |z|^{M_0 -1} ( |z|^{M_0+2 } + |z|^{2 }\| f \|
_{H^{-K',-S'}} +  \| f \| _{H^{-K',-S'}} ^3)    .
\end{eqnarray}

\item[(5)] There are constants   $c_{K',S'}$ and $c_{K, S,K',S'}$
such that \begin{equation}\label{lie.11.c}\begin{aligned}&
 |\Gamma (z,f,\|  f\| _2^2)| \leq c_{K',S'}  (\| \chi \| +\text{\eqref{lie.11.f}})  |z | ^{M_0-1}
 ( |z|+ \norma{f }
_{H^{-K',-S'}} ), \\ &  \|\mathcal{G}(z,f,\|  f\|
_2^2 )\|_{H^{K,S}}\leq c_{K, S,K',S'}   (\| \chi \| +\text{\eqref{lie.11.f}})   |z | ^{{M_0} } , \\
&  |\Gamma _{0}(z,f,\|  f\| _2^2)| \leq c_{K',S'} |z |^{M_0-1  } ( |z | +\| f   \|
_{H^{ -K',-S'}} )^2.
 \end{aligned}\end{equation}

\item[(6)] We have
\begin{equation} \label{eq:prime f1}
\begin{aligned} & e^{ \im \Gamma _0P_c(\omega _0) \sigma _3    }
=e^{ \im \Gamma _0  \sigma _3    } +T (\Gamma _0),
\end{aligned}
\end{equation}
where $T (r)\in C^\infty (\mathbb{R}, B ( H ^{-K',-S'},H ^{ K ,  S
}) ) $ for all $(K,S,K',S')$, with norm $ \| T (r) \| _{B ( H
^{-K',-S'},H ^{ K ,  S })} \le C (K,S,K',S') |r|.$ More
specifically, the range of $T(r)$ is $R(T(r))\subseteq L^2_d
(\mathcal{H})+L^2_d (\mathcal{H}^*),  $ $L^2_d$ defined two lines after
\eqref{eq:spectraldecomp}.

\end{itemize}
\end{lemma}

\proof   Claim (6) can be proved independently of the properties of $\Gamma_0$. Recall that
$P_c(\omega )=
1-P_d(\omega )$, see below \eqref{eq:spectraldecomp}, with $P_d(\omega )$ smoothing and of finite rank. Exploiting $\sigma _3P_d(\omega
)=P_d^*(\omega )\sigma _3$ it is elementary to prove
\begin{equation} \label{eq:prime f}
\begin{aligned} & e^{ \im \Gamma _0P_c(\omega _0) \sigma _3    }
=e^{ \im \Gamma _0  \sigma _3    } +T(\Gamma _0) \text{ with
}T(\Gamma _0)= - \im \sin \left (\Gamma _0\right ) P_d(\omega _0)
  \sigma _3+\\&
 +    \sum _{n=2}^\infty
\frac{ (\im \Gamma _0 )^n}{n!} \sum _{j=1}^{\left [\frac{n}{2}\right
]}\left(
\begin{matrix}
\left [\frac{n}{2}\right ]\\ j
\end{matrix} \right) K^j ( P_c(\omega _0) \sigma _3)
^{\varepsilon (n)},
\end{aligned}
\end{equation}
with $K=P_d(\omega _0)P_d^*(\omega _0)-P_d(\omega _0)-P_d^*(\omega
_0)$ and $\varepsilon (n) = \frac{1-(-1)^n}{2}$. $T(\Gamma _0)$ has
the properties of Claim (6).

In the sequel we prove Claims (1)--(5).
Set $\varrho = \| f \| _2^2$. For $b_{\mu \nu} '$ and
 $B_{\mu \nu}'$ derivatives with respect to $\varrho$, summing on repeated
 indexes, consider
\begin{equation}\label{eq:gam}\begin{aligned}&
\end{aligned} \gamma (z,f ,\varrho):=2(
   b_{\mu \nu} '(\varrho)z^\mu \overline{z}^\nu+ \langle
\sigma _1\sigma _3 B_{\mu \nu}'(\varrho), f\rangle z^\mu
\overline{z}^\nu ) .
\end{equation}
 For $\sigma _1f=\overline{f}$, then $\gamma
(z,f ,\varrho)\in \mathbb{R}$ by \eqref{chi.11}. We set up the following system:
\begin{equation}\label{auxHamEq1}\begin{aligned}& \im \dot z_j
=\sum _{|\mu + \nu |=M_0+1}\nu _j  \frac{z^\mu \overline{z}^\nu}{\overline{z}_j}b_{\mu
\nu}(\varrho)+\sum _{|\mu + \nu |=M_0 }\nu _j \frac{z^\mu
\overline{z}^\nu}{\overline{z}_j}\langle \sigma _1\sigma _3 B_{\mu
\nu}(\varrho), f\rangle  \\& \im \dot f =\sum _{|\mu + \nu |=M_0 }z^\mu \overline{z}^\nu
B_{\mu \nu}(\varrho) + \gamma (z,f,\varrho )  P_c(\omega _0 )\sigma
_3 f
\\& \dot \varrho  =- 2\im  \langle \sum _{|\mu + \nu |=M_0 } z^\mu \overline{z}^\nu
B_{\mu \nu}(\varrho) + \gamma (z,f,\varrho )  (P_c(\omega _0
)-P_c^*(\omega _0 ))\sigma _3 f, \sigma _1 f \rangle ,
\end{aligned}
\end{equation}
where in the last equation we exploited $\langle \sigma _3 f, \sigma
_1 f\rangle =0.$ By \eqref{chi.11} the flow leaves the set with
$\sigma _1f=\overline{f}$ and $\varrho \in \mathbb{R}$ invariant. In
particular, the set where $\varrho = \| f \| _2^2$ is invariant
under the flow of \eqref{auxHamEq1}. In a neighborhood of 0 the
lifespan of the solutions is larger than
 1.
 \eqref{lie.11.a} can always been written.  For $\gamma$ defined in \eqref{eq:gam}, we have

\begin{equation} \label{eq:Solsystemchi}
\begin{aligned} &
      f (t  ) = e^{- \im   \int _0^t
       \gamma ds P_c(\omega _0)  \sigma _3} f -\sum _{|\mu + \nu |=M_0}
       \im \int _0^t z^\mu \overline{z}^\nu
e^{  \im   \int _s^t
       \gamma ds' P_c(\omega _0) \sigma _3}
B_{\mu \nu} ds . \end{aligned}\nonumber
\end{equation}
  This yields  \eqref{lie.11.b}. We can always write
 \begin{equation} \label{eq:lie.rho1} \begin{aligned} &  \varrho '=\varrho +\Gamma _1
(z,f, \varrho )  .\end{aligned}
\end{equation}
This yields \eqref{lie.11.h}.
  Claims  (1)--(2)   follow  from the regularity of the
flow of \eqref{auxHamEq1} on the initial data. By \eqref{auxHamEq1} we get

 \begin{equation} \label{eq:closing1} \begin{aligned} &
  |\varrho (t) -  \varrho | \le C \sup _{0\le t'\le t}   |z(t')|^{M_0 -1}   ( |z(t')|^{M_0+2 }  +\\& + |z(t')|^{2 } \| f(t') \|
_{H^{-K',-S'}} +  \| f(t') \| _{H^{-K',-S'}} ^3)  .\end{aligned}
\end{equation}
Similarly we have
 \begin{equation} \label{eq:closing2} \begin{aligned} &
  |z (t) -  z | \le C \sup _{0\le t'\le t}   |z(t')|^{M_0 - 1}  \|   \chi (\varrho (t')) \|    ( |z(t')| +  \| f(t') \|
_{H^{-K',-S'}}  )  , \end{aligned}
\end{equation}
\begin{equation} \label{eq:closing3} \begin{aligned} &
\| \int _0^t z^\mu \overline{z}^\nu
e^{  \im   \int _s^t
       \gamma ds' P_c(\omega _0) \sigma _3}
B_{\mu \nu} ds\| _{H^{K,S}} \le   C \sup _{0\le t'\le t}|z(t')|^{M_0  } \|   \chi (\varrho (t')) \| .
 \end{aligned}
\end{equation}
Then $|z(t)|\approx |z|+ \| f  \|
_{H^{-K',-S'}}$ with in particular $|z(t)|\approx |z|$ when $M_0>1$.
By Claim (6) and by the fact that the exponent $\Gamma _0(z,f,\varrho )$ in \eqref{lie.11.b} is a uniformly bounded function, we get
$\| f (t) \|
_{H^{ -K',-S'}}\approx |z|+ \| f  \|
_{H^{-K',-S'}}$. Then
\begin{equation} \label{eq:closing4} \begin{aligned} & \left |
\|   \chi (\varrho (t')) \| -  \|   \chi (\varrho  ) \| \right |\le \text{\eqref{lie.11.f}}.
 \end{aligned}
\end{equation}
This implies that   the right hand sides of \eqref{eq:closing1}--\eqref{eq:closing3}  are bounded by the bounds of
$\Gamma _1$, $\Gamma  $ and $\mathcal{G}$ in the statement. This
yields the desired bounds on $\Gamma _1$, $\Gamma  $ and $\mathcal{G}$.
The bound on $\Gamma _0$ follows from
\begin{equation} \label{eq:closing5} \begin{aligned} &  |\int _0^t\gamma (t')  dt'|\le C   \sup _{0\le t'\le t} |z(t')|^{M_0  } ( |z(t')| +\| f (t') \|
_{H^{ -K',-S'}} ) \\& \le C_1 |z |^{M_0-1  } ( |z | +\| f   \|
_{H^{ -K',-S'}} )^2 .
 \end{aligned}
\end{equation}

 \qed

\subsection{Normal form}
\label{subsec:Normal form}
  Recall  the notation $\lambda
_j^0=\lambda _j(\omega _0)$ and $\delta _j= (\delta _{1j},...,
\delta _{mj})$, see before Lemma \ref{lem:ExpH}. Let $\mathcal{H}=\mathcal{H }_{\omega _{0}}
P_c(\mathcal{H} _{\omega _{0}} )$. For $r\ge 1$, using the
coefficients in \eqref{eq:ExpH2}  of the $H^{(r)}_2$ in Theorem
\ref{th:main}, let
\begin{equation} \label{eq:lambda}\lambda _j^{(r)}=
\lambda _j^{(r)} ( \| f\| _2^2 ) =\lambda _j (\omega _0) +  a
_{\delta _j\delta _j} ^{(r)}(\| f\| _2^2 ), \quad \lambda ^{(r)}=
(\lambda _1^{(r)}, \cdots, \lambda _m^{(r)}).\end{equation}
\begin{definition}
\label{def:normal form} A function $Z(z,f)$ is in normal form if it
is of the form
\begin{equation}
\label{e.12} Z=Z_0+Z_1
\end{equation}
where we have finite sums of the following types:
\begin{equation}
\label{e.12a}Z_1= \sum _{|\lambda ^0 \cdot(\nu-\mu)|>\omega _{0}}
z^\mu \overline{z}^\nu \langle  \sigma _1\sigma _3 G_{\mu \nu}(\| f
\| _2^2 ),f\rangle
\end{equation}
with $G_{\mu \nu}( x,\varrho )\in  C^\infty ( \mathbb{
R}_{\varrho},H_x^{K,S})$ for all $K$, $S$;
\begin{equation}
\label{e.12c}Z_0= \sum _{\lambda  ^0\cdot(\mu-\nu)=0} a_{\mu , \nu}
(\| f \| _2^2)z^\mu \overline{z}^\nu
\end{equation}

and $a_{\mu , \nu} (\varrho  )\in  C^\infty ( \mathbb{ R}_\varrho ,
\mathbb{C})$. We will always assume the symmetries
\eqref{eq:ExpHcoeff2}. \qed\end{definition}

For an $H_2^{(r)}$ as in \eqref{eq:ExpH2}    let
$H_2^{(r)}=D_2^{(r)} +(H_2^{(r)}-D_2^{(r)})$ where
\begin{equation}  \label{eq:Diag} D_2^{(r)}=
\sum _{j=1}^m
 \lambda  _j ^{(r)} ( \| f\| _2^2 )  |z^j|  + \frac{1}{2} \langle \sigma _3 \mathcal{H}_{\omega
_0} f, \sigma _1 f\rangle .
\end{equation}
In the following formulas we set $\lambda _j= \lambda  _j ^{(r)}$,
$\lambda = \lambda
  ^{(r)}$ and $D_2=D_2^{(r)}$.
We recall ($\lambda _j' (\varrho)$ is the derivative in $\varrho$)
that by \eqref{eq:PoissonBracket}, summing on repeated indexes,

 \begin{equation} \label{PoissBra1}
\begin{aligned} &\{ D_2, F  \} :=dD_2 (X_F)= \partial _j D_2
(X_F)_j+
\partial _{\overline{j}}
D_2 (X_F)_{\overline{j}} +\langle \nabla _fD_2, (X_F)_f\rangle \\ &=
-\im \partial _j D_2\partial _{\overline{j}} F+\im \partial
_{\overline{j}} D_2\partial _jF- \im  \langle \nabla _fD_2, \sigma
_3 \sigma _1\nabla _fF \rangle =\\&  \im \lambda _j  {z}_j\partial
_jF - \im \lambda _j \overline{z}_j
 \partial _{\overline{j}} F+\im \langle  \mathcal{H}f,
   \nabla _fF \rangle  +2\im \lambda _j'
(\| f \| _2^2) |z_j|^2 \langle   f , \sigma _3  \nabla _fF \rangle
.
\end{aligned}
\end{equation}
In particular, we have,  for $G=G(x)$, (we use   $\sigma _1\im
\sigma _2=\sigma _3$)
 \begin{equation} \label{PoissBra2} \begin{aligned} & \{
D_2, z^{\mu} \overline{z}^{\nu}   \} =  \im \lambda \cdot (\mu -
\nu) z^{\mu} \overline{z}^{\nu} ,\\& \{ D_2, \langle \sigma _1\sigma
_3 G ,f\rangle \} = -\im \langle f ,\sigma _1\sigma _3 \mathcal{H}G
\rangle -2 \, \im \sum _{j=1}^{m}\lambda ' _j |z_j|^2 \langle \sigma
_1 f,G\rangle ,
\\& \{ D_2, \frac{1}{2} \| f\| _2^2 \} =\{  D_2, \frac{1}{2} \langle
f , \sigma _1 f \rangle \} =\im \langle \mathcal{H}f,
 \sigma _1 f   \rangle = -\im \langle
 \beta ^\prime (\phi ^2  )\phi ^2 \sigma _3 f,
  f   \rangle .
\end{aligned}
\end{equation}
In the sequel we will assume (and prove) that $\| f\| _2$ is small.
We will consider only $|\mu  +\nu |\le 2N +3$. Then, $\lambda
^0\cdot(\mu-\nu)\neq 0$ implies $|\lambda ^0\cdot(\mu-\nu)|\ge c >0$
for some fixed $c$, and so we can assume also $|\lambda
 \cdot(\mu-\nu)|\ge c/2$. Similarly $|\lambda
^0\cdot(\mu-\nu)|<\omega  _{0} $ (resp.   $|\lambda
^0\cdot(\mu-\nu)|>\omega  _{0} $) will be assumed equivalent to
$|\lambda
 \cdot(\mu-\nu)|<\omega  _{0} $ (resp.   $|\lambda
 \cdot(\mu-\nu)|>\omega  _{0} $).

\begin{lemma}
\label{sol.homo} Consider

\begin{eqnarray}
\label{eq.stru.1} K =\sum _{|\mu +\nu |=M_0+1} k_{\mu\nu} (\| f \|
_2^2 ) z^{\mu} \overline{z}^{\nu} + \sum _{|\mu +\nu |=M_0} z^{\mu}
\overline{z}^{\nu}
 \langle  \sigma _1\sigma _3  K_{\mu   \nu
}(\| f \| _2^2)
  , f \rangle .
\end{eqnarray}
 Suppose that all the terms in \eqref{eq.stru.1} are not in normal
 form and that the symmetries
\eqref{eq:ExpHcoeff2} hold.   Consider
\begin{equation}
\label{solhomo} \begin{aligned} &\chi =  \sum_{|\mu +\nu
|=M_0+1}\frac{ k_{\mu\nu}(\| f \| _2^2)}{
  \im \lambda \cdot(\nu-\mu)} z^\mu \overline{z}^\nu \\& -
   \sum_{|\mu +\nu
|=M_0}
   z^\mu \overline{z}^\nu   \langle   \sigma _1\sigma _3
    \frac{1}{   \im  (\lambda \cdot
 (\mu -\nu ) -\mathcal{H}) }  K_{\mu   \nu
}(\| f \| _2^2)
  , f \rangle .\end{aligned}
\end{equation}
Then we have
\begin{equation}
\label{eq:homologicalEq} \left\{\chi , D_2  \right\} =K+L
\end{equation}
with,  summing on repeated indexes,
\begin{equation}
\label{eq:RestohomologicalEq} \begin{aligned} & L=   2
   \frac{k_{\mu\nu}'  }{(\mu -\nu )\cdot \lambda } z^{\mu}
\overline{z}^{\nu}   \langle
 \beta ^\prime (\phi ^2  )\phi ^2 \sigma _3 f,
   f   \rangle \\& +2 \lambda ' _j z^\mu \overline{z}^{\nu} |z_j|^2
   \left \langle \sigma _1 f , \frac{1}{(\mu -\nu )\cdot \lambda
   -\mathcal{H}}K_{\mu \nu } \right \rangle -\\&  2 \lambda ' \cdot(\mu -\nu )  z^\mu \overline{z}^{\nu} |z_j|^2
   \left \langle \sigma _1 f , \frac{1}{\left ( (\mu -\nu )\cdot \lambda
   -\mathcal{H}\right ) ^2}K_{\mu \nu } \right \rangle \langle
 \beta ^\prime (\phi ^2  )\phi ^2 \sigma _3 f,
   f   \rangle
\\& +2z^{\mu}
\overline{z}^{\nu} \left \langle  f ,\sigma _3 \sigma
_1\frac{1}{(\mu -\nu )\cdot \lambda
   -\mathcal{H}}K_{\mu \nu } '\right \rangle
 \langle
 \beta ^\prime (\phi ^2  )\phi ^2 \sigma _3 f,f\rangle   . \end{aligned}
\end{equation}
The coefficients in \eqref{solhomo}  satisfy \eqref{eq:ExpHcoeff2}.
\end{lemma}
\proof The proof follows by the tables \eqref{PoissBra2}, by the
product rule for the derivative and by the symmetry properties of
$\mathcal{H}$.   \qed

We split the proof of Theorem \ref{th:main} in two stages. We first
prove step $r=2$ of Theorem \ref{th:main}. We subsequently prove the
case  $r>2$.

\subsection{Proof of Theorem \ref{th:main}: the step $r=2$}
\label{subsec:step1} At this step,  our goal is to obtain a
hamiltonian similar to $H$, but with  $\widetilde{\resto ^{(1)}} =
0$. We will need to solve a nonlinear homological equation. We
consider a $\chi$ like in \eqref{chi.1} with $M_0=1$ satisfying
\eqref{chi.11}. We write
\begin{equation} \label{eq:ExpH11} \begin{aligned} & H \circ
\phi=d(\omega _0) -\omega _0\| u_0\| _2^2+ \psi (\|f'\| _2^2)
+(H_2^{(1)} +\widetilde{\resto ^{(1)}}   +\widetilde{\resto ^{(2)}})
 \circ \phi ,
\end{aligned}
\end{equation}
for $\phi$ the Lie transform of $\chi$. We  write
\eqref{lie.11.a}--\eqref{lie.11.b} as follows, where we sum on
repeated indexes and $\nabla _f$ does not act on $\|  f\| _2^2$:
\begin{eqnarray}
\label{lie.12.a}&    z'_j  -   z _j =
  \partial _k\Gamma _j(0,0,\|  f\| _2^2)z_k +
 \partial _{\overline{k}}\Gamma _j(0,0,\|  f\| _2^2)\overline{z}_k +\\&
  +
 \langle \nabla _f\Gamma _j(0,0,\|  f\| _2^2), f\rangle + r_j
   , \nonumber\\
\label{lie.12.b} & f' -  e^{ \im \Gamma _0 (z,f,\|  f\|
_2^2)P_c(\omega _0) \sigma _3}f =
  \partial _k\mathcal{G}( 0,0,\|  f\| _2^2)z_k +
 \partial _{\overline{k}}\mathcal{G}(0,0,\|  f\| _2^2)
 \overline{z}_k + r_f . \nonumber
\end{eqnarray}
By \eqref{lie.11.c}--\eqref{lie.11.d} the terms in
rhs\eqref{lie.12.a} satisfy (see \eqref{chi.12} for definition of $
\| \chi \| $)
\begin{equation}\label{lie.13.a}\begin{aligned} &
|\partial _k\Gamma _j |+\cdots \| \partial
_{\overline{k}}\mathcal{G}  \| _{H^{ K , S }} \le C \|  \chi \|
\\&
|r_j|+\| r_f\| _{H^{K,S}}\le C (|z|^2+\|  f\| _{H^{-K',-S'}} ^2).
\end{aligned}
\end{equation}  We   write  the $f'^{\otimes 2}$ in
 \eqref{eq:ExpH2resto}  schematically as
\begin{equation}  \label{eq:fotimes2}\begin{aligned} &
f '^{2}(x)= \sum _{|\mu +\nu |=2}A_{\mu \nu}(x,\|  f\| _2^2)
z^{\mu}\overline{z}^{\nu} +\sum _{|\mu +\nu |=1}
z^{\mu}\overline{z}^{\nu}
      \mathcal{A} _{\mu   \nu
}(\| f\| _2^2)(x)
    f (x) \\& +
    (e^{ \im \Gamma _0   \sigma _3}f + T(\Gamma _0)f)  ^{2}(x) +
    \varphi (x) r_j f(x)+r_f(x) f(x)+ \varphi (x) r_j ^2+r_f^2(x)
\end{aligned}
\end{equation}
where $ \varphi (x)$ represents an exponentially decreasing smooth
function. \eqref{lie.13.a} implies
\begin{equation}  \label{eq:fotimes21}\begin{aligned} &
  \sum _{\mu ,\nu} \| A_{\mu \nu}(x,\|  f\| _2^2)
\| _{H^{ K , S }} +\sum _{\mu ,\nu} \|
      \mathcal{A}_{\mu   \nu
}(\| f\| _2^2)
   \| _{H^{ K , S }} \le C \| \chi \| .
\end{aligned}
\end{equation}
We  consider
\begin{equation}  \label{eq:ExpH12}\begin{aligned} & H_2^{(1)}\circ
\phi+\widetilde{\resto ^{(1)}} \circ \phi +\int _{\mathbb{R}^3}
F_2(x,0 ,0, 0,\| f\| _2^2) f '^{\otimes 2}(x) dx\\&  + \langle
\nabla _f ^2 \widehat{\resto }^{(1)}_2 (0,0,\| f\| _2^2), f
'^{\otimes 2}\rangle .
\end{aligned}
\end{equation}
We will assume for the moment Lemma \ref{lem:1step1}:
\begin{lemma}
  \label{lem:1step1} The following difference is formed by terms which
   satisfy the properties stated for $\resto ^{(2)}$ in Theorem
   \ref{th:main}:
  \begin{equation} \label{eq:lem:1step1}
    \psi (\|f'\| _2^2)+ (H_2^{(1)} +\widetilde{\resto ^{(1)}}
      +\widetilde{\resto ^{(2)}})
 \circ \phi -\psi (\|f \| _2^2)-\text{\eqref{eq:ExpH12}} .
  \end{equation}
\end{lemma}
We postpone the proof of Lemma \ref{lem:1step1} and focus on
\eqref{eq:ExpH12} and on the choice of $\chi$.
\begin{lemma}
  \label{lem:2step1}  It is possible to choose $\chi$ such that
there exists $H_2^{(2)}$ as in (i) Theorem
   \ref{th:main} such that the difference
\eqref{eq:ExpH12}$-H_2^{(2)}$ is formed by terms which
   satisfy the properties   stated for $\resto ^{(2)}$ in Theorem
   \ref{th:main}.
\end{lemma}
\proof We have by   \eqref{chi.1} and using Definition
\ref{def:PoissonFunct}
\begin{equation}  \label{eq:ExpH13}\begin{aligned} & H_2^{(1)}\circ
\phi  = H_2^{(1)} +\int _0^1 \{ H_2^{(1)}, \chi \} \circ \phi _t dt
= H_2^{(1)}+ \\& \sum_{|\mu +\nu |=2}b_{\mu   \nu }(\| f\| _2^2)\int
_0^1 \{ H_2^{(1)}, z^{\mu} \overline{z}^{\nu}  \} \circ \phi _t dt
+\\&   \sum_{|\mu +\nu |=1 }
 \langle   \sigma _3\sigma _1B_{\mu   \nu
}(\| f\| _2^2)
  , \int _0^1 \{ H_2^{(1)},z^{\mu}
\overline{z}^{\nu}  f   \} \circ \phi _t dt \rangle  + \widetilde{R}
\end{aligned}
\end{equation}
with $ |\widetilde{R}|\le C (|z|+\| f \| _{H^{ -K',-S'}})^3,$
\eqref{PoissBra2}, Lemma \ref{lie_trans}. Then, by \eqref{PoissBra2}
for $\lambda =\lambda ^{(1)}$,  defined in \eqref{eq:lambda},
and substituting  $H_2^{(1)}= D_2^{(1)} + (H_2^{(1)}- D_2^{(1)})$ in the last two lines of \eqref{eq:ExpH13}, we get

\begin{equation} \label{eq:ExpH133}\begin{aligned} & H_2^{(1)}\circ
\phi      = H_2^{(1)}
+ \im \sum_{|\mu +\nu |=2}b_{\mu   \nu } \lambda \cdot (\mu - \nu)
z^{\mu} \overline{z}^{\nu}  \\& -\im \sum_{|\mu +\nu |=1} z^{\mu}
\overline{z}^{\nu} \langle f ,\sigma _1\sigma _3 (\mathcal{H}
-\lambda   \cdot (\mu - \nu) )B_{\mu   \nu }  \rangle +
{\widehat{R}},
\end{aligned} \nonumber
\end{equation}
with $D_2^{(1)}$ defined in \eqref{eq:Diag} and with, by
\eqref{PoissBra2}, \eqref{eq:ExpHcoeff1} and Lemma \ref{lie_trans},
\begin{equation} \label{eq:ExpH1330} |\widehat{R}|\le C (|z|+\| f \|
_{H^{ -K',-S'}})^3 +C \| \chi \| (\| \chi \| + \| f \| _{ 2}^2 )
(|z|+\| f \| _{H^{-K',-S'}})^2.  \end{equation} Similarly
\begin{equation} \label{eq:ExpH15}\begin{aligned}
&\widetilde{\resto ^{(1)}} \circ \phi   =\widetilde{\resto
^{(1)}}+\sum_{|\mu +\nu |=2}\widetilde{l}_{\mu   \nu }   z^{\mu}
\overline{z}^{\nu} + \sum_{|\mu +\nu |=1} z^{\mu} \overline{z}^{\nu}
\langle f ,\sigma _1\sigma _3 \widetilde{L}_{\mu \nu } \rangle
+\underline{R},
\end{aligned}
\end{equation}
\begin{equation} \label{eq:ExpH16}\begin{aligned}
  \text{with } |\widetilde{l}_{\mu   \nu }|+ \| \widetilde{L}_{\mu
\nu } \| _{H^{K,S}} \le C  \| \chi \| \, \| \widetilde{\resto
^{(1)}}\|
\end{aligned}
\end{equation}
\begin{equation}\label{eq:ExpH17} | \underline{{R}}|\le
 \text{rhs\eqref{eq:ExpH1330}}+\text{rhs\eqref{eq:ExpH16}}.
 \end{equation}
In \eqref{eq:ExpH12} we substitute $f'^{\otimes 2}$ using
\eqref{eq:fotimes2}. Then
\begin{equation}  \label{eq:ExpH18}\begin{aligned} &  \int
_{\mathbb{R}^3} F_2(x,0 ,0, 0,\| f\| _2^2) f '^{\otimes 2}(x) dx =
\widetilde{\chi } +\mathrm{R}
\end{aligned}
\end{equation}
with: $\widetilde{\chi }$ a polynomial like \eqref{chi.1} with
$M_0=1$ such that  $\| \widetilde{\chi } \| \le C\| f\| _2^2 \|
{\chi } \| $  by  claims (4) and (9) in Lemma \ref{lem:ExpH} and by
\eqref{eq:fotimes2}; $\widetilde{\chi }$  satisfies \eqref{chi.11}
by the fact that the rhs\eqref{eq:ExpH18} is real valued;
$\mathrm{R}$ formed by terms with the properties stated for $
{\resto ^{(2)}} $ in Theorem \ref{th:main}, see second line of
\eqref{eq:fotimes2}. By an argument  similar to the one for
\eqref{eq:ExpH18}, we have
\begin{equation}  \label{eq:ExpH182}\begin{aligned} &
 \langle
\nabla _f ^2 \widehat{\resto }^{(1)}_2 (0,0,\| f\| _2^2), f
'^{\otimes 2}\rangle = \widetilde{\chi } +\mathrm{R},
\end{aligned}
\end{equation}
with $\widetilde{\chi }$ and $\mathrm{R}$   different from the ones
in \eqref{eq:ExpH18}  but with the same properties. Then we have

\begin{equation} \label{eq:ExpH19}\begin{aligned}&
\text{\eqref{eq:ExpH12}}=H_2^{(1)} +\widetilde{\resto} ^{(2)}
+\widehat{\chi}     + \im \sum_{|\mu +\nu |=2}b_{\mu \nu } \lambda
\cdot (\mu - \nu) z^{\mu} \overline{z}^{\nu} -\\& -\im \sum_{|\mu
+\nu |=1} z^{\mu} \overline{z}^{\nu} \langle f ,\sigma _1\sigma _3
(\mathcal{H} -\lambda   \cdot (\mu - \nu) )B_{\mu   \nu }  \rangle +
\mathbf{R} ,
\end{aligned}
\end{equation}
where $\mathbf{R} $ satisfies  the properties stated for $ {\resto
^{(2)}}$ and $\widehat{\chi}$   is a  polynomial like
\eqref{chi.1}--\eqref{chi.11} with $M_0=1$ and ($\widehat{Z}$ and
$\widehat{K}$ will be defined in two lines)\begin{equation}
\label{eq:ExpH20}\begin{aligned}&\| \widehat{\chi } \| =  \|
\widehat{Z} \| + \| \widehat{K } \| \le C \| {\chi } \| (\| f\| _2^2
+\| \chi \|   + \| \widetilde{\resto} ^{(2)}\|   ) .
\end{aligned}\end{equation}
Here $\widehat{\chi }= \widehat{Z}+\widehat{K}$, where in
$\widehat{Z}=\sum   \widehat{b}_{\mu  \nu}(\| f\| _2^2) z^\mu
\overline{z}^\nu $ we sum over $|\mu +\nu |=2$, $\lambda ^0\cdot \mu
=\lambda ^0\cdot \nu$, i.e. in $\widehat{Z}$ we collect the null
form terms of $\widehat{\chi }$. We set
\begin{equation} \label{eq:ExpH2011}  H_2^{(2)}=
H_2^{(1)}+\widehat{Z} .
\end{equation}
Up to now $\chi$ is undetermined.  We choose $\chi$  with
coefficients $b_{\mu \nu}$ and $B_{\mu \nu } $ such that ${b}_{\mu
\nu}=0$ for $\lambda ^0\cdot \mu =\lambda ^0\cdot \nu$ and such that
the following system is satisfied:

\begin{equation} \label{eq:ExpH201}\begin{aligned}&
\widetilde{\resto ^{(1)}}
 +\widehat{K}     + \im \sum_{|\mu +\nu |=2}b_{\mu \nu}
\lambda \cdot (\mu - \nu) z^{\mu} \overline{z}^{\nu} -\\& -\im
\sum_{|\mu +\nu |=1} z^{\mu} \overline{z}^{\nu} \langle f ,\sigma
_1\sigma _3 (\mathcal{H} -\lambda   \cdot (\mu - \nu) )B_{\mu   \nu
}  \rangle =0.
\end{aligned}
\end{equation}
In coordinates, \eqref{eq:ExpH201} is
\begin{equation} \label{eq:ExpH21}\begin{aligned}&
a_{\mu \nu} ^{(1)} +\widehat{k}_{\mu \nu}  +\im b_{\mu \nu} \lambda
\cdot (\mu - \nu)=0 , \,   |\mu +\nu |=2,  \, \lambda ^0\cdot \mu
\neq \lambda ^0\cdot \nu , \\& G_{\mu \nu } +\widehat{K} _{\mu \nu}
  -  \im   (\mathcal{H}
-\lambda   \cdot (\mu - \nu) )B_{\mu   \nu }   =0 , \,   |\mu +\nu
|=1,
\end{aligned}
\end{equation}
where: $a_{\mu \nu}^{(1)}$ and $G_{\mu   \nu } $ are coefficients of
$\widetilde{\resto ^{(1)}} $, they are smooth functions of $\varrho
=\| f\| _2^2$, and are equal to 0 for $\varrho =0$; $\widehat{k}_{\mu \nu}\in \C$
and $\widehat{K}_{\mu   \nu } \in H ^{K,S}$ are coefficients of
$\widehat{K}$, and   are smooth functions of $\varrho =\| f\| _2^2$
and of the coefficients of $\chi$, where $b_{\mu \nu}\in \C$ and
$B_{\mu   \nu } \in H ^{K,S}$. By \eqref{eq:ExpH20}
\begin{equation} \label{eq:ExpH211} |\widehat{k}_{\mu \nu}|+ \|\widehat{K}_{\mu \nu}\|
_{H^{K,S}}\le  C \| {\chi } \| (\| f\| _2^2  +\| \chi \|   + \|
\widetilde{\resto ^{(1)}}  \|   ) .
\end{equation}
Then by the implicit function theorem we can solve the nonlinear
system \eqref{eq:ExpH21} with unknown $\chi$ obtaining (we consider
$b_{\mu \nu}$ only for $\lambda ^0\cdot \mu \neq \lambda ^0\cdot
\nu$)
\begin{equation} \label{eq:ExpH22}\begin{aligned}& |
b_{\mu \nu} +\frac{a_{\mu \nu}^{(1)}}{\im   \lambda \cdot (\mu -
\nu)} |+
 \| B_{\mu \nu} + \im R_{\mathcal{H}}(\lambda \cdot (\mu - \nu))
G_{\mu \nu}    \| _{H^{K,S}}\\& \le C \|\widetilde{\resto ^{(1)}} \|
(\| \widetilde{\resto ^{(1)}} \|  +\| f\| _2^2 ) .
\end{aligned}
\end{equation}
  Notice that  with the above choice of
  $\chi$ and with \eqref{eq:ExpH2011}, \eqref{eq:ExpH19}
yields

\begin{equation} \label{eq:ExpH190}\begin{aligned}&
\text{\eqref{eq:ExpH12}}=H_2^{(2)}   + \mathbf{R} ,
\end{aligned}
\end{equation}
where $\mathbf{R}$ has the properties stated for $ {\resto ^{(2)}} $
in Theorem \ref{th:main}. Hence Lemma \ref{lem:1step1} is proved.
\qed

\bigskip
{\it Proof of Lemma \ref{lem:1step1}}

By \eqref{lie.11.h}--\eqref{lie.11.f}, and with the big O smooth in
$z\in \C ^m$, $f\in H^{-K',-S'}_c$,

\begin{equation}  \label{eq:psi}
\psi (\|f'\| _2^2)  = \psi (\|f \| _2^2) +O\left ( |z|^{2 } \| f \|
_{H^{-K',-S'}} +   \| f \| _{H^{-K',-S'}} ^3 \right ).\end{equation}
The error term in \eqref{eq:psi} has the properties stated for $
{\resto ^{(2)}} $ in Theorem \ref{th:main}. We consider the terms
$\widetilde{\resto ^{(2)}} \circ \phi $. Terms, for $|\mu +\nu |  =
3$, like
\begin{equation} \label{eq:z3}z'^\mu \overline{z}'^\nu \int
_{\mathbb{R}^3}a (x,z',f',f'(x),\| f'\| _2^2 ) dx , \end{equation}
by \eqref{lie.11.a} and \eqref{lie.11.c} can be written as
\begin{equation}
\label{eq:z31}\begin{aligned} & (z ^\mu \overline{z} ^\nu  +O ((
|z|+ \norma{f } _{H^{-K',-S'}} )^3) )\int _{\mathbb{R}^3}a
(x,z',f',f'(x),\| f'\| _2^2 ) dx ,
\end{aligned}\end{equation} In the notation of Lemma \ref{lie_trans} we have
\begin{equation} \label{eq:z32}\begin{aligned}&
a( x,z ' ,f',f'(x ),\| f'\| _2^2) = a\big ( x,z  +
 \Gamma
,e^{ \im \Gamma _0  P_c(\omega _0) \sigma _3}f +\mathcal{G},\\&  e^{
\im \Gamma _0 \sigma _3}f (x)+[T(\Gamma _0) f](x)+\mathcal{G}(z,f,\|
f\| _2^2)(x) ,\| f \| _2^2+\Gamma _1 \big ) \\& = a( x,z   ,f ,f (x ),\|
f \| _2^2)+  O (  |z|+ \norma{f } _{H^{-K',-S'}}   ) .
\end{aligned}\end{equation}
The big O's in \eqref{eq:z31}--\eqref{eq:z32} are smooth in $z\in \C
^m$, $f\in H^{-K',-S'}_c$. Then \eqref{eq:z3} has the properties
stated for $ {\resto ^{(2)}} $ in Theorem \ref{th:main}. Similar
formulas can be used for
\begin{equation} \label{eq:zf3}\begin{aligned}&
\sum _{|\mu +\nu | =2  }z'^\mu \overline{z}'^\nu \int
_{\mathbb{R}^3} \left [ \sigma _1 \sigma _3G_{\mu \nu }
(x,z',f',f'(x), \| f'\| _2^2 )\right ]^*f'(x) dx +\\& \sum
_{j=3}^5\int _{\mathbb{R}^3} F_j(x,z' ,f', f'(x),\| f'\| _2^2)
f'^{\otimes j}(x) dx +
    \int _{\mathbb{R}^3}B (|f'(x)|^2/2 ) dx.
\end{aligned}\end{equation}
We treat with some detail these terms in the step $r>2,$ Subsection
\ref{subsec:step2}. Next we consider the term with $\int F_2
f'^{\otimes 2}(x) dx$. First of all, we can apply to $F_2$ an
analogue of \eqref{eq:z32} to obtain for $d=2$
\begin{equation} \label{eq:F32}\begin{aligned}&
F_d( x,z ' ,f',f'(x ),\| f'\| _2^2) = F_d( x,z  +
 \Gamma
,e^{ \im \Gamma _0  P_c(\omega _0) \sigma _3}f +\mathcal{G},\\&  e^{
\im \Gamma _0 \sigma _3}f (x)+[T(\Gamma _0) f](x)+\mathcal{G}(z,f,\|
f\| _2^2)(x) ,\| f \| _2^2+\Gamma _1 ) \\& = F_d ( x,0   ,0 ,f (x
),\| f \| _2^2)+  O (  |z|+ \norma{f } _{H^{-K',-S'}}   ) .
\end{aligned}\end{equation}
Then, modulo terms with the properties stated for $ {\resto ^{(2)}}
$ in Theorem \ref{th:main}, we get
\begin{equation} \label{eq:f2}\begin{aligned}&
 \int _{\mathbb{R}^3} F_2(x,0 ,0, f (x
),\| f \| _2^2) f'^{\otimes 2}(x) dx=
 \int _{\mathbb{R}^3} F_2(x,0 ,0, 0,\| f \| _2^2)
f'^{\otimes 2}(x) dx
\\&+ \int _{\mathbb{R}^3}
G_2(x, f (x),\| f \| _2^2) f(x)\otimes  f'^{\otimes 2} (x) dx   ,
\end{aligned}\end{equation}
where first term in rhs has been treated in Lemma \ref{lem:2step1}
and second term has the properties stated for $  \resto ^{(2)} _3 $
in Theorem \ref{th:main}. By a similar argument
\begin{equation} \label{eq:f22}\begin{aligned}&
\widehat{\resto} ^{(1)} _2 (z',f', \| f'\| _2 )- \langle \nabla _f
^2 \widehat{\resto }^{(1)}_2 (0,0,\| f\| _2^2), f '^{\otimes
2}\rangle
\end{aligned}\end{equation}
has the properties stated for $  \resto ^{(2)}   $ in Theorem
\ref{th:main}.
 \qed

We denote: $\chi _2=\chi$, $\Tr _2$ the   Lie transformation of
$\chi_2$, $Z^{(2)}=0 $.   $H_2^{(2)} $ has been defined in
\eqref{eq:ExpH2011}. We denote $\resto ^{(2)}=
\text{\eqref{eq:lem:1step1}} +\text{\eqref{eq:ExpH12}}- H_2^{(2)} $.
This $\resto ^{(2)}$ satisfies the conditions in Theorem
\ref{th:main}. This ends the proof of case $r=2$ in Theorem
\ref{th:main}.

\subsection{Proof of Theorem \ref{th:main}:  the step $r>2$}
\label{subsec:step2}

  Case $r=2$  has been  treated in Subsection \ref{subsec:step1}.
We have defined  $H_2^{(2)} $  in \eqref{eq:ExpH2011}.
 We proceed by
induction to complete the proof of Theorem \ref{th:main}. From the
argument below one can see that $H_2^{(r)} =H_2^{(2)} $ for all
$r\ge 2$. For $r\ge 2$, write  Taylor expansions
\begin{equation}
\label{r00}  \resto^{(r)}_{0}- \resto^{(r)}_{02} = \sum _{|\mu +\nu
|  = r+1 } z^\mu \overline{z}^\nu  \int _{\mathbb{R}^3}a_{\mu \nu
}^{(r)}(x,0,0,0,\| f\| _2^2 ) dx, \end{equation}
\begin{equation}
\label{r11}  \resto^{(r)}_{1}- \resto^{(r)}_{12} =\sum _{|\mu +\nu |
= r  }z^\mu \overline{z}^\nu \int _{\mathbb{R}^3} \left [ \sigma _1
\sigma _3G_{\mu \nu }^{(r)}(x,0,0,0, \| f\| _2^2 )\right ]^*f(x) dx.
\end{equation}
We have
 \begin{equation} \label{eq:r00} \begin{aligned} &
 \resto^{(r)}_{02} + \resto^{(r)}_{12} =\sum _{|\mu +\nu
|  = r+2 } z^\mu \overline{z}^\nu  \int
_{\mathbb{R}^3}\widetilde{a}_{\mu \nu }^{(r)}(x,z,f,0,\| f\| _2^2 )
dx+\\&   \sum _{|\mu +\nu |  = r+1  }z^\mu \overline{z}^\nu \int
_{\mathbb{R}^3} \left [ \sigma _1 \sigma _3\widetilde{G}_{\mu \nu
}^{(r)}(x,z,f,f(x), \| f\| _2^2 )\right ]^*f(x) dx + \\& \sum _{|\mu
+\nu |  = r   }z^\mu \overline{z}^\nu \int _{\mathbb{R}^3}
\widetilde{F}_{2 }^{(r)}(x,z,f,f(x), \| f\| _2^2 )
 \cdot  \left (f(x)\right ) ^{\otimes 2}  dx,
 \end{aligned}
  \end{equation}
with $\widetilde{a}_{\mu \nu }^{(r)}$ satisfying \eqref{eq:coeff a},
  $\widetilde{G}_{\mu \nu }^{(r)}$      \eqref{eq:coeff G}
and $\widetilde{F}_{2 }^{(r)}$ \eqref{eq:coeff F}. Since
$H^{(r)}=H\circ \Tr _r$ is real valued (because $H$ is real valued),
then both sides of equations \eqref{r00}--\eqref{eq:r00} are real
valued. In particular, $ {a}_{\mu \nu }^{(r)}$ and $ {G}_{\mu \nu
}^{(r)}$ satisfy \eqref{eq:ExpHcoeff2}. Set
\begin{equation}
\label{pr.3} \widetilde{K}_{r+1}:= \text{rhs\eqref{r00}}
+\text{rhs\eqref{r11}}.
\end{equation}
  Split $\widetilde{K}_{r+1}= K_{r+1}+ Z_{r+1}$ collecting inside
$Z_{r+1}$ all the terms of $\widetilde{K}_{r+1}$ in null form. The
coefficients of $ {K}_{r+1}$ and of  $Z_{r+1}$ satisfy
\eqref{eq:ExpHcoeff2}, by the argument just before \eqref{pr.3}. We
consider a (momentarily unknown) polynomial $\chi $ like
\eqref{chi.1}--\eqref{chi.11}, $M_0=r$. Denote by $ \phi $ its Lie
transformation. Let $(z',f')=\phi (z,f)$. For $d=2$, in the notation
of Lemma \ref{lie_trans} we have

\begin{equation} \label{binomials1}\begin{aligned}&
(\resto^{(r)}_d   - \widehat{\resto}^{(r)}_d)(z',f')  =
 \langle  F_d^{(r)}( z ' ,f',f'(\cdot ),\| f'\| _2^2),
(e^{ \im \Gamma _0 P_c(\omega _0)\sigma _3  }f +\mathcal{G}
)^{\otimes d} \rangle    .
\end{aligned}\end{equation}
Then rhs\eqref{binomials1}$=$
\begin{equation} \label{binomials}\begin{aligned}&
     =  \sum_{j=0}^{d}\left(
\begin{matrix}
d\\ j
\end{matrix}
\right) \langle  F_d^{(r)}( z ' ,f',f'(\cdot ),\| f'\| _2^2),
\mathcal{G}
 ^{\otimes (d-j)}\otimes [e^{ \im \Gamma _0 P_c(\omega _0)
 \sigma _3   }f]^{\otimes j}\rangle =\\& \sum_{j=0}^{d}\left(
\begin{matrix}
d\\ j
\end{matrix}
\right) \sum_{\ell =0}^{j}\left(
\begin{matrix}
j\\ \ell
\end{matrix}
\right)\langle  F_d^{(r)}( \cdots ), \mathcal{G}
 ^{\otimes (d-j)}\otimes [T (\Gamma _0) f]^{\otimes (j-\ell)}
 \otimes [e^{ \im \Gamma _0
 \sigma _3   }f]^{\otimes \ell }\rangle  .
\end{aligned}\end{equation}
In the notation of Lemma \ref{lie_trans} we have for $d=2$
\begin{equation} \label{Fd}\begin{aligned}&
F_d^{(r)}( z ' ,f',f'(x ),\| f'\| _2^2)(x)=\\& F_d^{(r)}\big ( z  +
 \Gamma
,e^{ \im \Gamma _0  P_c(\omega _0) \sigma _3}f +\mathcal{G}, e^{ \im
\Gamma _0 \sigma _3}f (x)+[T(\Gamma _0) f](x) ,\| f \| _2^2+\Gamma
_1 \big )(x) .
\end{aligned}\end{equation}
Then
\begin{equation} \label{Fd2}\begin{aligned}&
F_2^{(r)}( z ' ,f',f'(x ),\| f'\| _2^2)(x)= F_2^{(r)}( 0 ,0, f(x)
,\| f \| _2^2  )(x) +\\&    O (  |z|+ \norma{f } _{H^{-K',-S'}}   )
= F_2^{(r)}( 0 ,0,0 ,\| f \| _2^2  )(x) +\\&    G ( 0 ,0, f(x) ,\| f
\| _2^2  )(x) f(x)+  O (  |z|+ \norma{f } _{H^{-K',-S'}}   ) ,
\end{aligned}\end{equation}
where the big O are smooth in $z\in \C^m$ and $f\in H^{-K',-S'}$
with values in $H^{K,S} (\mathbb{R}^3,   B   (
 (\mathbb{C}^2)^{\otimes 2},\mathbb{C} )$ and where $G$ has values
 in $H^{K,S} (\mathbb{R}^3,   B   (
 (\mathbb{C}^2)^{\otimes 3},\mathbb{C} )$ and satisfies estimates
 \eqref{eq:coeff F}. So the last line of \eqref{Fd2}
 when plugged in  \eqref{binomials} for $d=2$ yields terms
 with the properties of $\sum _{d=0}^3\resto _d ^{(r+1)}.$
We focus now on the first term in the
 rhs of \eqref{Fd2}. Schematically, in analogy to  \eqref{eq:fotimes2}
we write
 \begin{equation}  \label{eq:fotimes22}\begin{aligned} &
f '^{2}(x)= \sum _{|\mu +\nu |=r} z^{\mu}\overline{z}^{\nu}
      \mathcal{A}_{\mu   \nu
}(\| f\| _2^2)(x)
    f (x) \\&
+\sum _{|\mu +\nu |=2r}A_{\mu \nu}(x,\|  f\| _2^2)
z^{\mu}\overline{z}^{\nu}    +
    (e^{ \im \Gamma _0   \sigma _3}f + T(\Gamma _0)f)  ^{2}(x) \\& +
    \varphi (x) r_j f(x)+r_f(x) f(x)+ \varphi (x) r_j ^2+r_f^2(x),
\end{aligned}
\end{equation}
where we have \eqref{eq:fotimes21} and $|r_j|+\| r_f\| _{H^{K,S}}\le
C (|z| +\|  f\| _{H^{-K',-S'}}  ) ^{r+1}.$ Then
\begin{equation}  \label{eq:ExpH181}\begin{aligned} &  \int
_{\mathbb{R}^3} F_2^{(r)}(x,0 ,0, 0,\| f\| _2^2) f '^{\otimes 2}(x)
dx = \widetilde{\chi } _1 +\mathrm{R}_1
\end{aligned}
\end{equation}
with: $\mathrm{R}_1$ formed by terms obtained by the last two lines
of \eqref{eq:fotimes22} has  the properties stated for $ {\resto
^{(r+1)}} $ in Theorem \ref{th:main}; $\widetilde{\chi }_1$ a
polynomial like \eqref{chi.1} with $M_0=r$ arising from the first
line of rhs of \eqref{eq:fotimes22}   is such that   $\| \widetilde{\chi }_1\| \le C\| f\| _2^2 \|  \chi  \|$ by the inductive  hypothesis
$F_2^{(r)}(x,0 ,0, 0,0 )=0$ in (iv.2-5) Theorem \ref{th:main} and by \eqref{eq:fotimes21};
$\widetilde{\chi }_1$ satisfies \eqref{chi.11} because each side in
\eqref{eq:ExpH181} is real valued. We have
\begin{equation}  \label{eq:quadr}\begin{aligned} &
\widehat{\resto }^{(r)}_2(z',f',\| f'\| _2^2) =
 \langle
\nabla _f ^2 \widehat{\resto }^{(1)}_2 (0,0,\| f\| _2^2), f
'^{\otimes 2}\rangle +\\& + (\widehat{\resto }^{(r)}_2(z',f',\| f'\|
_2^2)-\langle \nabla _f ^2 \widehat{\resto }^{(r)}_2 (0,0,\| f\|
_2^2), f '^{\otimes 2}\rangle ),
\end{aligned}
\end{equation}
where   the second line on rhs of \eqref{eq:quadr} yields terms which
have the properties of elements of ${\resto }^{(r+1)}$. We have
\begin{equation}  \label{eq:ExpH183}\begin{aligned} &
 \langle
\nabla _f ^2 \widehat{\resto }^{(r)}_2 (0,0,\| f\| _2^2), f
'^{\otimes 2}\rangle = \widetilde{\chi }_2 +\mathrm{R}_2
\end{aligned}
\end{equation}
where $\widetilde{\chi }_2$ and $\mathrm{R}_2$  have the same
properties of $\widetilde{\chi }_1$ and $\mathrm{R}_1$ in
\eqref{eq:ExpH183}. Split $H_2^{(r)}=D_2^{(r)}
+(H_2^{(r)}-D_2^{(r)})$ for $D_2^{(r)}$ in \eqref{eq:Diag}. Then
\begin{equation}
\label{eq:ExpH1831} \left\{ \chi  ,H_2 ^{(r)} -D_2^{(r)}\right\}
=\widetilde{\chi }_3 +\mathrm{R}_3
\end{equation}
where $\widetilde{\chi }_3$ and $\mathrm{R}_3$  have the same
properties of $\widetilde{\chi }_1$ and $\mathrm{R}_1$ in
\eqref{eq:ExpH183}. Set $ \widetilde{\chi }=\sum
_{j=1}^{3}\widetilde{\chi }_j$. Split now $\widetilde{\chi
}=\widetilde{Z}+\widehat{K}$ collecting in $\widetilde{Z}$ the null
form terms in $\widetilde{\chi }$. Then we choose  the yet unknown
$\chi$ such that its coefficients $b_{\mu \nu}$ and $B_{\mu   \nu }
$ satisfy the system

\begin{equation} \label{eq:ExpH20111}\begin{aligned}&
\widetilde{K}_{r+1} +\widehat{K}     + \im \sum_{|\mu +\nu
|=2}b_{\mu \nu} \lambda \cdot (\mu - \nu) z^{\mu} \overline{z}^{\nu}
-\\& -\im \sum_{|\mu +\nu |=1} z^{\mu} \overline{z}^{\nu} \langle f
,\sigma _1\sigma _3 (\mathcal{H} -\lambda   \cdot (\mu - \nu)
)B_{\mu   \nu }  \rangle =0.
\end{aligned}
\end{equation}
Notice that for $\widehat{K}\equiv 0$ system \eqref{eq:ExpH20111}
would be linear and admit  exactly one solution. By $\|
\widetilde{\chi } \| \le C\| f\| _2^2 \| {\chi } \| $ we get $\|
\widehat{K} \| \le C\| f\| _2^2 \| {\chi } \| $. So by the implicit
function theorem there exists exactly one solution of
\eqref{eq:ExpH20111}. This solution is close to the solution of system
\eqref{eq:ExpH20111} when $\widehat{K}\equiv 0$. Furthermore,
this system has solution $\chi_{r+1}=\chi$ which satisfies
\eqref{eq:ExpHcoeff2}, or what is the same, \eqref{chi.11}. For
$L_{r+1}$ of   type \eqref{eq:RestohomologicalEq}, $\chi_{r+1}$
satisfies

\begin{equation}
\label{HomEqmain} \left\{ \chi_{r+1} ,H_2 ^{(r)}\right\}
=\widetilde{K}_{r+1} +\widehat{K}+ L_{r+1}.
\end{equation}  Call $\phi_{r+1}= \phi$  the Lie transform
of $\chi_{r+1}$.   For $ \Tr_{r+1}=\Tr_r\circ\phi_{r+1}$ set
\begin{equation}
\label{Hamr} H^{(r+1)}:=H^{(r)}\circ\phi_{r+1}=
H\circ(\Tr_r\circ\phi_{r+1}) = H\circ \Tr_{r+1}.
\end{equation}
Since $\chi_{r+1}$  satisfies \eqref{eq:ExpHcoeff2}, $H^{(r+1)}$ is
well defined and real valued. Split
\begin{eqnarray}
\label{n1.1} & H^{(r)}\circ\Tra_{r+1}=  H_2^{(r)}+Z^{(r)} +
Z_{r+1}+\widetilde{Z}
\\ & \label{n1.3}
 +  (Z^{(r)}\circ\Tra _{r+1}-Z^{(r)}) +(\widetilde{Z} \circ\Tra
_{r+1} -\widetilde{Z}  ) \\ & \label{n1.6} + (\widetilde{K}_{r+1}
+\widehat{K})\circ\Tra_{r+1}-\widetilde{K}_{r+1} -\widehat{K}
\\ &
\label{n1.4}  + H_2^{(r)}\circ\Tra_{r+1}- \left(H_2^{(r)}+\left\{
H_2^{(r)}, \chi_{r+1}\right\}\right)
\\ &
\label{n1.5}  + (\resto^{(r)}_{02}+\resto^{(r)}_{12})\circ\Tra_{r+1}
\\ &
\label{n1.7} +    \sum _{d=3}^{5}(\resto ^{(r)}_d -\widehat{\resto}
^{(r)}_d) \circ\Tra_{r+1}+ \widehat{\resto} ^{(r)}_d \circ\Tra_{r+1}
\\ & \label{n1.9}  +
(\resto ^{(r)}_2 - \widetilde{Z}-\widehat{K})\circ\Tra_{r+1}
\\ & \label{n1.8}   +    \psi \circ\Tra_{r+1}+
  \resto ^{(r)}_6\circ\Tra_{r+1}
 \ .
\end{eqnarray}
Define $H^{(r+1)}_2=H^{(r )}_2 $  (this proves  $H^{(r )}_2=H^{(2
)}_2 $) and $Z^{(r+1)}:=Z^{(r)}+Z_{r+1}+\widetilde{Z}$. Its
coefficients satisfy \eqref{eq:ExpHcoeff2} (because $H^{(r+1)}$ is
real valued) and it is a normal form. We have already discussed that
\eqref{n1.9}  has the properties stated for $\resto ^{(r+1)}$. By
expansions \eqref{binomials1}--\eqref{Fd} we get that the first
summation  in \eqref{n1.7}      has the properties stated for
$\resto ^{(r+1)}$. By an analogous argument,   terms $
\widehat{\resto} ^{(r)}_d
 (z',f')$  have the properties stated for $\resto ^{(r+1)}$.
We have, for $T=T(\Gamma _0)$,
\begin{equation} \label{resto61}\begin{aligned}& | f' (x)|^2=
| f  (x)|^2 +\mathcal{E}(x)\text{ with } \mathcal{E}(x):=2 (T(\Gamma
_0)f(x))^*\sigma _1 e^{ \im \Gamma _0 \sigma _3}f(x)\\& + | T(\Gamma
_0)f(x)|^2+ 2\mathcal{G} ^* (x)\sigma _1 e^{ \im \Gamma _0 \sigma
_3}f(x)+2\mathcal{G} ^* (x)\sigma _1 T(\Gamma _0)f(x)+| \mathcal{G}
(x)|^2.
\end{aligned}\end{equation}
Then
\begin{equation} \label{resto6}\begin{aligned}& \resto^{(r)}_6
\circ \phi _{r+1}  =
 \int _{\mathbb{R}^3}  B ( | f' (x)| ^2/2) dx=
 \int _{\mathbb{R}^3}  B ( | f  (x)| ^2/2) dx\\& +\frac{1}{2}
 \int _{\mathbb{R}^3}  dx\,  \mathcal{E}(x)\int _0^1 B '( | f  (x)| ^2/2
 +s\, \mathcal{E}(x)/2) ds.
\end{aligned}\end{equation}
The last line in \eqref{resto6}   has the properties stated for
${\resto}^{(r+1)} -{\resto}^{(r+1)}_6$ by Lemma \ref{lie_trans}. By
\eqref{eq:r00} and by the fact that $\widetilde{a}_{\mu \nu }^{(r)}$
satisfies \eqref{eq:coeff a},
  $\widetilde{G}_{\mu \nu }^{(r)}$      \eqref{eq:coeff G}
and $\widetilde{F}_{2 }^{(r)}$ \eqref{eq:coeff F},  the terms
$\resto^{(r)}_{02}+\resto^{(r)}_{12}$  has the properties stated for
$\sum _{d=0}^2{\resto}^{(r+1)}_d$. The same conclusion holds for
\eqref{n1.5}.
  By Lemma
\ref{lie_trans} and by an analogue of \eqref{eq:psi},  we have that
$\psi \circ \phi _r=\psi + \widetilde{\psi} $ where
$\widetilde{\psi}$  has the properties stated for $\sum
_{d=1}^3\resto^{(r+1)}_d$ by \eqref{lie.11.f}. We have
\begin{equation} \label{eq:Zcirc}\begin{aligned} &  Z^{(r)}\circ\Tra_{r+1} -Z^{(r)}=
 \int_0^1  \{ Z^{(r)} ,   \chi_{r+1}  \}
  \circ\phi_{r+1}^t dt.
\end{aligned}  \end{equation}
We have
\begin{equation} \label{eq:chiZ}\begin{aligned} &  \left |
 \{\chi_{r+1},Z^{(r)} \} \right | \le C (|z|^{r+2} +|z|^{r+1}
 \| f\| _{H^{-K',-S'}} ).
\end{aligned}   \end{equation}
By   \eqref{eq:chiZ} we conclude that \eqref{eq:Zcirc}  has the
properties stated for $\resto^{(r+1)}$. The same  is true for the
other terms in \eqref{n1.3}--\eqref{n1.6}.   We have, for
$H_2=H_2^{(r)}$,

\begin{equation} \label{eq:Hcirc}\begin{aligned} &
  H_2\circ\phi_{r+1}-(H_2+\left\{ H_2,   \chi_{r+1}  \right\})=
  \int_0^1 \frac{t^2}{2!} \left\{
\left\{ H_2,   \chi_{r+1} \right\},   \chi_{r+1}
\right\}\circ\phi_{r+1}^t dt \\& =- \int_0^1 \frac{t^2}{2!}\left\{
K_{r+1}+\widehat{K} +L_{r+1},   \chi_{r+1} \right\}\circ\phi_{r+1}^t
dt.
\end{aligned}   \end{equation}
Then $\left |
 \{  K_{r+1}+\widehat{K}
+L_{r+1} ,   \chi_{r+1} \} \right | \le \text{rhs \eqref{eq:chiZ} }$
implies that \eqref{eq:Hcirc}  has the properties stated for
$\resto^{(r+1)}$.

 \qed

 \section{Dispersion}
 \label{sec:dispersion}
We apply Theorem \ref{th:main} for $r=2N+1 $  (recall $N=N_1$ where
$N_j\lambda _j<\omega _0 <(N_j+1)\lambda _j).$ In the rest of the
paper we work with the hamiltonian $H^{(r)}$. We will drop the upper
index. So we will set $H=H^{(r)}$, $H_2=H_2^{(r)}$, $\lambda
_j=\lambda _j^{(r)}$,    $\lambda  =\lambda  ^{(r)}$,
$Z_a=Z_a^{(r)}$ for $a=0,1$ and $\resto =\resto ^{(r)}$. In
particular we will denote by $G_{\mu \nu}$ the coefficients $G_{\mu
\nu}^{(r)}$ of $Z_1^{(r)}$. We will show:
\begin{theorem}\label{proposition:mainbounds} There is a fixed
$C >0$ such that for $\varepsilon _0>0$ sufficiently small and for
$\epsilon \in (0, \varepsilon _0)$ we have
\begin{eqnarray}
&   \|  f \| _{L^r_t( [0,\infty ),W^{ 1 ,p}_x)}\le
  C \epsilon \text{ for all admissible pairs $(r,p)$}
  \label{Strichartzradiation}
\\& \| z ^\mu \| _{L^2_t([0,\infty ))}\le
  C \epsilon \text{ for all multi indexes $\mu$
  with  $\lambda\cdot \mu >\omega _0 $} \label{L^2discrete}\\& \| z _j  \|
  _{W ^{1,\infty} _t  ([0,\infty )  )}\le
  C \epsilon \text{ for all   $j\in \{ 1, \dots , m\}$ }
  \label{L^inftydiscrete} .
\end{eqnarray}
\end{theorem}
Estimate \eqref{L^inftydiscrete} is a consequence of the classical
proof of orbital stability in Weinstein \cite{W1}. Notice that
\eqref{NLS}  is time reversible, so in particular
\eqref{Strichartzradiation}--\eqref{L^inftydiscrete} are true over
the whole real line. The proof, though, exploits that $t\ge 0$,
specifically when for $\lambda \in \sigma _c(\mathcal{H})$  we
choose $R_{\mathcal{H}}^+(\lambda )=R_{\mathcal{H}} (\lambda +\im 0
)$ rather than $R_{\mathcal{H}}^-(\lambda )=R_{\mathcal{H}} (\lambda
-\im 0 )$ in formula \eqref{eq:g variable}. See the discussion on
p.18 \cite{SW3}.

The proof of
 Theorem
\ref{proposition:mainbounds}  involves a  standard continuation
argument. We assume

\begin{eqnarray}
&   \|  f \| _{L^r_t([0,T],W^{ 1 ,p}_x)}\le
  C _1\epsilon \text{ for all admissible pairs $(r,p)$} \label{4.4a}
\\& \| z ^\mu \| _{L^2_t([0,T])}\le
 C_2 \epsilon \text{ for all multi indexes $\mu$
  with  $\omega \cdot \mu >\omega _0 $} \label{4.4}
\end{eqnarray}
for fixed sufficiently large  constants $C_1$, $C_2$ and then we
prove that  for $\epsilon $ sufficiently small, \eqref{4.4a} and
\eqref{4.4} imply the same estimate but with $C_1$, $C_2$ replaced
by $C_1/2$, $C_2/2$. Then \eqref{4.4a} and \eqref{4.4} hold   with
$[0,T]$ replaced by $[0,\infty )$.

The proof consists in three main steps.
\begin{itemize}
\item[(i)] Estimate $f$ in terms of $z$.
\item[(ii)] Substitute the variable $f$  with a
new "smaller" variable $g$ and find smoothing estimates for $g$.
\item[(iii)] Reduce the system for $z$ to a closed system involving
only the $z$ variables, by insulating the  part of $f$  which
interacts with $z$, and by decoupling the rest (this reminder is
$g$). Then clarify the nonlinear Fermi golden rule.
\end{itemize}
The first  two steps are the same of \cite{cuccagnamizumachi}. The
only novelty of the  proof with respect to \cite{cuccagnamizumachi}
 is step (iii), specifically the part on the Fermi golden rule. At issue is the non negativity of some crucial
 coefficients in the equations of $z$. This point is solved
 using the same ideas in  Lemma 5.2 \cite{bambusicuccagna}.
 The fact that they are not 0 is assumed by hypothesis (H11).
 The fact that if not 0 they are positive, is proved here.

Step (i) is encapsulated by the following proposition:

\begin{proposition}\label{Lemma:conditional4.2} Assume \eqref{4.4a}
  and \eqref{4.4}. Then there exist constants $C=C(C_1,C_2), K_1$,
  with $K_1$ independent of $C_1$, such that, if
  $C(C_1,C_2) \epsilon  $  is sufficiently small,   then we have
\begin{eqnarray}
&   \|  f \| _{L^r_t([0,T],W^{ 1 ,p}_x)}\le
  K_1  \epsilon \text{ for all admissible pairs $(r,p)$}\ .
  \label{4.5}
\end{eqnarray}
\end{proposition}
Consider $Z_1$ of the form \eqref{e.12a}. Set:
 \begin{equation}\label{eq:G^0} G_{\mu \nu}^0=G_{\mu
\nu}(\| f \| _2^2 ) \text{ for $\| f \| _2^2=0$;     $\lambda
^0_j=\lambda _j(\omega _0)$}.\end{equation} Then we have (with
finite sums and with the derivative in the variable $\| f \| _2^2$ performed w.r.t. the $\| f \| _2^2$ arguments explicitly emphasized in Theorem \ref{th:main})

\begin{equation}\label{eq:f variable} \begin{aligned}  &\im \dot f -
\mathcal{H}f - 2    (\partial _{ \| f \| _2^2} H)
P_c(\omega _0)\sigma _3 f = \sum _{|\lambda  ^0
\cdot(\nu-\mu)|>\omega _0} z^\mu \overline{z}^\nu
  G_{\mu \nu}^0  \\&  +
  \sum _{|\lambda ^0 \cdot(\nu-\mu)|>\omega _0} z^\mu
\overline{z}^\nu  (G_{\mu \nu} -  G_{\mu \nu}^0)    +\sigma _3
\sigma _1  \nabla _f \resto -  2    (\partial _{ \| f \| _2^2}  \resto)
P_c(\omega _0)\sigma _3 f   .
\end{aligned}\end{equation}
The proof of Proposition \ref{Lemma:conditional4.2} is standard and
is an easier version of the arguments in \S 4 in
\cite{cuccagnamizumachi}. The dominating  term in the rhs of
\eqref{eq:f variable} is the one  on the first line, with
contribution to $f$   bounded by $C(C_2) \epsilon $ by the
endpoint Strichartz estimate and by \eqref{4.4} (we recall   that
the third term in the lhs, in part becomes a phase through an
integrating factor, in part goes on the rhs: see
\cite{cuccagnamizumachi}; this trick is due to \cite{BP2}). Notice
also, that Theorem \ref{proposition:mainbounds} implies   by the
arguments on pp. 67--68 in \cite{cuccagnamizumachi}

 \begin{equation}\label{scattering11}  \lim_{t\to +\infty}
\left \|  e^{\im \theta (t) \sigma _3}f (t) -
 e^{ \im t \Delta   \sigma_3}{f}_+  \right \|_{H^1}=0
 \end{equation}
for a  $  f_+\in H^1$ with $ \| {f}_+    \|_{H^1}\le C \epsilon$ and
for    $ \theta(t) =    t\omega _0+2\int _0^t (\partial _{ \| f \| _2^2} H) (t')dt' .$   We claim that $\theta(t)=\vartheta (t)-\vartheta (0)$.
 This claim, Theorem \ref{th:main}, Theorem \ref{proposition:mainbounds}
and \eqref{scattering11} imply Theorem \ref{theorem-1.2}.  To prove the claim
  we substitute the
last system of coordinates in \eqref{system2} to obtain
 \begin{equation}\label{scattering12}   \im \dot f -\mathcal{H}f-(\dot \vartheta -\omega _0) P_c(\omega _0)\sigma _3 f =G
 \end{equation}
where $G$  is a functional with values in $  \in C(\R,  L^1_ {x}) $. The two equations are equivalent. This implies    $G=\text{rhs\eqref{eq:f variable}}
 $ and  $\dot \vartheta -\omega _0= 2     \partial _{ \| f \| _2^2}  H. $
This yields the claim  $\theta(t)=\vartheta (t)-\vartheta (0)$.

Step (ii) in the proof of Theorem \ref{proposition:mainbounds}
consists in introducing the variable

\begin{equation}
  \label{eq:g variable}
g=f+ \sum _{|\lambda ^0 \cdot(\mu-\nu)|>\omega _0} z^\mu
\overline{z}^\nu
   R ^{+}_{\mathcal{H}}  (\lambda  ^0 \cdot(\mu-\nu) )
    G_{\mu \nu}^0 .
\end{equation}
Substituting the new variable  $g$ in \eqref{eq:f variable}, the
first line on the rhs of  \eqref{eq:f variable} cancels out. By an
easier version of Lemma 4.3 \cite{cuccagnamizumachi} we have:

\begin{lemma}\label{lemma:bound g} For $\epsilon$ sufficiently small
 and for $C_0=C_0(\mathcal{H})$  a fixed constant, we have
\begin{equation} \label{bound:auxiliary}\| g
\| _{L^2_tL^{2,-S}_x}\le C_0 \epsilon + O(\epsilon
^2).\end{equation}
\end{lemma}
As in \cite{cuccagnamizumachi}, the part of $f$ which couples
nontrivially with $z$ comes from  the polynomial in $z$  contained
in  \eqref{eq:g variable}. $g$ and $z$ are decoupled.

\subsection{The Fermi golden rule}
\label{subsec:FGR}

We proceed as in the related parts in
\cite{bambusicuccagna,cuccagnamizumachi}. The only difference with
    \cite{cuccagnamizumachi} is that the preparatory work in Theorem
    \ref{th:main} makes transparent the positive semidefiniteness
    of the crucial coefficients.

Set $R_{\mu \nu }^+=R_{ \mathcal{H} }^+ (\lambda ^0\cdot (\mu   -\nu
) ).$ We will have $\lambda _j^0=\lambda _j (\omega _0)$ and
$\lambda _j=\lambda _j (\| f\| _2^2)$ as in Section
\ref{subsec:Normal form}.
  $|\lambda _j^0- \lambda _j|\lesssim  C_1^2\epsilon ^2$
by \eqref{4.4a}, so  in the sequel we can assume that   $\lambda ^0$
satisfies the same inequalities of $\lambda  .$ We substitute
\eqref{eq:g variable} in $\im \dot z_j = \frac{\partial}{\partial
\overline{z}_j} H^{(r)}$ obtaining
\begin{equation}\label{eq:FGR0} \begin{aligned} & \im \dot z _j
=
\partial _{\overline{z}_j}(H_2+Z_0) +   \sum  _{
 |\lambda  \cdot (\mu   -\nu )| > \omega
_0 }  \nu _j\frac{z ^\mu
 \overline{ {z }}^ { {\nu} } }{\overline{z}_j}  \langle g ,
 \sigma _1\sigma _3 G
_{\mu \nu }\rangle       +
\partial _{  \overline{z} _j}  \resto
 \\ &  - \sum  _{ \substack{| \lambda  \cdot
(\alpha    -\beta )|> \omega _0
\\
 |\lambda  \cdot (\mu -\nu )|> \omega
_0 }}  \nu _j\frac{z ^{\mu +\alpha }  \overline{{z }}^ { {\nu}
+\beta}}{\overline{z}_j} \langle   R_{   \alpha \beta}^+G^0 _{ \alpha
\beta },\sigma _1 \sigma _3G  _{\mu \nu }\rangle     .
\end{aligned}  \end{equation}
We rewrite this as
\begin{eqnarray} \label{equation:FGR1}& \im \dot z _j= \partial _{\overline{z}_j}(H_2+Z_0) +
     \mathcal{E}_j
\\ &  \label{equation:FGR12} -\sum  _{ \substack{ \lambda
\cdot  \beta    > \omega _0
\\
 \lambda   \cdot   \nu > \omega
_0  \\ \lambda \cdot  \beta -\lambda _k    < \omega _0   \, \forall
\, k \, \text{ s.t. } \beta _k\neq 0\\ \lambda  \cdot  \nu -\lambda
_k <\omega _0  \, \forall \, k \, \text{ s.t. } \nu _k\neq 0}} \nu
_j\frac{ \overline{{z }}^ {\nu +\beta } }{\overline{z}_j}\langle
R_{ 0 \beta}^+
 { G} _{   0\beta  }^0, \sigma _1 \sigma _3G ^0_{0 \nu }\rangle
 \\ &  \label{equation:FGR13} -\sum  _{ \substack{ \lambda
\cdot  \alpha    > \omega _0
\\
 \lambda  \cdot   \nu > \omega
_0  \\ \lambda \cdot  \alpha -\lambda _k   < \omega _0   \, \forall
\, k \, \text{ s.t. } \alpha _k\neq 0\\ \lambda \cdot  \nu -\lambda
_k   <\omega _0 \, \forall \, k \, \text{ s.t. } \nu _k\neq 0}} \nu
_j\frac{z ^{ \alpha } \overline{{z }}^ {\nu }
}{\overline{z}_j}\langle R_{ \alpha 0 }^+
 G _{   \alpha 0}^0, \sigma _1 \sigma _3G^0 _{0 \nu }\rangle  .
\end{eqnarray}
Here the elements in \eqref{equation:FGR12} will be eliminated
through a new change of variables. $\mathcal{E}_j$ is a reminder
term defined by

\begin{equation}   \begin{aligned} & \mathcal{E}_j:=
\text{rhs\eqref{eq:FGR0}} -\text{\eqref{equation:FGR12}}- \text{\eqref{equation:FGR13}}   .
\end{aligned} \nonumber
\end{equation}
Set
\begin{equation}\label{equation:FGR2}  \begin{aligned}   &
\zeta _j =z _j -\sum  _{ \substack{ \lambda \cdot  \beta    > \omega
_0
\, , \,
 \lambda \cdot   \nu > \omega
_0  \\  \lambda \cdot  \beta -\lambda _k   < \omega _0   \, \forall
\, k \, \text{ s.t. } \beta _k\neq 0\\ \lambda \cdot  \nu -\lambda
_k  <\omega _0   \, \forall \, k \, \text{ s.t. } \nu _k\neq 0}}
\frac{ \nu _j}{\lambda ^0 \cdot (\beta  +\nu) } \frac{ \overline{{z
}}^ {\nu +\beta } }{\overline{z}_j}\langle R_{   0 \beta}^+
 { G} _{  0 \beta  }^0, \sigma _1 \sigma _3G _{0 \nu }^0\rangle \\&
+ \sum  _{ \substack{ \lambda \cdot  \alpha    > \omega _0
\, , \,
 \lambda \cdot   \nu > \omega _0\\
 \lambda ^{0}\cdot  \alpha   \neq
 \lambda ^{0}\cdot   \nu
 \\ \lambda
\cdot  \alpha -\lambda _k   < \omega _0   \,
 \forall \, k \, \text{ s.t.
} \alpha _k\neq 0\\ \lambda \cdot  \nu -\lambda _k  <\omega _0   \,
\forall \, k \, \text{ s.t. } \nu _k\neq 0}}\frac{ \nu _j}{\lambda
^0 \cdot (\alpha  - \nu) } \frac{z ^{ \alpha } \overline{ z}^ { \nu
}}{\overline{z}_j} \langle R_{ \alpha 0 }^+ G^0 _{ \alpha 0},\sigma
_1 \sigma _3 G _{0 \nu }^0\rangle \end{aligned}
\end{equation}
Notice that  in \eqref{equation:FGR2}, by $\lambda \cdot \nu >
\omega _0$, we have $| {\nu} |>1$. Then by \eqref{4.4}

\begin{equation}  \label{equation:FGR3} \begin{aligned}   & \| \zeta  -
 z  \| _{L^2_t} \le C \epsilon \sum _{\substack{ \lambda
  \cdot \alpha    > \omega
_0   \\ \lambda \cdot  \alpha -\lambda _k   < \omega _0   \, \forall
\, k \, \text{ s.t. } \alpha _k\neq 0 }} \| z ^{\alpha }\| _{L^2_t}
\le CC_2M\epsilon ^2\\&  \| \zeta  -
 z \| _{L^\infty _t} \le C ^3\epsilon ^3
\end{aligned}
\end{equation}
with $C$ the constant in \eqref{L^inftydiscrete} and $M$ the number
of terms in the rhs. In the new variables \eqref{equation:FGR1} is
of the form

\begin{equation} \label{equation:FGR4} \begin{aligned} &
   \im \dot \zeta
 _j=
\partial _{\overline{\zeta}_j}H_2 (\zeta , f ) +
\partial _{\overline{\zeta}_j}Z_0 (\zeta ,  f )+  \mathcal{D}_j
 \\ &   -\sum  _{ \substack{ \lambda ^0
\cdot  \alpha =\lambda ^0\cdot   \nu   > \omega _0
   \\ \lambda
\cdot  \alpha -\lambda _k   < \omega _0   \, \forall \, k \, \text{
s.t. } \alpha _k\neq 0\\ \lambda \cdot  \nu -\lambda _k  <\omega _{0}  \,
\forall \, k \, \text{ s.t. } \nu _k\neq 0}} \nu _j \frac{\zeta ^{
\alpha } \overline{ \zeta}^ { \nu }}{\overline{\zeta}_j} \langle R_{
\alpha 0 }^+ G ^0_{ \alpha 0},\sigma _1 \sigma _3 G ^0_{0 \nu
}\rangle .
\end{aligned}
\end{equation}
From these equations by $\sum _j \lambda _j ^0 ( \overline{\zeta}_j
\partial _{\overline{\zeta}_j}(H_2+Z_0) - {\zeta}_j
\partial _{ {\zeta}_j}(H_2+Z_0)    ) =0$  we get

\begin{equation} \label{eq:FGR5} \begin{aligned}
 &\partial _t \sum _{j=1}^m \lambda  _j ^0
 | \zeta _j|^2  =  2  \sum _{j=1}^m \lambda _j^0\Im \left (
\mathcal{D}_j \overline{\zeta} _j \right ) -\\&    -2 \sum _{
\substack{ \lambda ^0 \cdot  \alpha =\lambda ^0\cdot   \nu   >
\omega _0
   \\ \lambda
\cdot  \alpha -\lambda _k   < \omega _0   \, \forall \, k \, \text{
s.t. } \alpha _k\neq 0\\ \lambda \cdot  \nu -\lambda _k  <\omega _0
\, \forall \, k \, \text{ s.t. } \nu _k\neq 0}} \lambda ^0\cdot \nu
\Im \left ( \zeta ^{ \alpha } \overline{\zeta }^ { \nu  } \langle
R_{ \alpha 0}^+  G_{ \alpha 0}^0,  \sigma _1\sigma _3  G ^0 _{0\nu
}\rangle \right ) .
\end{aligned}
\end{equation}
  We have the following lemma,   whose proof (we skip) is similar to Appendix B
\cite{bambusicuccagna}.
\begin{lemma}
\label{lemma:FGR1} Assume  inequalities \eqref{4.4}. Then for a
fixed constant $c_0$ we have
\begin{eqnarray}\label{eq:FGR7} \sum _j\|\mathcal{D}_j \overline{\zeta} _j\|_{
L^1[0,T]}\le (1+C_2)c_0 \epsilon ^{2}
 . \end{eqnarray}
\end{lemma}
   For the sum in the second line of \eqref{eq:FGR5}
we get

\begin{equation} \label{eq:FGR8} \begin{aligned} & 2\sum _{r>\omega
_0 } r
    \Im \left   \langle R_{
\mathcal{H}}^+ (r )\sum _{  \lambda ^0\cdot \alpha =r }\zeta ^{
\alpha } G _{  \alpha 0}^0, \sigma _1 \sigma _3\sum _{  \lambda
^0\cdot \nu =r}\overline{\zeta} ^{ \nu } G ^0_{0\nu
 }  \right \rangle      =\\&  2\sum _{r>\omega
_0 } r     \Im \left    \langle R_{ \mathcal{H}}^+ (r )\sum _{
\lambda ^0\cdot \alpha =r }\zeta ^{ \alpha } G _{  \alpha 0}^0,
\sigma _3\overline{\sum _{  \lambda ^0\cdot \alpha =r}\zeta ^{
\alpha } G ^0_{  \alpha 0} }\right \rangle    ,
\end{aligned}
\end{equation}
where we have used $G_{\mu \nu }^0=-\sigma _1 \overline{G^0}_{  \nu
\mu} $. Then we have the key structural result of this paper.

\begin{lemma}
\label{lemma:FGR8} We have rhs\eqref{eq:FGR8}$\ge 0.$
\end{lemma}
\proof First of all, it is not restrictive to assume $ G ^0_{  \alpha 0} = P_c(\omega _0) G ^0_{  \alpha 0}  $. We have  $ G ^0_{  \alpha 0}  \in \mathcal{S}(\R ^3, \C ^2)$ for all $\alpha$.  For
   $W(\omega ) =\lim_{t\to\infty}e^{-it \mathcal{H}_\omega }e^{it\sigma_3
(-\Delta+\omega  )}$, there exist  $ F_{  \alpha }  \in W^{k,p}(\R ^3, \C ^2)$
 for all $k\in \R$ and $p\ge 1$ with
 $ G ^0_{  \alpha 0} =W(\omega _0)F_\alpha $, \cite{Cu1}.
By standard theory,   $R_{ \mathcal{H}}^+ (r )G ^0_{  \alpha 0}\in L^ {2,-s}(\R ^3, \C ^2)$ for any $s>1/2 $ and $r>\omega _0$. Let $\mathbf{G}=\sum _{
\lambda ^0\cdot \alpha =r}\zeta ^{ \alpha } G^0 _{ \alpha 0}$ and  $\mathbf{F}=\sum _{
\lambda ^0\cdot \alpha =r}\zeta ^{ \alpha } F _{ \alpha }$. Let $^t{\mathbf{F}}=(\mathbf{F}_1,\mathbf{F}_2)  $. Then, see Lemma 4.1  \cite{Cu2},

\begin{equation} \label{eq:FGR81} \begin{aligned} &
\Im \left    \langle R_{ \mathcal{H}}^+ (r )  \mathbf{G},
\sigma _3\overline{  \mathbf{G} }\right \rangle =\lim _{\varepsilon \searrow 0} \Im  \left    \langle R_{ \mathcal{H}}  (r +\im \varepsilon ) \mathbf{G},
\sigma _3\overline{ \mathbf{G} }\right \rangle \\& =\lim _{\varepsilon \searrow 0} \Im \left    \langle R_{  \sigma _3(-\Delta +\omega _0)}  (r +\im \varepsilon )  \mathbf{F} ,
\sigma _3\overline{  \mathbf{F}}\right \rangle \\& =\lim _{\varepsilon \searrow 0} \Im \left    \langle R_{    -\Delta  }  (r-\omega _0  +\im \varepsilon )  \mathbf{F}_1 ,
 \overline{  \mathbf{F}_1}\right \rangle \\& =\lim _{\varepsilon \searrow 0} \int _{\R ^3} \frac{\varepsilon}{(\xi ^2- (r-\omega _0))^2+\varepsilon ^2  }  |\widehat{\mathbf{F}}_1 (\xi )|^2 d\xi \ge 0.
\end{aligned}
\end{equation}
\qed

\bigskip
Now we will assume the following hypothesis.

\begin{itemize}
\item[(H11)]  We assume
that for some fixed constants for any vector $\zeta \in
\mathbb{C}^n$ we have:
\begin{equation} \label{eq:FGR} \begin{aligned} &  \sum _{
\substack{ \lambda ^0 \cdot  \alpha =\lambda ^0\cdot   \nu   >
\omega _0
   \\ \lambda
\cdot  \alpha -\lambda _k   < \omega _0   \, \forall \, k \, \text{
s.t. } \alpha _k\neq 0\\ \lambda \cdot  \nu -\lambda _k  <\omega _0
\, \forall \, k \, \text{ s.t. } \nu _k\neq 0}} \lambda ^0\cdot \nu
\Im \left ( \zeta ^{ \alpha } \overline{\zeta }^ { \nu  } \langle
R_{ \alpha 0}^+  G_{ \alpha 0}^0,  \sigma _1\sigma _3  G ^0 _{0\nu
}\rangle \right )
   \\&
 \approx \sum _{ \substack{ \lambda ^0\cdot  \alpha
> \omega _0
\\
   \lambda ^0
\cdot  \alpha -\lambda _k  ^0  < \omega _0 \, \forall \, k \, \text{
s.t. } \alpha _k\neq 0}}  | \zeta ^\alpha  | ^2 .
\end{aligned}
\end{equation}

\end{itemize}

 By (H11) we have

\begin{equation} \label{eq:FGR10} \begin{aligned} &
2\sum _j \lambda _j^0\Im \left ( \mathcal{D}_j   \overline{\zeta} _j
\right )\gtrsim \partial _t \sum _j \lambda _j^0| \zeta _j|^2  +  \\&
      \sum _{ \substack{ \lambda ^0\cdot  \alpha
> \omega _0
\\
   \lambda ^0
\cdot  \alpha -\lambda _k  ^0  < \omega _0 \, \forall \, k \, \text{
s.t. } \alpha _k\neq 0}}  | \zeta ^\alpha  | ^2
\end{aligned}
\end{equation}
Then, for $t\in [0,T]$ and assuming Lemma \ref{lemma:FGR1} we have

\begin{eqnarray}& \sum _j \lambda _j ^0 | \zeta
_j(t)|^2 +\sum _{ \substack{ \lambda ^0\cdot  \alpha
> \omega _0
\\
   \lambda ^0
\cdot  \alpha -\lambda _k ^0  < \omega   _0\, \forall \, k \, \text{
s.t. } \alpha _k\neq 0}}  \| \zeta ^\alpha \| _{L^2(0,t)}^2\lesssim
\epsilon ^2+ C_2\epsilon ^2.\nonumber
\end{eqnarray}
By \eqref{equation:FGR3} this implies $\| z ^\alpha \|
_{L^2(0,t)}^2\lesssim \epsilon ^2+ C_2\epsilon ^2$ for all the above
multi indexes. So, from  $\| z ^\alpha \| _{L^2(0,t)}^2\lesssim
  C_2^2\epsilon ^2$ we conclude $\| z ^\alpha \|
_{L^2(0,t)}^2\lesssim  C_2\epsilon ^2$. This means that we can take
$C_2\approx 1$. This yields Theorem \ref{proposition:mainbounds}.

\begin{remark}
\label{rem:genericity} Notice that by $r>\omega _0$, \eqref{eq:FGR}
  appears generic. We do not try to prove this point. It
should not be hard, see for example the genericity result
Proposition 2.2 \cite{bambusicuccagna}.
\end{remark}

\begin{remark}
\label{rem:integrability} In general we expect Hypothesis (H11), or
higher order versions, to hold. Specifically, if at some step of the
normal form argument (H11) fails because some of the inequalities as
in Lemma \ref{lemma:FGR8} is an equality,
  one can continue the normal form procedure and obtain some steps
later a new version  of (H11). This will yield an analogue of
Theorem \ref{proposition:mainbounds}, with \ref{L^2discrete}
replaced by a similar but weaker inequality. We could have stated
(H11) and proved Theorem \ref{proposition:mainbounds} in this more
general form, but this would have complicated further the
presentation.
\end{remark}
\begin{remark}
\label{rem:excited states} If instead of ground states we consider standing
waves with nodes, and if $\dim N_g(\mathcal{H}_{\omega})=2 $ with \eqref{eq:Kernel}, if we assume (H1)--(H11) with \eqref{eq:1.2} in (H5)  replaced by $\frac d {d\omega } \| \phi _ {\omega }\|^2_{L^2(\R^3)}\neq 0$,
if we assume $\sigma (\mathcal{H}_{\omega})\subset \R$, then
by \cite{Cu3} it is possible to prove that the hamiltonian $K$ in Lemma \ref{lem:K} has quadratic part
$$K_2=\sum _{j=1}^{m} \gamma _j \lambda _j (\omega  ) |z_j|^2+
\frac{1}{2} \langle \sigma _3 \mathcal{H}_{\omega  } f, \sigma _1
f\rangle $$
with $\gamma _j$ equal to either 1 or $-1$ and with at least one
$\gamma _j=-1$ (in other words the energy has a saddle
at $\phi_{\omega  }$ in the surface formed by the $  u$ with $\| u \| _{L^2}
=\| \phi_{\omega  } \| _{L^2}$). Then a simple elaboration of the proof of the present paper, along the lines of sections 3 or 4 in  \cite{Cu3},
can be used to strengthen Theorem 3.2  \cite{Cu3}
showing that $\phi_{\omega  } $ is orbitally unstable. Furthermore,  following the argument in \cite{Cu3},
it can be shown that if a solution $u(t)$  remains close to ground states as $t\nearrow +\infty$ (resp. $t\searrow -\infty$), it actually scatters to ground states, that is   \eqref{scattering} and \eqref{Strichartz} with the plus (resp. minus) sign.
\end{remark}

DISMI University of Modena and Reggio Emilia, Via Amendola 2, Pad.
Morselli, Reggio Emilia 42122, Italy.

{\it E-mail Address}: {\tt cuccagna.scipio@unimore.it}

\end{document}